\documentclass[11pt,a4paper,reqno]{amsart}
\usepackage[T1]{fontenc}
\usepackage{mathrsfs}
\usepackage{amsfonts,amssymb,amsmath,amsgen,amsthm}
\usepackage[pdftex]{graphicx,color,epsfig}%les deux lignes a rajouter pour le graphique
\usepackage{amsbsy}
\usepackage{mathrsfs}
\usepackage{bbm}
\usepackage{titletoc}
\usepackage{enumerate}
\usepackage{hyperref}
\usepackage{subfigure}
\usepackage[subfigure]{ccaption}
\theoremstyle{plain}
\newtheorem{theorem}{Theorem}[section]

\newtheorem{lemma}[theorem]{Lemma}

\newtheorem{proposition}[theorem]{Proposition}

\theoremstyle{remark}
\newtheorem{remark}[theorem]{Remark} 
\newtheorem*{notation}{Notation}
\newtheorem{example}[theorem]{Example}

\catcode`@=11
\def\section{\@startsection{section}{1}%
  \z@{1.5\linespacing\@plus\linespacing}{.5\linespacing}%
  {\normalfont\bfseries\large\centering}}
\catcode`@=12
\def\C{{\mathbb C}}% complex numbers
\def\R{{\mathbb R}}% real numbers
\def\N{{\mathbb N}}% nonnegative integers
\def\Z{{\mathbb Z}}% integers
\def\T{{\mathbb T}}% integers
\def\O{\mathcal O}

\newcommand{\dt}{\ensuremath{\Delta t}}
\newcommand{\dx}{\ensuremath{\Delta x}}

\def\({\left(}
\def\){\right)}
\def\<{\left\langle}
\def\>{\right\rangle}
\def\le{\leqslant}
\def\ge{\geqslant}

\def\d{{\partial}}
\def\eps{\varepsilon}

\def\eik{\phi_{\rm eik}}

\DeclareMathOperator{\RE}{Re}
\DeclareMathOperator{\IM}{Im}
\DeclareMathOperator{\DIV}{div}

\numberwithin{equation}{section}

\newcommand{\be}{\begin{equation}}
\newcommand{\ee}{\end{equation}}
\newcommand{\bea}{\begin{eqnarray}}
\newcommand{\eea}{\end{eqnarray}}
\newcommand{\bee}{\begin{eqnarray*}}
\newcommand{\eee}{\end{eqnarray*}}

\def\ni{\noindent}
\def\bs{\bigskip}
\def\ms{\medskip}

\def\eps{\varepsilon}

\def\pref#1{{\rm \ref{#1}}}

\def\pa{\partial}
\def\na{\nabla}

\begin{document}

\title[AP scheme for NLS]{An asymptotic preserving scheme based on a new formulation for NLS  in the semiclassical limit}  
\author[C. Besse]{Christophe Besse}
\email{Christophe.Besse@math.univ-lille1.fr}
\address{Laboratoire Paul Painlev\'e, Universit\'e de Lille 1}
\author[R. Carles]{R\'emi Carles}
\email{Remi.Carles@math.cnrs.fr}
\address{CNRS \& Univ. Montpellier 2, Math\'ematiques}
\author[F. M\'ehats]{Florian M\'ehats}
\email{florian.mehats@univ-rennes1.fr}
\address{IRMAR, Universit\'e de Rennes 1 and INRIA, IPSO Project}

\begin{abstract}
  We consider the semiclassical limit for the nonlinear
Schr\"odinger equation. We introduce a phase/amplitude representation 
given by a system similar to the hydrodynamical formulation, whose 
novelty consists in including some asymptotically vanishing viscosity. 
We prove that the system is always locally well-posed in a class of 
Sobolev spaces, and globally well-posed for a fixed positive Planck 
constant in the one-dimensional case. We propose a second order
numerical scheme which is asymptotic preserving. Before singularities appear in 
the limiting Euler equation, we recover the quadratic physical
observables as well as the wave function with mesh size and time step 
independent of the Planck constant. This approach is 
also well suited to the linear Schr\"odinger equation. 
\end{abstract}
\maketitle

\section{Introduction and main results}
\label{sec:intro}

\subsection{Motivation}
\label{sec:motiv}
We consider the cubic, defocusing, nonlinear Schr\"odinger equation (NLS)
\begin{equation}
  \label{eq:nls}
  i\eps\d_t u^\eps +\frac{\eps^2}{2}\Delta u^\eps = |u^\eps|^2
  u^\eps,\quad (t,x)\in \R_+\times \R^d,
\end{equation}
with WKB type initial data
\begin{equation}
  \label{eq:ci}
  u^\eps(0,x) = a_0(x)e^{i\phi_0(x)/\eps},
\end{equation}
the functions $a_0$ and $\phi_0$ being real-valued.
The aim of this article is to construct an {\em Asymptotic Preserving}
(AP) numerical scheme for this equation in the semiclassical limit
$\eps\to 0$. We seek a scheme which provides an approximation of the
solution $u^\eps$ with an accuracy, for fixed numerical parameters
$\Delta t$, $\Delta x$, that is not degraded as the scaled Planck
constant $\eps$ goes to zero. In other terms, such a scheme is
consistent with \eqref{eq:nls} for all fixed $\eps>0$, and when $\eps
\to 0$ converges to a consistent approximation of the limit equation,
which is the compressible, isentropic Euler equation
\eqref{eq:euler}. The difficulty here is that, when $\eps$ is small,
the solution $u^\eps$ becomes highly oscillatory with respect to the
time and space variables, and converges to its limit only in a weak
sense. In order to follow these oscillations without reducing $\Delta
t$ and $\Delta x$ at the size of $\eps$, which may be computationally
demanding, our construction relies on a fluid reformulation of
\eqref{eq:nls}, well adapted to the semiclassical limit.

In kinetic theory, plasma physics and radiative transfer, the interest
in AP schemes has considerably grown in the recent years since the
pioneer works  on the subject \cite{coronperthame,GPT,jin,klar}. Among
a long list of works on this topic, we can mention
\cite{JPT,gossetoscani,buetdespres,lemoumieussens,bennoune,crispeldegondvignal,filbetjin,CGLV,degond2,lafittesamaey,degond3}. For
the Schr\"odinger equation, fewer works exist. In the stationary case,
one can cite \cite{nba} which proposes a WKB-type transformation in
order to filter out the oscillations in space. In the time-dependent
case, the usual numerical methods for solving the nonlinear
Schr\"odinger equation -- finite-difference schemes
\cite{wu,DFP,akrivis,akrivis2}, time-splitting methods
\cite{PM,weideman,BBD}, relaxation scheme \cite{besse} or even
symplectic methods \cite{sanzserna} -- are designed in order to
guarantee the convergence of the discrete wave function $u^\eps$ when
$\eps>0$ is fixed. As it is analyzed for instance in \cite{MPP,MPPS}
for the case of finite-differences or in \cite{klein} for Runge-Kutta
schemes, these methods suffer from the oscillations in the
semiclassical regime if the time and space steps are not shrunk. Note
however that, among these methods, the status of splitting methods is
particular. Indeed, in the linear case, it is shown in
\cite{BaoJinMark} (see also the recent review article \cite{JMS}) by
Wigner techniques that, if the space step $\Delta x$ has to follow the
parameter $\eps$ when it becomes small, the time step can be chosen
independently of $\eps$ (before the appearance of caustics) if we are
only interested in getting a good approximation of the quadratic
observables associated to $u^\eps$. Still, none of these methods is
Asymptotic Preserving in the sense defined above. 

The quantum fluid reformulation of the Schr\"odinger equation provides
a natural framework for the construction of AP numerical methods. The
Madelung transform \cite{madelung} is a polar decomposition of the
solution of \eqref{eq:nls}, written as 
$$
u^\eps(t,x)=a^\eps(t,x)e^{i\phi^\eps(t,x)/\eps},
$$
where the amplitude $a^\eps$ and the phase $\phi^\eps$ are
real-valued. Inserting this \emph{ansatz} into \eqref{eq:nls} yields
the system of quantum hydrodynamics (QHD), or Madelung equations,
satisfied by $\rho^\eps=(a^\eps)^2$ and $v^\eps=\na \phi^\eps$: 
\begin{equation}\label{eq:qhd}
  \left\{
    \begin{aligned}
    &\d_t \rho^\eps + \DIV \(\rho^\eps v^\eps \)=0,\qquad  &&\rho^\eps_{\mid
      t=0}=(a_0)^2, \\
      & \d_t v +v^\eps\cdot \nabla v^\eps+
      \nabla \rho^\eps  =\frac{\eps^2}{2}\na\left(\frac{\Delta \sqrt{\rho^\eps}}{\sqrt{\rho^\eps}}\right),\qquad &&
      v^\eps_{\mid t=0}=\nabla\phi_0.
      \end{aligned}
\right.
\end{equation}
This system takes the form of the compressible Euler equation (with
the pressure law $p(\rho)=\rho^2/2$), with an additional term
testifying for the quantum character of the equation, the so-called
Bohm potential $\frac{\eps^2}{2}\na\left(\frac{\Delta
    \sqrt{\rho^\eps}}{\sqrt{\rho^\eps}}\right)$. Two comments are in
order. First, the term in $\eps$ does not appear any more as a
singular perturbation and, in the limit $\eps\to 0$, the quantum
pressure disappears, leading formally to the compressible, isentropic
Euler equation 
\begin{equation}\label{eq:euler}
  \left\{
    \begin{aligned}
    &\d_t \rho + \DIV \(\rho v \)=0,\qquad  &&\rho_{\mid
      t=0}=(a_0)^2, \\
      & \d_t v +v\cdot \nabla v+
      \nabla \rho  = 0,\qquad &&
      v_{\mid t=0}=\nabla\phi_0,
      \end{aligned}
\right.
\end{equation}
which is the semiclassical limit of our problem (actually valid before
the appearance of shocks). Second, the fluid-like form of this system
enables to consider many numerical methods originated from
computational fluid mechanics. 

For the linear Schr\"odinger equation, the Madelung formulation is at
the heart of the method of quantum trajectories (or Bohmian dynamics),
see \cite{wyatt}, which are based on a Lagrangian resolution of this
system. In Eulerian coordinates, the work \cite{DGM} also exploits
this formulation to construct Asymptotic-Preserving schemes for
\eqref{eq:qhd}. Unfortunately, the main drawback of the QHD system is
the form of the quantum potential, which becomes singular as the
density $\rho^\eps$ vanishes (see \cite{CaDaSa12} for a recent
survey). Hence, these methods do not provide AP 
schemes in the presence of vacuum. Another point of view was recently
developed in \cite{remibijan} for the nonlinear equation
\eqref{eq:nls}, taking advantage of another fluid reformulation of the
equation --\,due to Grenier \cite{Grenier98}\,-- which tolerates the
presence of vacuum. In the next Subsection~\ref{sec:grenier}, we
present in detail this formulation relying on a complex amplitude
version of Madelung transform. However, as we will see, the drawback
of this reformulation  is that this system may develop singularities
even for fixed $\eps>0$, whereas the solution of NLS remains
smooth. Then, in the last Subsection~\ref{sec:model} of this
introduction, we present a modification in order to remedy this
problem and preserve the regularity of the solution. The new QHD model
that we propose contains a viscosity term in the equation for the
velocity, but is still equivalent to the original equation
\eqref{eq:nls}. 

\subsection{Backgrounds on NLS and its semiclassical limit}
\label{sec:grenier}

In this section, we recall some known results on the Cauchy problem for NLS \eqref{eq:nls} and on its semiclassical limit as $\eps\to0$. For details on this subject, we refer to the textbooks \cite{CazCourant} and \cite{bookRemi}.

The following result is standard in the case $\eps=1$, and can easily be deduced when $\eps\in (0,1]$, by scaling arguments.
\begin{proposition}\label{prop:nls}  Let $d\le 3$, and $u_0^\eps\in H^1(\R^d)$.
 \begin{enumerate}[(i)]
  \item There exists a unique solution $u^\eps\in
    C(\R;H^1(\R^d))\cap L^{8/d}_{\rm loc}(\R;W^{1,4}(\R^d))$ to
    \eqref{eq:nls}, such that $u_{\mid t=0}^\eps=u_0^\eps$. Moreover,
    the following conservations hold:
    \begin{align*}
      \text{Mass:}\quad & \frac{d}{dt}\|u^\eps(t)\|_{L^2}^2 =0.\\
\text{Momentum:}\quad& \frac{d}{dt}\IM\int_{\R^d} \overline{u^\eps}(t,x)\nabla
u^\eps(t,x)dx =0.\\
\text{Energy:}\quad & \frac{d}{dt}\(\|\eps\nabla
u^\eps(t)\|_{L^2}^2+\|u^\eps(t)\|_{L^4}^4\)=0. 
    \end{align*}
\item If in addition $u_0^\eps\in H^k(\R^d)$ for $k\in \N$, $k\ge 2$,
  then $u^\eps\in C(\R;H^k(\R^d))$. 
  \end{enumerate}
\end{proposition}
Note that, if $u^\eps = a^\eps e^{i\phi^\eps/\eps}$ with $\phi^\eps\in \R$, the above three conservation laws become
\begin{align}
      & \frac{d}{dt}\|a^\eps(t)\|_{L^2}^2 =0.\label{eq:masse}\\
& \frac{d}{dt}\(\IM\int_{\R^d} \overline{a^\eps}(t,x)\nabla
a^\eps(t,x)dx+\frac{1}{\eps}\int_{\R^d}|a^\eps(t,x)|^2 \nabla
\phi^\eps(t,x)dx\)=0.\notag\\
& \frac{d}{dt}\(\|\eps\nabla
a^\eps(t)+ia^\eps(t)\nabla\phi^\eps(t)\|_{L^2}^2+\|a^\eps(t)\|_{L^4}^4\)=0.
\label{eq:energy}   
\end{align}
%\begin{remark}
%  {\tt pour NLS focalisant, les estimations d'energie sont moins
%    interessantes. Le probleme est plus grave : l'equation d'Euler
%    \eqref{eq:eulersym} devient elliptique, et ca change beaucoup de
%    choses (le cadre fonctionnel doit etre adapte : il faut travailler
%    dans l'analytique)}
%\end{remark}

Let us now describe the limit $\eps\to 0$ for \eqref{eq:nls}. The
toolbox for studying semiclassical Schr\"odinger equations contains a
variety of methods, depending on the quantities we are interested
in. If one is  only interested in a description of the dynamics of
quadratic observables, such as the mass, current or energy densities,
the Wigner transform is well adapted and has been applied to linear
Schr\"odinger equations and Schr\"odinger-Poisson systems
\cite{GMMP,lionspaul,ZhangSIMA}; on the other hand, it is not adapted to the
study of NLS \cite{CFMS09}. Still for a description of quadratic observables, the modulated energy method \cite{brenier,zhang,LinZhang,ACMA} enables to treat the nonlinear equation \eqref{eq:nls}. Yet, for a pointwise description of the wavefunction $u^\eps$, WKB techniques are more convenient. We refer to \cite{bookRemi} for a presentation of WKB analysis of Schr\"odinger equations. 

In the case of \eqref{eq:nls}, also referred to as supercritical
geometric optics in this context, the justification of the WKB
expansion was done by Grenier \cite{Grenier98}. Let us briefly present
his idea. As we said above, the main inconvenience of the QHD system
\eqref{eq:qhd} is the singularity of the Bohm potential at the zeroes
of the function $\rho^\eps$ (\emph{i.e.} at the zeroes of the wavefunction
$u^\eps$). Looking again for solutions $u^\eps$ under the form  
\be
\label{madelung}
u^\eps(t,x)=a^\eps(t,x)e^{i\phi^\eps(t,x)/\eps},
\ee
but allowing now $a^\eps$ to take complex values, Grenier proposed to define $a^\eps$ and $\phi^\eps$ as the solutions of the system
\begin{equation}\label{eq:systemmanuel}
  \left\{
    \begin{aligned}
     &\d_t a^\eps +\nabla \phi^\eps \cdot \nabla a^\eps
      +\frac{1}{2}a^\eps \Delta \phi^\eps =i\frac{\eps}{2}\Delta
      a^\eps,\qquad  && a^\eps_{\mid
      t=0}=a_0,\\
      & \d_t \phi^\eps +\frac{1}{2}|\nabla \phi^\eps|^2+
      |a^\eps|^2  = 0,\qquad &&
      \phi^\eps_{\mid t=0}=\phi_0.  \end{aligned}
\right.
\end{equation} 
This system, which is formally equivalent to \eqref{eq:nls}, can also
be expressed in terms of the unknowns $a^\eps$ and
$v^\eps=\na\phi^\eps$. Differentiating with respect to $x$ the second equation in
\eqref{eq:systemmanuel}, we find, using the fact 
that $v^\eps$ is irrotational:
\begin{equation}\label{eq:systemmanuel2}
  \left\{
    \begin{aligned}
    & \d_t a^\eps +v^\eps \cdot \nabla a^\eps
      +\frac{1}{2}a^\eps \DIV v ^\eps =i\frac{\eps}{2}\Delta
      a^\eps,\qquad  && a^\eps_{\mid
      t=0}=a_0,\\
      & \d_t v^\eps +v^\eps\cdot \nabla v^\eps+
      \nabla |a^\eps|^2  = 0,\qquad &&
      v^\eps_{\mid t=0}=\nabla\phi_0.
    \end{aligned}
\right.
\end{equation}
The important remark made in \cite{Grenier98} is to notice that 
the above 
system is  hyperbolic symmetric, perturbed by a
skew-symmetric term. In the limit case $\eps=0$, and if $a_0$ is
real-valued and nonnegative, 
\eqref{eq:systemmanuel2} is the Euler equation \eqref{eq:euler}
written in symmetric form (see e.g. \cite{MUK86,JYC90}):
\begin{equation}\label{eq:eulersym}
  \left\{
    \begin{aligned}
    & \d_t a +v \cdot \nabla a
      +\frac{1}{2}a \DIV v =0,\qquad  && a_{\mid
      t=0}=a_0,\\
      & \d_t v +v\cdot \nabla v+
      \nabla |a|^2  = 0,\qquad && 
      v_{\mid t=0}=\nabla\phi_0.
    \end{aligned}
\right.
\end{equation}
This remark enables to obtain, in small time (but independent of
$\eps\in [0,1]$), uniform estimates of the solution $(a^\eps,v^\eps)$
of \eqref{eq:systemmanuel2} in Sobolev norms and prove the convergence
of these quantities to the solution $(a,v)$ of \eqref{eq:eulersym}, as
long as this solution exists and is smooth. From the numerical point
of view, the regular structure of \eqref{eq:systemmanuel2}, in
particular when $a^\eps$ vanishes or has very small values, is an
advantage compared to \eqref{eq:qhd}, and is at the basis of the
schemes constructed in \cite{remibijan}. 

\subsection{New model and main results}
\label{sec:model}

A drawback of \eqref{eq:systemmanuel2} is that large time existence results (for
fixed $\eps>0$) do not seem available. In particular, it is not known
whether \eqref{eq:systemmanuel2} still has a smooth solution past the
time where the solution to \eqref{eq:eulersym}
has ceased to be smooth (such a time necessarily exists if $\phi_0$ and
$a_0$ are compactly supported, regardless of their size, from
\cite{MUK86} --- see also \cite{Xin98}). A practical consequence of this fact is that a numerical method based on \eqref{eq:systemmanuel2}, such as the one proposed in \cite{remibijan}, may not approximate correctly the original equation \eqref{eq:nls} on arbitrary time intervals, for fixed $\eps>0$.

To overcome this issue, we use the degree of freedom given by
\eqref{madelung} in a manner which is slightly different from the
approach introduced by E.~Grenier. From \eqref{madelung}, we have
\begin{align*}
  i\eps\d_t u^\eps +\frac{\eps^2}{2}\Delta u^\eps - |u^\eps|^2
  u^\eps &
= \ -\(\d_t \phi^\eps + \frac{1}{2}|\nabla \phi^\eps|^2+
      |a^\eps|^2\)a^\eps e^{i\phi^\eps/\eps} \\ 
&+ i\eps\( \d_t a^\eps +\nabla \phi^\eps \cdot \nabla a^\eps
      +\frac{1}{2}a^\eps \Delta \phi^\eps -i\frac{\eps}{2}\Delta
      a^\eps\)e^{i\phi^\eps/\eps} .
\end{align*}
Introduce the viscous term $\eps^2 a^\eps e^{i\phi^\eps/\eps}\Delta
\phi^\eps$, and reorder the terms so as to consider the new system
\begin{equation}\label{eq:bcm0}
  \left\{
    \begin{aligned}
   & \d_t a^\eps +\nabla \phi^\eps \cdot \nabla a^\eps
      +\frac{1}{2}a^\eps \Delta \phi^\eps =i\frac{\eps}{2}\Delta
      a^\eps-i\eps a^\eps\Delta \phi^\eps,\qquad  && a^\eps_{\mid
      t=0}=a_0,    \\
   & \d_t \phi^\eps +\frac{1}{2}|\nabla \phi^\eps|^2+
      |a^\eps|^2  = \eps^2\Delta \phi^\eps,\qquad &&
      \phi^\eps_{\mid t=0}=\phi_0.
    \end{aligned}
\right.
\end{equation} 
Proceeding like above, we get the following new system:
\begin{equation}\label{eq:bcm}
  \left\{
    \begin{aligned}
      & \d_t a^\eps +v^\eps \cdot \nabla a^\eps
      +\frac{1}{2}a^\eps \DIV v ^\eps =i\frac{\eps}{2}\Delta
      a^\eps -i\eps a^\eps \DIV v ^\eps,\qquad  &&a^\eps_{\mid
      t=0}=a_0^\eps,\\
      & \d_t v^\eps +v^\eps\cdot \nabla v^\eps+
      2 \RE \(\overline a^\eps \nabla a^\eps\)  = \eps^2\Delta v^\eps,\qquad &&
      v^\eps_{\mid t=0}=\nabla\phi_0.
    \end{aligned}
\right.
\end{equation}
We recover physical observables, respectively particle, current and energy densities, thanks to the following formula
\begin{equation}\label{eq:observables}
    \begin{aligned}
      & \rho^\eps=|a^\eps|^2\\
& \mathbf{j}^\eps=\eps \IM(\overline{a^\eps} \nabla a^\eps)+\rho^\eps v^\eps, \\
& e^\eps=\left | \eps \nabla a^\eps+ i a^\eps v^\eps \right |^2 + |a^\eps|^4.
    \end{aligned}
\end{equation}
\begin{remark}[Linear case]
  In the case of a linear Schr\"odinger equation
  \begin{equation*}
    i\eps\d_t u^\eps +\frac{\eps^2}{2}\Delta u^\eps = V u^\eps,
  \end{equation*}
we will see that the introduction of this artificial diffusion is even
more striking than in the nonlinear framework; see
Subsection~\ref{sec:linear}. 
\end{remark}
\begin{remark}[Irreversibility]
Our choice for introducing the viscous term in \eqref{eq:bcm0} makes
the system irreversible (in time), while \eqref{eq:nls} is reversible. To solve the Schr\"odinger equation backward in time, it suffices to change the signs to
\begin{equation*}
  \left\{
    \begin{aligned}
   & \d_t a^\eps +\nabla \phi^\eps \cdot \nabla a^\eps
      +\frac{1}{2}a^\eps \Delta \phi^\eps =i\frac{\eps}{2}\Delta
      a^\eps+i\eps a^\eps\Delta \phi^\eps,\qquad  && a^\eps_{\mid
      t=0}=a_0,    \\
   & \d_t \phi^\eps +\frac{1}{2}|\nabla \phi^\eps|^2+
      |a^\eps|^2  = -\eps^2\Delta \phi^\eps,\qquad &&
      \phi^\eps_{\mid t=0}=\phi_0.
    \end{aligned}
\right.
\end{equation*} 
\end{remark}
Let us now present our main theoretical results on this system. Our first result concerns the local in time well-posedness of the Cauchy problem in any dimension, for fixed values of $\eps\ge 0$.
\begin{theorem}\label{theo:local}
  Let $(a_0,\phi_0)\in H^{s}(\R^d)\times H^{s+1}(\R^d)$, with $s>d/2+1$. Then the following holds.\\
  (i) The system \eqref{eq:bcm} admits a unique maximal solution
  $(a^\eps,v^\eps)\in C([0,T_{\rm max}^\eps);H^s\times H^s)$.\\
  (ii) There exists $T>0$ independent of $\eps\in [0,1]$ such that $T_{\rm max}^\eps\ge T$. Moreover, the $L^\infty([0,T];H^s\times H^s)$ norm of $(a^\eps,v^\eps)$ is bounded uniformly in $\eps\in [0,1]$.\\
(iii)  Defining $\phi^\eps \in C([0,T^\eps_{\rm max});H^{s+1})$ by
 \begin{equation}\label{eq:phi}
    \phi^\eps(t,x) = \phi_0(x) -\int_0^t \(\frac{1}{2}|v^\eps(s,x)|^2
    +|a^\eps(s,x)|^2-\eps^2\DIV v^\eps(s,x)\)ds,
  \end{equation}
  then the function
  \begin{equation*}
    u^\eps = a^\eps e^{i\phi^\eps/\eps}\in C([0,T_{\rm max}^\eps);H^{s})
  \end{equation*}
is the unique solution to \eqref{eq:nls} satisfying \eqref{eq:ci}. 
\end{theorem}
\begin{remark}\label{rem:solmax}
In fact, the proof of this theorem provides a criteria of global existence for the solution to \eqref{eq:bcm}. Indeed, we shall prove that the lifespan $T_{\rm max}^\eps$ is independent of $s>d/2+1$ and that 
  \begin{equation}
  \label{lifespan}
 T_{\rm max}^\eps<\infty\Longrightarrow   \int_0^{T_{\rm max}^\eps}
   \( \|(a^\eps,v^\eps)(t)\|_{W^{1,\infty}}+\|a^\eps(t)\|_{L^\infty}^2\)dt=\infty. 
  \end{equation}
\end{remark}
Our second result states the convergence of the solution to \eqref{eq:bcm} towards the solution of the Euler equation, as $\eps\to 0$.
\begin{notation}
 For two positive numbers $\alpha^\eps$ and $\beta^\eps$,
the notation $ \alpha^\eps\lesssim \beta^\eps$
means that there exists $C>0$ \emph{independent of} $\eps$ such that
for all $\eps\in (0,1]$, $\alpha^\eps\le C\beta^\eps$.  
\end{notation}
\begin{theorem}\label{theo:cveuler}
  Let $s>d/2+3$, and $(a_0,\phi_0)\in H^{s}\times H^{s+1}$. Let $T>0$
  such that the Euler equation \eqref{eq:euler} has a unique solution
  $(\rho,v)\in C([0,T];H^s\times H^s)$. Then \eqref{eq:eulersym} has a
  unique solution $(a,v)\in C([0,T];H^s\times H^s)$. Moreover, for
  $\eps>0$ sufficiently small, we have $T^\eps_{\rm max}\ge T$ and
  \begin{equation*}
    \|(a^\eps,v^\eps)-(a,v)\|_{L^\infty([0,T];H^{s-2}\times H^{s-2})}\lesssim
    \eps. 
  \end{equation*}
\end{theorem}
Theorems~\ref{theo:local} and \ref{theo:cveuler} are not different
from the corresponding results obtained by Grenier in \cite{Grenier98}
on his system \eqref{eq:systemmanuel2}. We now state a global
existence result for fixed $\eps>0$ in dimension one, that is not
available on \eqref{eq:systemmanuel2}. This result takes advantage of
the viscous term in the second equation of \eqref{eq:bcm}. 
\begin{theorem}\label{theo:global}
  Suppose $d=1$ and let $(a_0,\phi_0)\in H^{s}\times H^{s+1}$, with $s\ge 2$. Then for $\eps>0$ fixed, the solution to
  \eqref{eq:bcm} is global in time, in the sense that $T_{\rm
    max}^\eps=\infty$ in  Theorem~\pref{theo:local}. 
\end{theorem}

\bs
This article is organized as follows. Section \ref{section:theo} is
devoted to the proofs of our three theorems: the local existence
result in Subsection \ref{sec:local}, the semiclassical limit in
Subsection \ref{sec:cvEuler} and the global existence result in
Subsection \ref{sec:global}. Section \ref{sec:numerique} is of
different nature and concerns numerics. We describe a second-order
Asymptotic Preserving numerical scheme for NLS, based on the
reformulation \eqref{eq:bcm}, and we give the results of numerical
experiments, in dimensions 1 and 2. We show that our scheme is indeed
AP before the formation of singularities in the Euler equation, which
means that the time and space steps $\Delta t$ and $\Delta x$ can be
taken independently of $\eps$. After the formation of singularities,
we show that, in order to get a good approximation of the solution, we
need to take $\Delta t$ and $\Delta x$ of order $\mathcal O(\eps)$. At
the end of the article, in Section \ref{sec:extension}, we sketch two
extensions of our results: we discuss on the case of the linear
Schr\"odinger equation, then we discuss on other nonlinearities.

%For two positive numbers $a^\eps$ and $b^\eps$, the notation $ a^\eps\lesssim b^\eps$ means that there exists $C>0$ \emph{independent of} $\eps$ such that for all $\eps\in (0,1]$, $a^\eps\le Cb^\eps$.  

%\input{nls}

%\input{local}

\section{Proofs of the main theorems}
\label{section:theo}

This section is devoted to the proofs of our main results stated in Theorems \ref{theo:local}, \ref{theo:cveuler} and \ref{theo:global}.

\subsection{Local existence}
\label{sec:local}
In this paragraph, we show the local existence of a smooth maximal solution to \eqref{eq:bcm}.
\begin{proof}[Proof of Theorem \pref{theo:local}]
{\em Items (i) and (ii): local existence and uniform estimates.} Following \cite{Grenier98}, introduce the vector-valued unknown function
\begin{equation*}
  U^\eps = 
    \begin{pmatrix}
       \RE a^\eps \\
       \IM a^\eps \\
       v^\eps_1 \\
      \vdots \\
       v^\eps_d
    \end{pmatrix}
= 
    \begin{pmatrix}
       a_1^\eps \\
       a_2^\eps \\
       v^\eps_1 \\
      \vdots \\
       v^\eps_d
    \end{pmatrix}
\in \R^{{d+2}}.
\end{equation*}
In terms of this unknown function, \eqref{eq:bcm} reads
\begin{equation}
  \label{eq:systhyp}
  \partial_t U^\varepsilon +\sum_{j=1}^d
  A^j(U^\varepsilon)\partial_j U^\varepsilon 
  = \frac{\varepsilon}{2} L 
  U^\varepsilon + \eps^2 DU^\eps + \eps
  \sum_{j=1}^dB^j(U^\eps)\partial_j U^\eps ,
\end{equation}
where the $(k,\ell)_{1\le k,\ell\le d+2}$ elements of the matrices $A^j(U) \in {\mathcal M}_{d+2}(\R)$ are given by
\begin{equation*}
\left\{\begin{array}{l}
A^j_{k,k}(U)= U_{j+2}\quad \mbox{for } k=1,\ldots,d+2\\[2mm]
A^j_{1,j+2}(U)=U_1/2,\quad A^j_{2,j+2}(U)=U_2/2\\[2mm]
A^j_{j+2,1}(U)=2U_1,\quad A^j_{j+2,2}(U)=2U_2\\[2mm]
A^j_{k,\ell}(U)=0\quad \mbox{otherwise.}
\end{array}\right.
\end{equation*}
The linear operator $L$
is given by
\begin{equation*}
  L = 
  \begin{pmatrix}
   0  &-\Delta &0& \dots & 0   \\
   \Delta  & 0 &0& \dots & 0  \\
   0& 0 &&0_{d\times d}& 
   \end{pmatrix}.
\end{equation*}
A first important remark is that even though $L$ is a differential
operator of order two, it causes no loss of regularity in the energy
estimates, since it is skew-symmetric. Also, \eqref{eq:systhyp} is
hyperbolic symmetric (or symmetrizable). Indeed, let  
\begin{equation}\label{eq:symetriseur}
  S=\begin{pmatrix}
     I_2 & 0\\
     0& \frac{1}{4}I_d
    \end{pmatrix}.
\end{equation}
This matrix is symmetric and positive, and
$SA^j$ is symmetric for all $k$
\begin{equation*}
  SA^j(U)\in {\mathcal S}_{d+2}(\R),\quad \forall U\in
  \R^{d+2}.
\end{equation*}
The diffusive term $D$ is given by
\begin{equation*}
  D =
  \begin{pmatrix}
    0 & 0 & 0&\dots & 0\\
0& 0 &0&\dots & 0\\
0 & 0& &I_{d\times d}\Delta & 
  \end{pmatrix}.
\end{equation*}
Finally, the matrices $B^j$, defined analogously to $A^j$, are given by
\begin{equation*}
\left\{\begin{array}{l}
B^j_{1,j+2}(U)=U_2,\quad B^j_{2,j+2}(U)=-U_1\\[2mm]
B^j_{k,\ell}(U)=0\quad \mbox{otherwise.}
\end{array}\right.
\end{equation*}
When the right-hand side of \eqref{eq:systhyp} is zero, local
existence of a unique solution in $H^s$ with $s>d/2+1$ is a consequence of
standard quasilinear analysis for hyperbolic symmetric systems; see
e.g. \cite{Taylor3}. The important point in energy estimates consists
in computing
\begin{equation*}
  \frac{d}{dt}\<\Lambda^s (SU^\eps),\Lambda^s U^\eps\>,\quad
  \Lambda=(1-\Delta)^{1/2}, 
\end{equation*}
and taking advantage of symmetry to obtain, thanks to tame estimates,
\begin{equation*}
  \frac{d}{dt}\<\Lambda^s (SU^\eps),\Lambda^s U^\eps\>\le
  C\|U^\eps\|_{W^{1,\infty}}\<\Lambda^s (SU^\eps),\Lambda^s U^\eps\>,
\end{equation*}
for some constant $C$ depending only on $A$ and $d$.

In the case of the complete system \eqref{eq:systhyp}, we first recall
that as noticed in \cite{Grenier98}, the term $L$ does not alter the
above energy estimate, since $SL$ is skew-symmetric, so we have
\begin{align*}
  \frac{d}{dt}\<\Lambda^s (SU^\eps),\Lambda^s U^\eps\>& \le
  C\|U^\eps\|_{W^{1,\infty}}\<\Lambda^s (SU^\eps),\Lambda^s U^\eps\> +
  2\eps^2\<S\Lambda^s DU^\eps, \Lambda^s U^\eps\> \\
&\quad+ 2\eps\sum_{j=1}^d
  \<S\Lambda^s \(B^j(U^\eps)\d_j U^\eps\),\Lambda^s U^\eps\>. 
\end{align*}
In terms of $(v^\eps,a^\eps)$, the first new term reads
\begin{equation*}
  2\eps^2\<S\Lambda^s DU^\eps, \Lambda^s U^\eps\>=\frac{\eps^2}{2}\sum_{k=1}^d\int
  \Delta \(\Lambda^s v^\eps_k \)\Lambda^s v^\eps_k = -\frac{\eps^2}{2}\int
  \lvert \nabla \Lambda^s v^\eps \rvert^2.
\end{equation*}
For the second new term, we have
\begin{equation*}
  2\eps\<S\Lambda^s \(B^j(U^\eps)\d_j U^\eps\),\Lambda^s U^\eps\> = 
\frac{\eps}{2}\int \Lambda^s\(a_2^\eps \d_j v_j^\eps\)\Lambda^s a_1^\eps- \frac{\eps}{2}\int
\Lambda^s\(a_1^\eps \d_j v_j^\eps\)\Lambda^s a_2^\eps .
\end{equation*}
By Kato--Ponce commutator estimates \cite{katoponce},
\begin{equation*}
  \left\lVert \Lambda^s\(a_k^\eps \d_j v_j^\eps\)
    - a_k^\eps \Lambda^s \d_j v_j^\eps \right\rVert_{L^2} \lesssim
  \|\nabla a_k^\eps\|_{L^\infty} \|v_j^\eps\|_{H^s}+
  \|a_k^\eps\|_{H^s}\|\nabla v^\eps\|_{L^\infty}. 
\end{equation*}
Then using Cauchy--Schwarz and Young inequalities, 
\begin{equation*}
  \left\lvert\int   a_k^\eps \Lambda^s \d_j v_j^\eps \Lambda^s
    a_\ell^\eps\right\rvert\le  
 \|\nabla \Lambda^s v^\eps\|_{L^2}\|a_k^\eps\Lambda^s
    a_\ell^\eps\|_{L^2} \le \delta \eps \|\nabla \Lambda^s
    v^\eps\|_{L^2}^2 + \frac{1}{4\delta\eps}\|a_k^\eps\Lambda^s
    a_\ell^\eps\|_{L^2}^2. 
\end{equation*}
Choosing $\delta>0$ sufficiently small and independent of $\eps$, we
infer
\begin{equation}\label{eq:apriorigen}
  \frac{d}{dt}\<\Lambda^s (SU^\eps),\Lambda^s U^\eps\>+
  \frac{\eps^2}{4}\|\nabla \Lambda^s 
    v^\eps\|_{L^2}^2  \lesssim
  \|U^\eps\|_{W^{1,\infty}}\|U^\eps\|_{H^s}^2 +
  \|a^\eps\|_{L^\infty}^2 \|a^\eps\|_{H^s}^2.
\end{equation}
The local existence result then follows easily by adapting the standard arguments from \cite{Taylor3}. In view of the energy estimate \eqref{eq:apriorigen}, we also infer \eqref{lifespan}. The lifespan $T^\eps_{\rm max}$ of the maximal solution to \eqref{eq:bcm} must be expected to actually depend on $\eps$. For instance, if $\phi_0,a_0\in C_0^\infty(\R^d)$, then from
  \cite{MUK86}, $T^0_{\rm max}<\infty$, regardless of the dimension $d$. On
  the other hand, we prove further (proof of Theorem \ref{theo:global}) that if $d=1$, $T^\eps_{\rm max}=\infty$ for all $\eps>0$.
  
  \bs
  \ni
  {\em Item (iii): $u^\eps$ is solution of the original NLS problem.}  Define $\phi^\eps$ by \eqref{eq:phi}. Let us prove that $(a^\eps,\phi^\eps)\in C([0,T_{\rm max}^\eps);H^s\times H^{s+1})$ solves the system \eqref{eq:bcm0}. Since $(a^\eps,v^\eps)\in C([0,T_{\rm max}^\eps);H^s\times H^s)$, we
readily have $\phi^\eps \in C([0,T_{\rm
  max}^\eps);H^{s-1})$. Moreover, $v^\eps$ is irrotational: indeed,
$\mathrm{curl }\, v^\eps$ satisfies a homogeneous equation (see
e.g. \cite[p.~291]{BCD11}),  it is zero at time $t=0$, so
$\mathrm{curl }\, v^\eps\equiv 0$. 
Set $\Omega^\eps = D v^\eps -\nabla v^\eps$, where $D v^\eps$ stands
for the Jacobian matrix of $v^\eps$, and $\nabla v^\eps$ stands for
its transposed matrix. It solves
\begin{equation*}
  \d_t \Omega^\eps +v^\eps\cdot \nabla \Omega^\eps +\Omega^\eps\cdot D
  v^\eps + \nabla v^\eps \cdot \Omega^\eps = \eps^2 \Delta \Omega^\eps. 
\end{equation*}
As a consequence, $\nabla |v^\eps|^2 = 2v^\eps\cdot \nabla v^\eps$. 
We then deduce from \eqref{eq:phi} and from the second equation of
\eqref{eq:bcm}  
\begin{equation*}
  \d_t\(\nabla \phi^\eps -v^\eps\)= \nabla \d_t \phi^\eps -\d_t v^\eps =0,
\end{equation*}
so we infer from the initial condition $v^\eps_{\mid
  t=0}=\nabla\phi_0=\na\phi^\eps_{\mid t=0}$ that $v^\eps = \nabla
\phi^\eps$. Therefore, $\nabla 
\phi^\eps\in C([0,T_{\rm max}^\eps);H^s)$, hence 
$$\phi^\eps \in
C([0,T_{\rm max}^\eps);H^{s+1}).$$
Replacing $v^\eps$ with
$\nabla \phi^\eps$ in \eqref{eq:phi}  shows that $\phi^\eps$ solves the
first equation in \eqref{eq:bcm0}, and the second equation in
\eqref{eq:bcm0} is obviously satisfied. 

Now, define $u^\eps$ by \eqref{madelung}. The property $u^\eps \in C([0,T^\eps];H^{s})$ is straightforward,
  since $H^s(\R^d)$ is an algebra. It is clear that $u^\eps$ solves
  \eqref{eq:nls} and \eqref{eq:ci}. Since 
$$C([0,T^\eps];H^s(\R^d))\subset
  C([0,T^\eps];H^1(\R^d))\cap L^{8/d}([0,T^\eps];W^{1,4}(\R^d)),$$
  uniqueness stems from Proposition~\ref{prop:nls}. 
\end{proof}

\subsection{Convergence towards the Euler equation}
\label{sec:cvEuler}
In this subsection, we study the semiclassical limit of \eqref{eq:bcm}.
\begin{proof}[Proof of Theorem \pref{theo:cveuler}]
  The proof proceeds as for
  \cite[Theorem~1.2]{Grenier98}. First, Theorem~\ref{theo:local}
  yields a solution $(a,v)\in C([0,T^0_{\rm max});H^s\times 
  H^s)$. Necessarily, $T^0_{\rm max}>T$ (where $T$ is an existence
  time for the Euler equation), for if we had $T^0_{\rm 
    max}\le T$, then by uniqueness for \eqref{eq:euler},
  $(|a|^2,v)=(\rho,v)\in L^\infty([0,T];H^s\times H^s)$. The
  equation for $a$ in \eqref{eq:eulersym} 
  then shows that $a\in L^\infty([0,T];H^{s-1})$, hence $(a,v)\in
  L^\infty([0,T];H^{s-1}\times H^{s})$, which is impossible if $ T^0_{\rm
    max}\le T$: this yields the first part of the theorem.

To prove the error estimate, introduce the vector-valued unknown $U$,
associated to $(a,v)$. We know that $U \in 
  C([0,T];H^{s})$. Consider the error $W^\eps = U^\eps
  -U$. It solves
  \begin{equation*}
    \d_t W^\eps + \sum_{j=1}^d \( A^j(U^\eps)\d_j U^\eps
    -A^j(U)\d_j U\) =\frac{\eps}{2}LU^\eps+\eps^2 DU^\eps + \eps
  \sum_{j=1}^dB^j(U^\eps)\partial_j U^\eps .   
  \end{equation*}
Rewrite this equation as:
\begin{align*}
\d_t W^\eps + \sum_{j=1}^d A^j(U^\eps)\d_j W^\eps
 & =-\sum_{j=1}^d \(A^j(U^\eps) -A^j(U)\)\d_j U+
  \frac{\eps}{2}LW^\eps+\frac{\eps}{2} LU \\
&\quad +\eps^2 DW^\eps+\eps^2DU + \eps
  \sum_{j=1}^dB^j(U^\eps)\partial_j U^\eps .     
\end{align*}
The operator on the left hand side is the same operator as in
\eqref{eq:systhyp}.  The term $LW^\eps$ is not present in the energy
estimates, since it is skew-symmetric. The term $\eps LU$ is
considered as a source term: it is of order $\eps$, uniformly in
$L^\infty([0,T];H^{s-2})$. Next, the first term on the right hand side is a
semi-linear perturbation:
\begin{align*}
  \left\| \(A^j(U^\eps) -A^j(U)\)\d_j U \right\|_{H^{s-2}} &\le 
\left\| \(A^j(U^\eps) -A^j(U)\)\right\|_{H^{s-2}} \left\|U
\right\|_{H^{s-1}}\\
&\lesssim \left\| \(A^j(W^\eps+U) -A^j(U)\)\right\|_{H^{s-2}} \\
&\le C\( \|W^\eps\|_{L^\infty}, \|U\|_{H^{s-1}}\) \|W^\eps\|_{H^{s-2}},
\end{align*}
where we have used tame estimates. Finally, we know that 
$W^\eps$ is bounded in $L^\infty([0,T]\times\R^d)$, as the
difference of two bounded terms. The terms involving $D$ and $B$ are
treated in a similar way, by adapting the estimates used in the proof
of Theorem~\ref{theo:local}. We end up with
\begin{align*}
  \frac{d}{dt}\<S \Lambda^{s-2} W^\eps,\Lambda^{s-2}
  W^\eps\>+\frac{\eps^2}{4}\|\nabla \Lambda^{s-2} W^\eps\|_{L^2}^2& \lesssim  
\eps\|W^\eps\|_{H^{s-2}} +
  \|W^\eps\|_{H^{s-2}}^2\\
&\lesssim \eps^2 + \<S \Lambda^{s-2} W^\eps,\Lambda^{s-2} W^\eps\> ,  
\end{align*}
and the result follows from Gronwall lemma. 
\end{proof}

\subsection{Global existence result in one dimension}
\label{sec:global}
 In this section, we consider \eqref{eq:bcm} with $\eps>0$ fixed, in dimension $d=1$. To prove a global existence result, we shall need finer estimates than \eqref{eq:apriorigen}, in particular for low order derivatives. 

\begin{proof}[Proof of Theorem \pref{theo:global}] For $s=0$ and
  $s=1$, \eqref{eq:bcm} can be refined, and this will be useful to
  prove a global existence result. 

\subsubsection*{$L^2$ estimates}

Multiplying the second equation in \eqref{eq:bcm} by
$\overline{a^\eps}$, integrating in space, and considering the real
part yields:
\begin{equation*}
  \frac{d}{dt}\|a^\eps(t)\|_{L^2}^2=0.
\end{equation*}
We have simply recovered the conservation of mass \eqref{eq:masse}. 
\smallbreak

Multiplying now the first equation in \eqref{eq:bcm} by $v^\eps$ and
integrating in space, we take advantage of the dimensional one:
\be
\frac{1}{2}\frac{d}{dt}\|v^\eps(t)\|_{L^2}^2 +
\eps^2\|\d_x
  v^\eps(t)\|_{L^2}^2 \lesssim
   \int |a^\eps(t,x)|^2|\d_x v^\eps(t,x)|dx.\label{d=1}
   \ee

\subsubsection*{$\dot H^1$ estimates}
Differentiating \eqref{eq:bcm} in space, multiplying then the first
equation by $\d_x v^\eps$, the second by $\d_x \overline{a^\eps}$,
and integrating, we find (thanks to the symmetrizer):
\begin{align}
  \frac{d}{dt}&\( \|\pa_x v^\eps(t)\|_{L^2}^2 + 4\|\pa_x
  a^\eps(t)\|_{L^2}^2\) +\eps^2\|\pa^2_x v^\eps(t)\|_{L^2}^2\nonumber\\
 & \lesssim   \(\|\pa_x
  v^\eps(t)\|_{L^\infty}+\|a^\eps(t)\|_{L^\infty}+\|a^\eps(t)\|_{L^\infty}^2\)\(
  \|v^\eps(t)\|_{H^1}^2 + \|a^\eps(t)\|_{H^1}^2 \).\label{H1point}
\end{align}

\bs
Now, we have the tools to prove Theorem \ref{theo:global}.  In view of Remark~\ref{rem:solmax}, it suffices to prove that for
  any fixed $\eps>0$,  
  \begin{equation*}
    v^\eps,\d_x v^\eps,a^\eps,\d_xa^\eps , |a^\eps|^2\in L^1_{\rm
      loc}\([0,\infty);L^\infty(\R^d)\).  
  \end{equation*}
Since $\eps>0$ is fixed, we shall omit it in the notations, and fixed
it equal to $1$ in \eqref{eq:bcm}. Using Item  {\em (iii)} of Theorem \ref{theo:local}, and in view of \eqref{eq:energy}, we
have the {\it a priori} estimate
\begin{equation*}
  \|a(t)\|_{L^4}\lesssim 1,\quad \forall t\ge 0. 
\end{equation*}
The refined $L^2$ estimate \eqref{d=1} for $v$ and
Cauchy--Schwarz inequality yield
\begin{equation*}
  \frac{1}{2}\frac{d}{dt}\|v(t)\|_{L^2}^2 +
\|\d_x  v(t)\|_{L^2}^2 \lesssim \|a(t)\|_{L^4}^2 \|\d_x
v(t)\|_{L^2}\lesssim \|\d_x v(t)\|_{L^2}.
\end{equation*}
Using Young's inequality
\begin{equation*}
   \|\d_x v(t)\|_{L^2} \le \delta  \|\d_x v(t)\|_{L^2}^2 +\frac{1}{4\delta},
\end{equation*}
we infer, for $\delta>0$ sufficiently small,
\begin{equation*}
 \frac{d}{dt}\|v(t)\|_{L^2}^2 +
\|\d_x  v(t)\|_{L^2}^2 \lesssim 1. 
\end{equation*}
Therefore, 
\begin{equation*}
  \|v(t)\|_{L^2}^2 +\int_0^t \|\d_x  v(\tau)\|_{L^2}^2d\tau \lesssim
  1+t. 
\end{equation*}
The Gagliardo--Nirenberg inequality yields
\begin{align*}
  \int_0^t \|v(\tau)\|_{L^\infty}d\tau& \lesssim \int_0^t
  \|v(\tau)\|_{L^2}^{1/2}\|\d_x v(\tau)\|_{L^2}^{1/2}d\tau \lesssim
  (1+t)^{1/4} \int_0^t\|\d_x v(\tau)\|_{L^2}^{1/2}d\tau \\
& \lesssim  (1+t)^{1/4}  \int_0^t\(1+\|\d_x
v(\tau)\|_{L^2}^{2}\)d\tau\lesssim (1+t)^{5/4}, 
\end{align*}
where we have used the Young inequality. In view of \eqref{eq:energy}, we
infer
\begin{equation*}
  \|\d_x a(t)\|_{L^2}\lesssim 1 + \|a(t)v(t)\|_{L^2}\lesssim 1 +
  \|v(t)\|_{L^\infty},
\end{equation*}
and the Gagliardo--Nirenberg inequality yields
\begin{equation*}
  \|a(t)\|_{L^\infty} \lesssim \| a(t)\|_{L^2}^{1/2}\|\d_x
  a(t)\|_{L^2}^{1/2} \lesssim 1 + \|v(t)\|_{L^\infty}^{1/2}. 
\end{equation*}
We infer 
\begin{equation*}
  \int_0^t \(\|a(\tau)\|_{L^\infty}+  \|a(\tau)\|_{L^\infty}^2\)
d\tau\lesssim (1+t)^{5/4}.
\end{equation*}
It only remains to prove that $\d_x v, \d_x a\in L^1_{\rm
  loc}([0,\infty);L^\infty(\R)).$ The $\dot H^1$ estimate \eqref{H1point} yields
\begin{align*}
 \int_0^t \|\d_x^2 &v(\tau)\|_{L^2}^2d\tau \lesssim 1+ \int_0^t\(
 \|\d_x v(\tau)\|_{L^\infty} +1 +
 \|a(\tau)\|_{L^\infty}^2\)\(\|a(\tau)\|_{H^1}^2
 +\|v(\tau)\|_{H^1}^2\)d\tau \\
&\lesssim 1 + \int_0^t\(
 \|\d_x v(\tau)\|_{L^\infty} +1 +
 \|a(\tau)\|_{L^\infty}^2\)\(1+\|v(\tau)\|_{L^\infty}^2
 +\|v(\tau)\|_{H^1}^2\)d\tau \\
&\lesssim 1 + \int_0^t\(
 \|\d_x v(\tau)\|_{L^\infty} +1 +
 \|a(\tau)\|_{L^\infty}^2\)\(1
 +\|v(\tau)\|_{H^1}^2\)d\tau,
\end{align*}
from Sobolev embedding. Using an integration by parts and the Young
inequality, we find
\begin{align*}
\|\d_xv(\tau)\|_{L^2}^2&\le \|v(\tau)\|_{L^2}\|\d_x^2 v(\tau)\|_{L^2}\\
&\le \frac{\delta}{ \|\d_x v(\tau)\|_{L^\infty} +1 +
 \|a(\tau)\|_{L^\infty}^2} \|\d_x^2 v(\tau)\|_{L^2}^2\\
&\quad  + 
\frac{ \|\d_x v(\tau)\|_{L^\infty} +1 +
 \|a(\tau)\|_{L^\infty}^2}{4\delta}\|v(\tau)\|_{L^2}^2.
\end{align*}
Taking $\delta>0$ sufficiently small, we infer:
\begin{align*}
  \int_0^t \|\d_x^2 v(\tau)\|_{L^2}^2d\tau &\lesssim 1 + 
(1+t)\int_0^t\(
 \|\d_x v(\tau)\|_{L^\infty} +1 +
 \|a(\tau)\|_{L^\infty}^2\)d\tau\\
&\quad + (1+t)\int_0^t\(
 \|\d_x v(\tau)\|_{L^\infty} +1 +
 \|a(\tau)\|_{L^\infty}^2\)^2d\tau\\
&\lesssim 1+t^2+ (1+t)\int_0^t \|\d_x v(\tau)\|_{L^\infty}^2d\tau +
(1+t)\int_0^t \|a(\tau)\|_{L^\infty}^4d\tau .
\end{align*}
We also have
\begin{align*}
  \int_0^t \|a(\tau)\|_{L^\infty}^4d\tau &\lesssim t + \int_0^t
  \|v(\tau)\|_{L^\infty}^2d\tau \lesssim t + \int_0^t
  \|v(\tau)\|_{L^2}\|\d_xv(\tau)\|_{L^2} d\tau\\ 
&\lesssim t + \sqrt{1+t}\int_0^t
  \|\d_xv(\tau)\|_{L^2} d\tau\lesssim (1+t)^{3/2},
\end{align*}
and so
\begin{equation*}
 \int_0^t \|\d_x^2 v(\tau)\|_{L^2}^2d\tau \lesssim (1+t)^{5/2} +
 (1+t)\int_0^t \|\d_x v(\tau)\|_{L^\infty}^2d\tau .  
\end{equation*}
Now Gagliardo--Nirenberg and Cauchy--Schwarz inequalities yield
\begin{align*}
 \int_0^t \|\d_x v(\tau)\|_{L^\infty}^2d\tau &\lesssim 
\int_0^t \|\d_x v(\tau)\|_{L^2}\|\d_x^2 v(\tau)\|_{L^2}d\tau \\
& \lesssim \(\int_0^t \|\d_x v(\tau)\|_{L^2}^2d\tau\)^{1/2}\(\int_0^t
\|\d_x^2 v(\tau)\|_{L^2}^2d\tau\)^{1/2} \\
&\lesssim \sqrt{1+t}\((1+t)^{5/2} + (1+t)\int_0^t \|\d_x
v(\tau)\|_{L^\infty}^2d\tau\)^{1/2} \\
&\lesssim (1+t)^{7/4} + (1+t)\(\int_0^t \|\d_x
v(\tau)\|_{L^\infty}^2d\tau\)^{1/2}. 
\end{align*}
The Young inequality implies
\begin{equation*}
  \int_0^t \|\d_x v(\tau)\|_{L^\infty}^2d\tau \lesssim (1+t)^2,
\end{equation*}
hence  $\d_x v\in L^2_{\rm loc}([0,\infty);L^\infty)\subset L^1_{\rm
  loc}([0,\infty);L^\infty)$. 

Since $\d_x a \in L^1_{\rm loc}([0,\infty);L^2)$,
Gagliardo--Nirenberg shows that it suffices to prove that $\d_x^2 a
\in L^1_{\rm loc}([0,\infty);L^2)$ to conclude that $\d_x a \in
L^1_{\rm loc}([0,\infty);L^\infty)$. For that, we use again  Item  {\em (iii)} of Theorem \ref{theo:local}, from which we now that $a e^{i\phi}
\in L^1_{\rm 
  loc}([0,\infty);H^2)$. Differentiating $ a e^{i\phi}$ twice, we
have, since $\d_x \phi=v$,
\begin{align*}
  \|\d_x^2 a(t)\|_{L^2}&\lesssim \|a e^{i\phi}\|_{H^2} + \|v \d_x
  a\|_{L^2} + \|a\d_x v \|_{L^2} + \|v^2 a\|_{L^2}\\
&\lesssim \|a e^{i\phi}\|_{H^2} + \|v\|_{L^\infty}\|\d_xa\|_{L^2} +
\|a \|_{L^\infty}\|\d_x v \|_{L^2}+ \|v\|_{L^\infty}^2.  
\end{align*}
\mbox{}From the above estimate, we infer $\d_x^2 a
\in L^1_{\rm loc}([0,\infty);L^2)$, which completes the proof of the
theorem.  
\end{proof}

\begin{remark}[Higher dimension]
  Even though we have obviously taken advantage of the one-dimensional
  setting to use the embedding $H^1\subset L^\infty$, the true reason
  why Theorem~\ref{theo:global} is limited to $d=1$ lies
  elsewhere. The refined $L^2$ estimate \eqref{d=1} becomes, if $d\ge
  2$,
\begin{equation*}
\frac{1}{2}\frac{d}{dt}\|v^\eps(t)\|_{L^2}^2 +
\eps^2\|\nabla 
  v^\eps(t)\|_{L^2}^2 \lesssim
\int \(|a^\eps(t,x)|^2+|v^\eps(t,x)|^2\)|\nabla v^\eps(t,x)|dx.
\end{equation*}
Note the appearance of the new term $\int |v^\eps(t,x)|^2|\nabla
v^\eps(t,x)|dx$: it correspond to the fact that $\int v^\eps \cdot
\(v^\eps\cdot \nabla v^\eps\)$ is zero if $d=1$, but not if $d\ge
2$ (in general). This aspect seems to be the only real limitation to extend
Theorem~\ref{theo:global} to $d\ge 2$. 
\end{remark}

\section{Numerics}
\label{sec:numerique}

In this section, we define a second order numerical scheme for
\eqref{eq:bcm}. Since $\eps$ is no longer a singular perturbation
parameter in this reformulation of NLS, this scheme is naturally
Asymptotic Preserving for our original problem in the  semiclassical
limit $\eps\to 0$, as long as the solution of the Euler equation
remains smooth. Since we solve our problem in dimension 1 on a bounded
interval and in dimension 2 on a rectangle, we add periodic boundary
conditions to \eqref{eq:nls} or to \eqref{eq:bcm}, both formulations
remaining equivalent. 

\subsection{The AP numerical scheme}
The nonlinear system to be solved reads
\begin{equation}
  \label{eq:syst1D}
  \left \{
    \begin{array}{ll}
\displaystyle      \partial_t a +\DIV (av) +(i \eps-\frac{1}{2}) a
\DIV v=i\frac{\eps}{2} \Delta a\ ;\qquad & a_{\mid t=0}=a_0,\\ 
\displaystyle \partial_t v+\nabla
\left(\frac{|v|^2}{2}+|a|^2\right)=\eps^2 \Delta v\ ;\qquad & v_{\mid
  t=0}=v_0. 
    \end{array}
\right .
\end{equation}
The semi-discretization in time is realized through a Strang splitting
scheme. Let us denote by $\Delta t_n=t_n-t_{n-1}$ the variable time
step, such that $t_n=\sum_{k=1}^{n} \Delta t_k$. On each time step
$[t^n,t^{n+1}]$, we split \eqref{eq:syst1D} into two subsystems and
apply the second-order Strang splitting algorithm. 

\ms
\ni
Step 1 for $t \in [t_n, t_n+\frac{\dt_{n+1}}{2}]$:
\begin{equation*}
  \left \{
    \begin{array}{ll}
\displaystyle      \partial_t a_1 =i\frac{\eps}{2} \Delta  a_1\ ; &
a_{1\mid t=t_n}=a_{\mid t=t_n},\\[2mm] 
\displaystyle \partial_t v_1=\eps^2\Delta  v_1\ ; & v_{1\mid t=t_n}=v_{\mid t=t_n}.
    \end{array}
\right .  
\end{equation*}

\ms
\ni
Step 2 for $t \in [t_n, t_n+\dt_{n+1} ]$:
\begin{equation*}
  \left \{
    \begin{array}{ll}
\displaystyle      \partial_t a_2 +\DIV (a_2v_2) +(i \eps-\frac{1}{2}) a_2 \DIV v_2=0\ ; & a_{2\mid t=t_n}=a_{1\mid t=t_n+\frac{\dt_{n+1}}{2}},\\[2mm]
\displaystyle \partial_t v_2+\nabla_x \left(\frac{|v_2|^2}{2}+|a_2|^2\right)=0\ ; & v_{2\mid t=t_n}=v_{1\mid t=t_n+\frac{\dt_{n+1}}{2}}.
    \end{array}
\right .  
\end{equation*}

\ms
\ni
Step 3 for $t \in [t_n, t_n+\frac{\dt_{n+1}}{2}]$:
\begin{equation*}
  \left \{
    \begin{array}{ll}
\displaystyle      \partial_t a_3 =i\frac{\eps}{2} \Delta  a_3\ ; & a_{3\mid t=t_n}=a_{2\mid t=t_{n+1}},\\[2mm]
\displaystyle \partial_t v_3=\eps^2 \Delta  v_3\ ; & v_{3\mid t=t_n}=v_{2\mid t=t_{n+1}}.
    \end{array}
\right .  
\end{equation*}
Then, $(a_{3},v_{3})_{\mid t=t_n+\frac{\dt_{n+1}}{2}}$ is an approximation of $(a,v)_{\mid t=t_{n+1}}$, solution to \eqref{eq:syst1D}.\\

This splitting scheme enables to decouple the fluid part (step 2) from
the parabolic and Schr\"odinger parts (step 1 and 3), which allows a
standard treatment of each step. Note that for $\eps=0$, we recover
exactly the Euler equations \eqref{eq:eulersym}. We denote by
$(a^{n,\eps},v^{n,\eps})$ the numerical solution at time $t_n$. In the
spirit of \cite{HLR13}, the following result can be proven. 
\begin{proposition}\label{prop:un}
  Under the assumptions of Theorem~\pref{theo:local}, there exist
  $\eps_0>0$ and 
$C,c_0$ independent of $\eps\in [0,\eps_0]$ such that for all $n\in \N$ such that $t_n\in [0,T]$, where $T$ is as in   Theorem~\pref{theo:local} (ii), for  all $\Delta t_n\in (0,c_0]$, 
  \begin{align*}
   & \left\| \left\lvert a^{n,\eps}\right\rvert^2
     -\rho^\eps(t_n)   \right\|_{L^1(\R^d)\cap L^\infty(\R^d)}\le C (\max_{n}{\Delta t_n})^2,\\ 
&\left\| \eps \IM(\overline{a^{n,\eps}} \nabla a^{n,\eps})+\left\lvert a^{n,\eps}\right\rvert^2 v^{n,\eps} - \mathbf{j}^\eps(t_n) \right\|_{L^1(\R^d)\cap
  L^\infty(\R^d)}\le C (\max_{n}{\Delta t_n})^2.
 \end{align*}
\end{proposition}

We assume, without loss of generality, that we are in a periodic
framework for the space variable. This allows to solve steps 1 and 3
in a spectral way by using fast Fourier transform, whereas for step 2,
we use a Lax-Wendroff scheme (with directional splitting in dimension
2). The extension to Dirichlet or Neumann boundary conditions would
require to consider sine or cosine transforms in place of FFT. The
total scheme is consistent with \eqref{eq:syst1D}, and second-order in
time and space. Since the Lax-Wendroff scheme is explicit, a
Courant-Friedrich-Lewy condition is necessary for stability: 
$$
  \dt_n \le \mathrm{CFL}\frac{\Delta x}{\max_j{\(|v_j^n|+|a_j^n|\)}}.
$$
where $v^n_j$ denotes the approximation to $v(t_n,x_j)$, $x_j$ being a
node of the mesh discretization of the computational domain
$\Omega$, and $\Delta x$ is the space discretization parameter. In all
our numerical experiments, we take $\mathrm{CFL}=0.8$.

\subsection{Numerical experiments in dimension 1}
\begin{figure}[!htbp]
  \centering
  \subfigure[Particle density $\rho_{\mathrm{ref}}$]{\includegraphics[width=.46\textwidth]{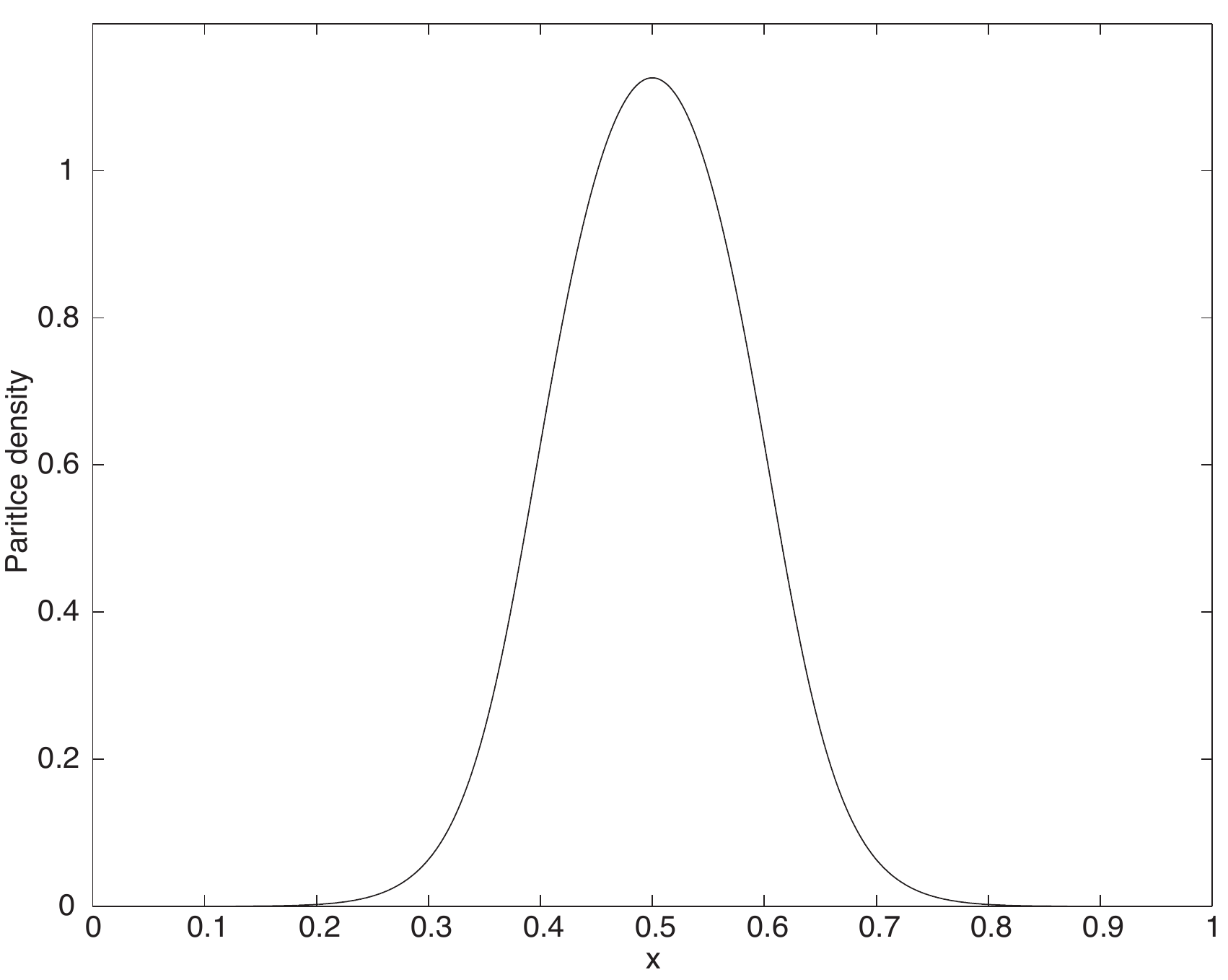}}\qquad
  \subfigure[Current density ${\mathbf j}_{\mathrm{ref}}$]{\includegraphics[width=.46\textwidth]{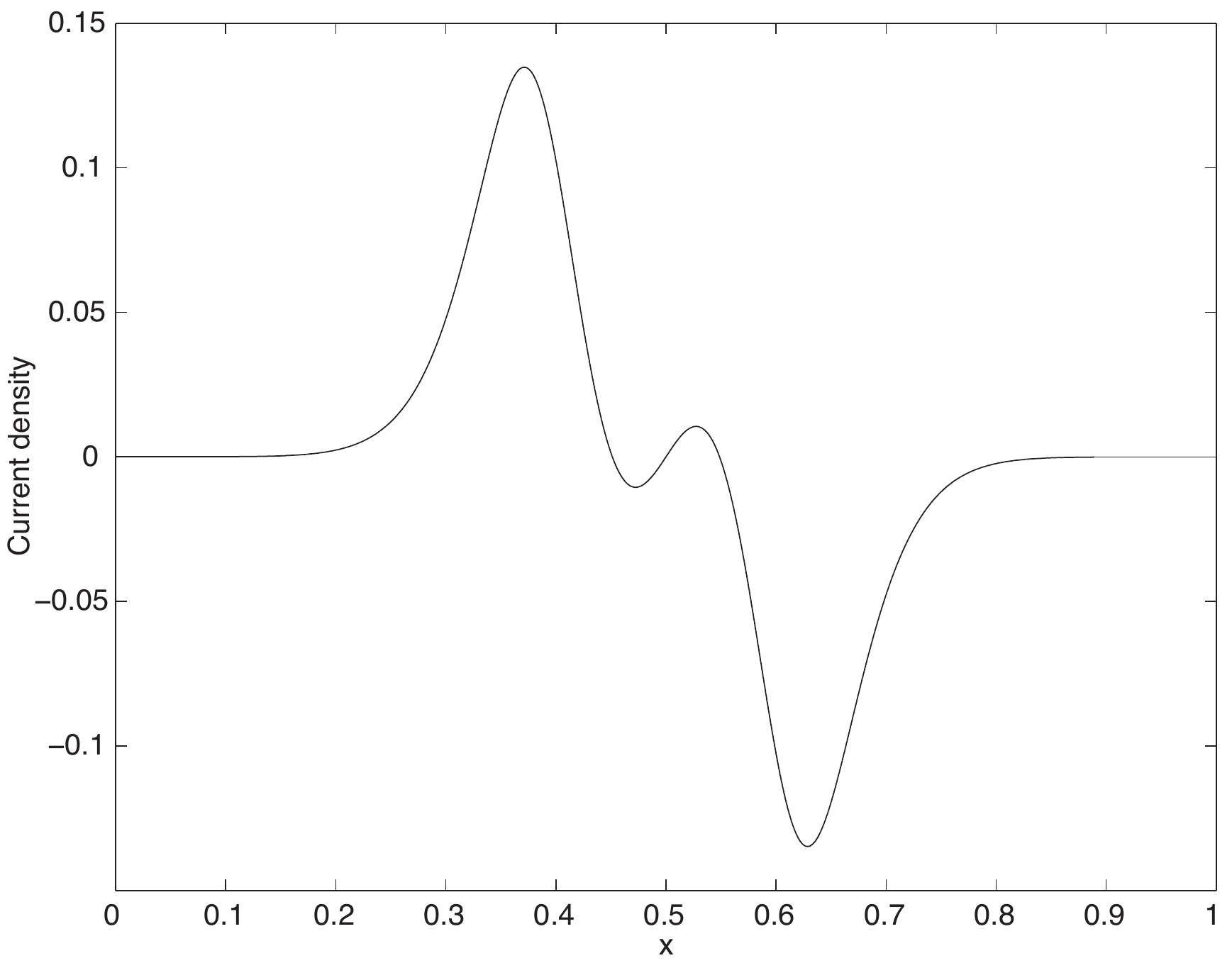}}
  
  \caption{Particle and current densities at time $T^*=0.05$ and for $\eps=10^{-4}$.}
  \label{fig:dc_t05}
\end{figure}

\begin{figure}[!htbp]
  \centering
  \subfigure[Particle density $\rho_{\mathrm{ref}}$]{\includegraphics[width=.46\textwidth]{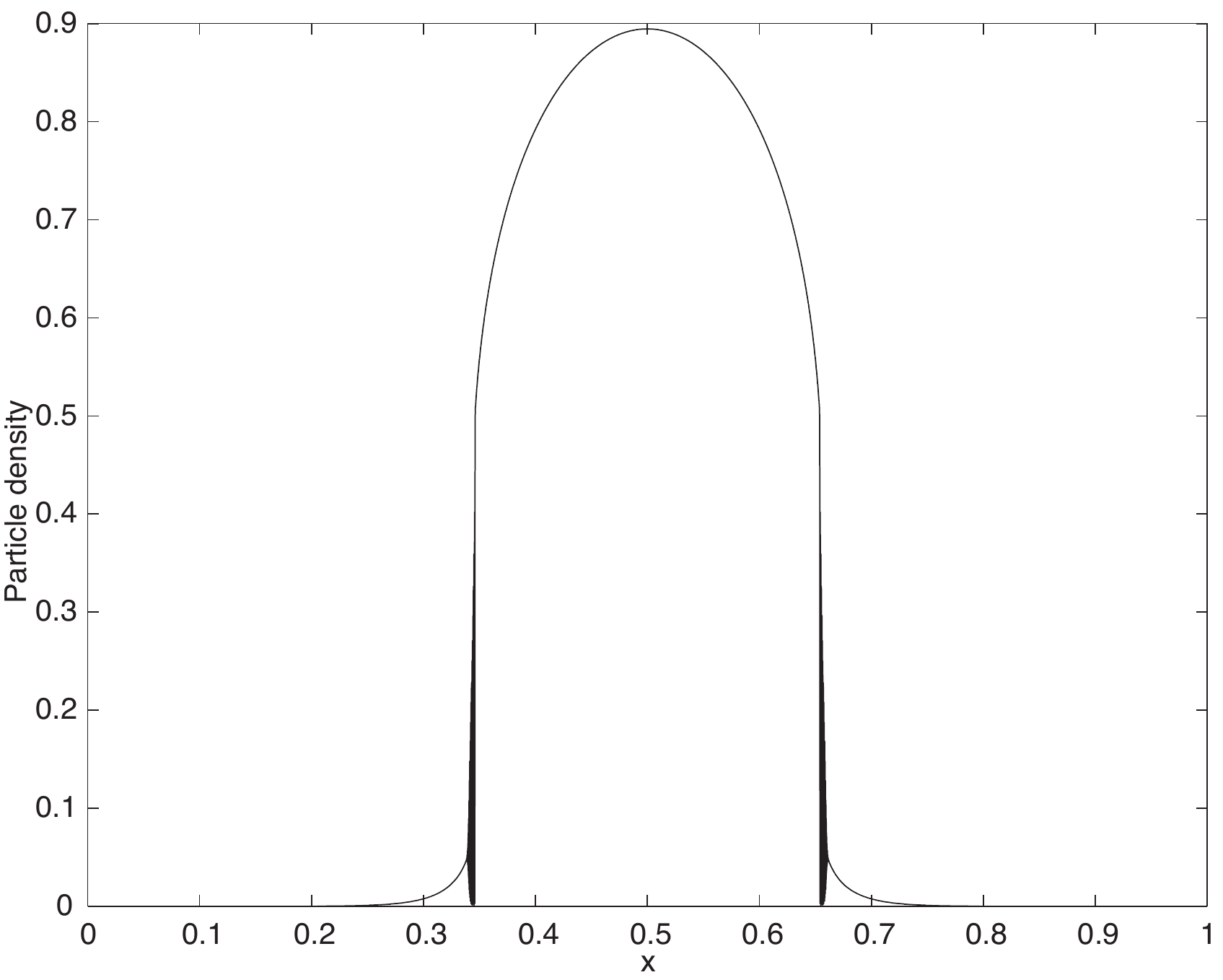}}\qquad
  \subfigure[Current density ${\mathbf j}_{\mathrm{ref}}$]{\includegraphics[width=.46\textwidth]{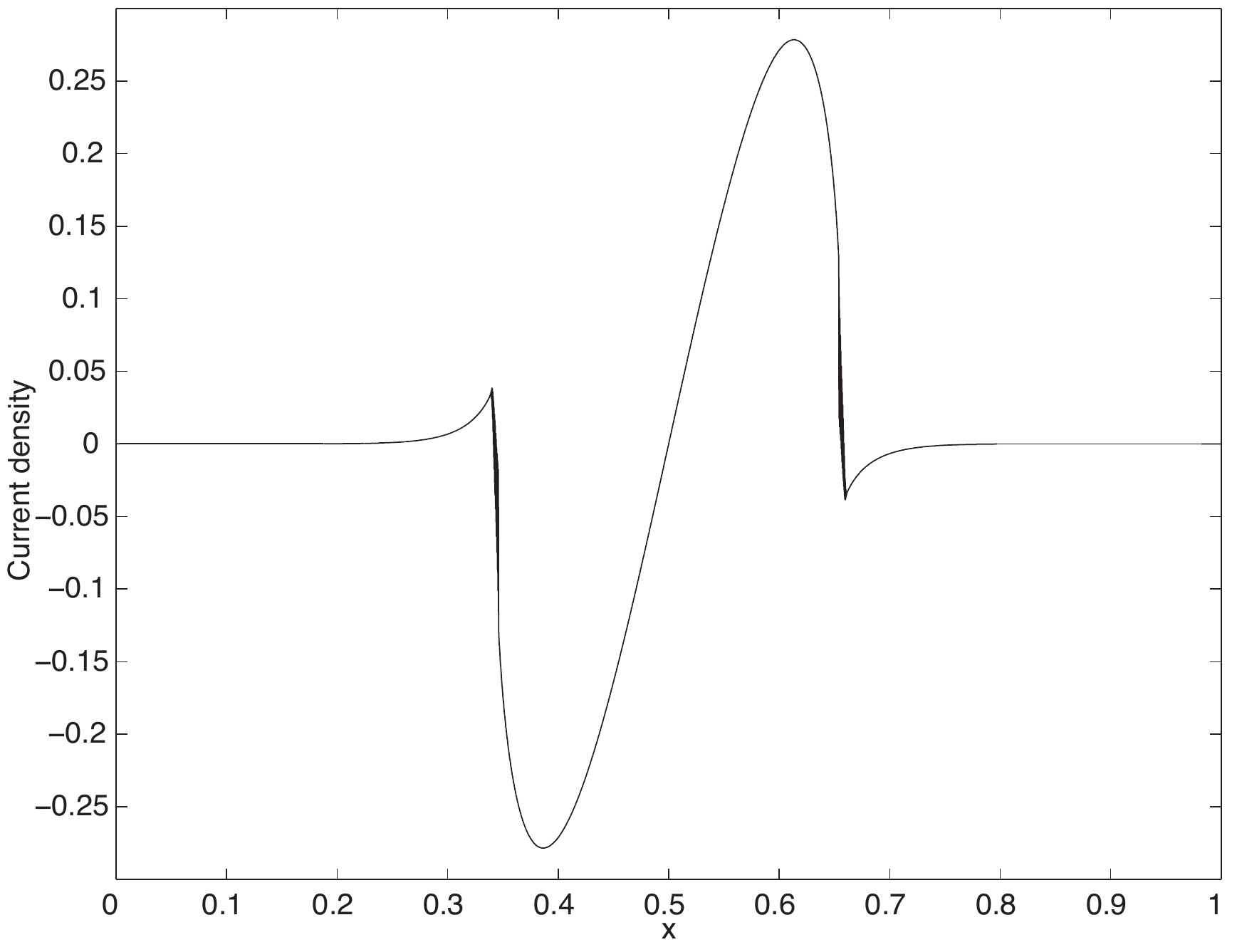}}
  
  \caption{Particle and current densities at time $T^*=0.13$ and for $\eps=10^{-4}$.}
  \label{fig:dc_t13}
\end{figure}

\begin{figure}[!htbp]
  \centering
  \subfigure[Particle density $\rho_{\mathrm{ref}}$]{\includegraphics[width=.46\textwidth]{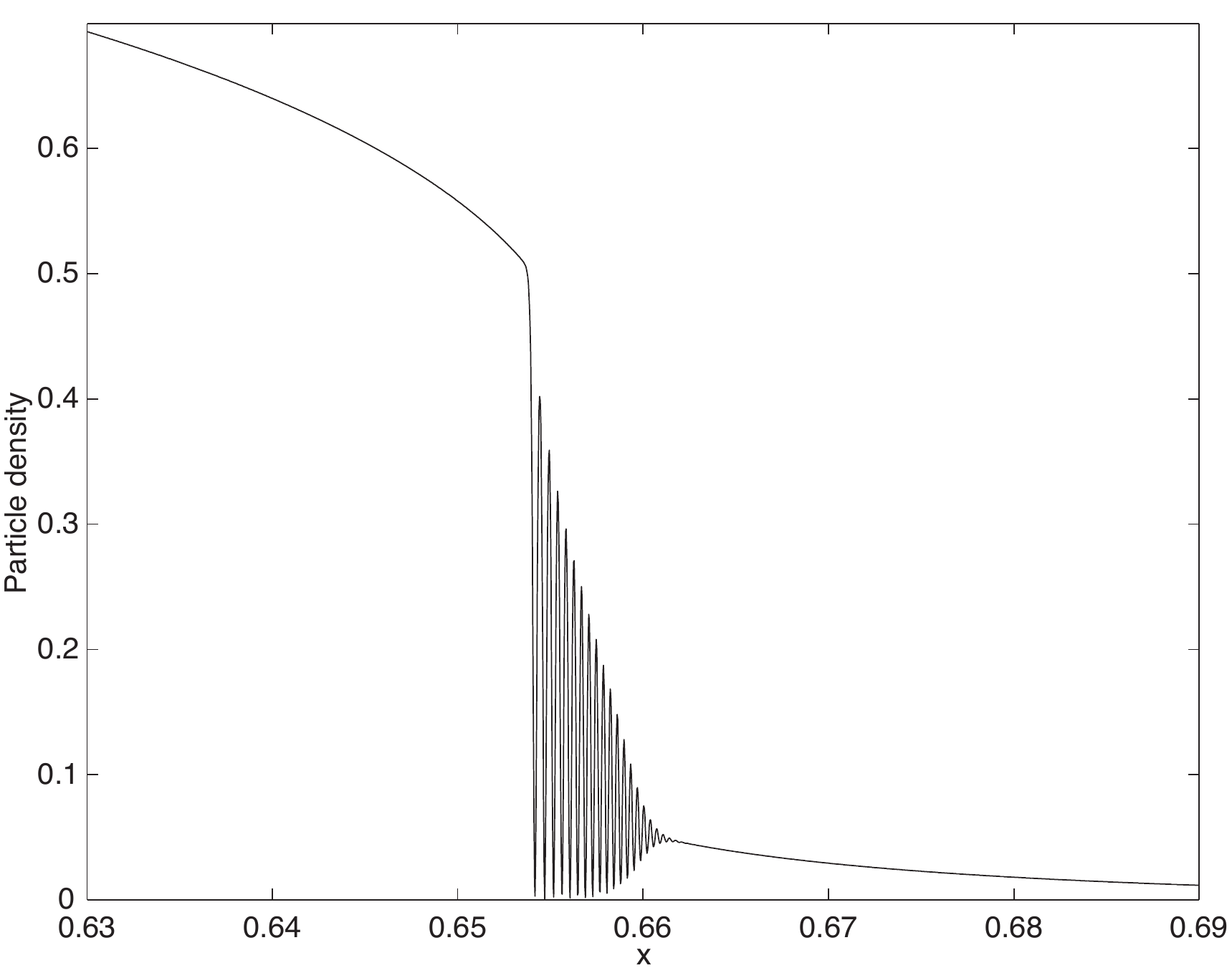}}\qquad
  \subfigure[Current density ${\mathbf j}_{\mathrm{ref}}$]{\includegraphics[width=.46\textwidth]{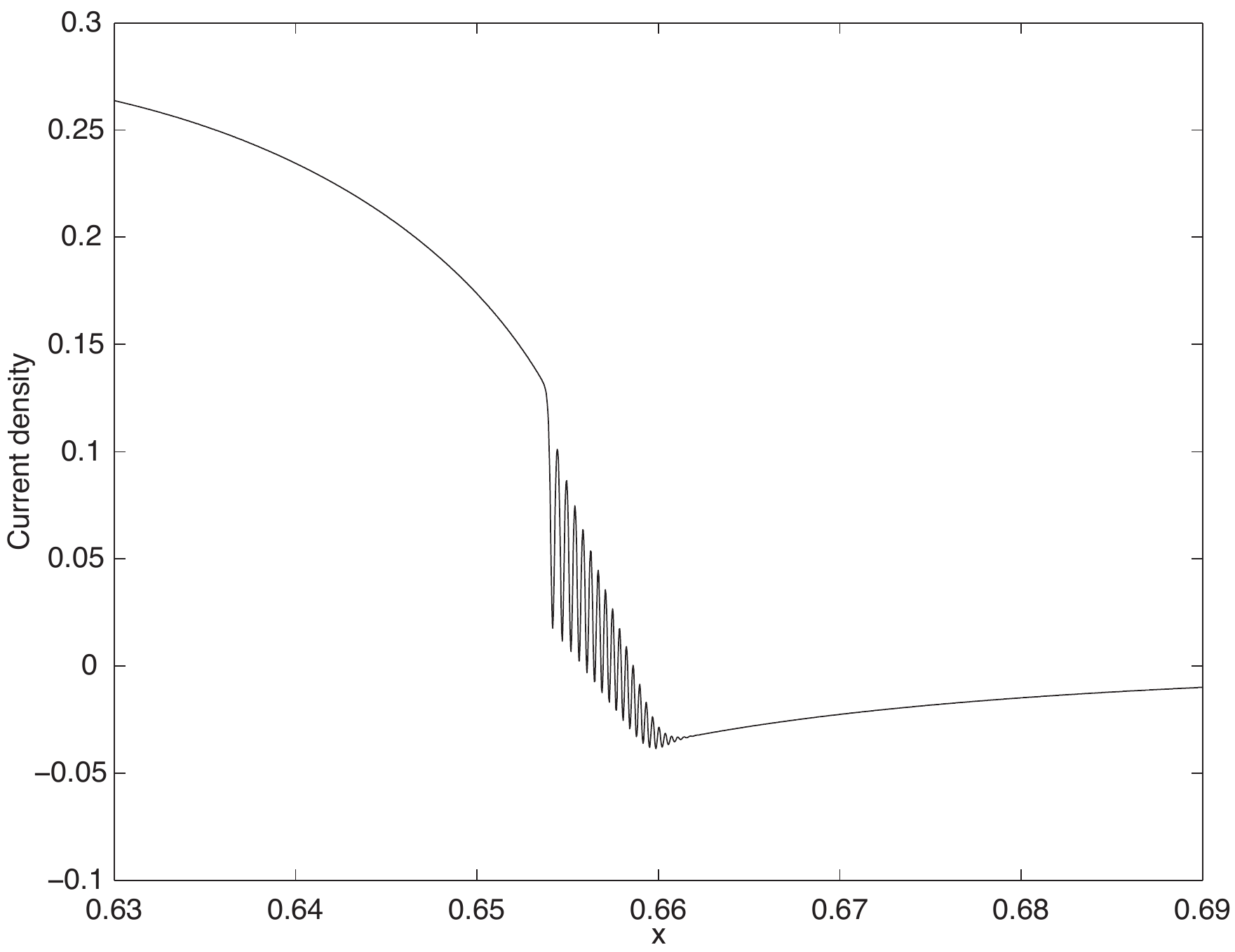}}
  
  \caption{Zoom around singularity of density and current at time $T^*=0.13$ and for $\eps=10^{-4}$.}
  \label{fig:dc_t13zoom}
\end{figure}

In order to analyze the behavior of the numerical solutions
(\ref{eq:syst1D}), we proceed like in \cite{BJM2}, and compare the
behavior of numerical physical observables (\ref{eq:observables})
computed by our AP numerical scheme to a reference solution. Since exact
solutions to (\ref{eq:nls}) with initial data (\ref{eq:ci}) are not
available, we numerically generated one thanks to the usual splitting
method described in \cite{BJM2}, applied to
\eqref{eq:nls}-\eqref{eq:ci} with a meshing strategy that ensures that
the time step $\Delta t \lesssim \eps$ and the space step $\Delta x
\lesssim \eps$. We consider the initial datum 
$$
\begin{array}{l}
a_0(x)=e^{-25(x-0.5)^2},\\
v_0(x)=-\frac{1}{5}\partial_x \ln{(e^{5(x-0.5)}+e^{-5(x-0.5)})},
\end{array}
$$
with $x \in [-1/2,3/2]$. In all our simulations, the smallest value of
$\eps$ is $10^{-4}$. We then discretize the $x$-interval with
$J=2^{15}$ subintervals and the time step is $\Delta t=\eps/100$. We
refer to the reference solution in a generic way by the notation
$u_{\rm ref}^\eps$. With the considered initial datum, new oscillations
appear in the reference solution passed a time $T^* \sim 0.10$,
meaning that singularities were created in the limit Euler system
(\ref{eq:eulersym}). In Figures~\ref{fig:dc_t05} and
\ref{fig:dc_t13}, we plot   
particle and current densities $\rho^\eps_{\mathrm{ref}}$ and
${\mathbf j}^\eps_{\mathrm{ref}}$ for $\eps=10^{-4}$ respectively at
time $t=0.05$ and $t=0.13$, before and past formation of shocks. We
present on Figure~\ref{fig:dc_t13zoom} a zoom close to singularities
which shows oscillations on physical observables. 

% Exp\'erience num\'erique 1D : on r\'ealise un test \'equivalent \`a celui pr\'esent dans le papier \cite{BaoJinMark}. Pour cela, la variable d'espace $x$ est choisie dans $[x_l,x_r]=[0,1]$ et 
% $$
% \begin{array}{l}
% a_0(x)=e^{-25(x-0.5)^2},\\
% v_0(x)=-\frac{1}{5}\partial_x \ln{(e^{5(x-0.5)}+e^{-5(x-0.5)})}.
% \end{array}
% $$
% Ne disposant pas de solution exacte, on en g\'en\`ere une en calculant une solution num\'erique par la m\'ethode de \cite{BaoJinMark} avec un pas d'espace $\Delta x=2^{-15}$ et un pas de temps $\Delta t=\eps/100$. Nous d\'enoterons par $u_{\rm ref}$ cette solution de r\'ef\'erence dans la suite de ce document. Dans la limite $\eps \to 0$, le syst\`eme converge vers le syst\`eme d'Euler (\ref{eq:eulersym}) que nous r\'esolvons par une m\'ethode standard. Nous repr\'esentons sur les figures \ref{fig:dc_t05}-\ref{fig:dc_t13zoom} la densit\'e $\rho=|a|^2$ et le courant $J=\eps \mathrm{Im}(\overline{a} \partial_x a)+\rho v$ pour $\eps=10^{-4}$ aux temps $T^*=0.05$ et $T^*=0.13$, soit respectivement avant et apr\`es la formation d'un choc pour les \'equations d'Euler.
\begin{figure}[!htbp]
  \centering

  \subfigure[Error w.r.t $J=2^M$]{\includegraphics[width=.46\textwidth]{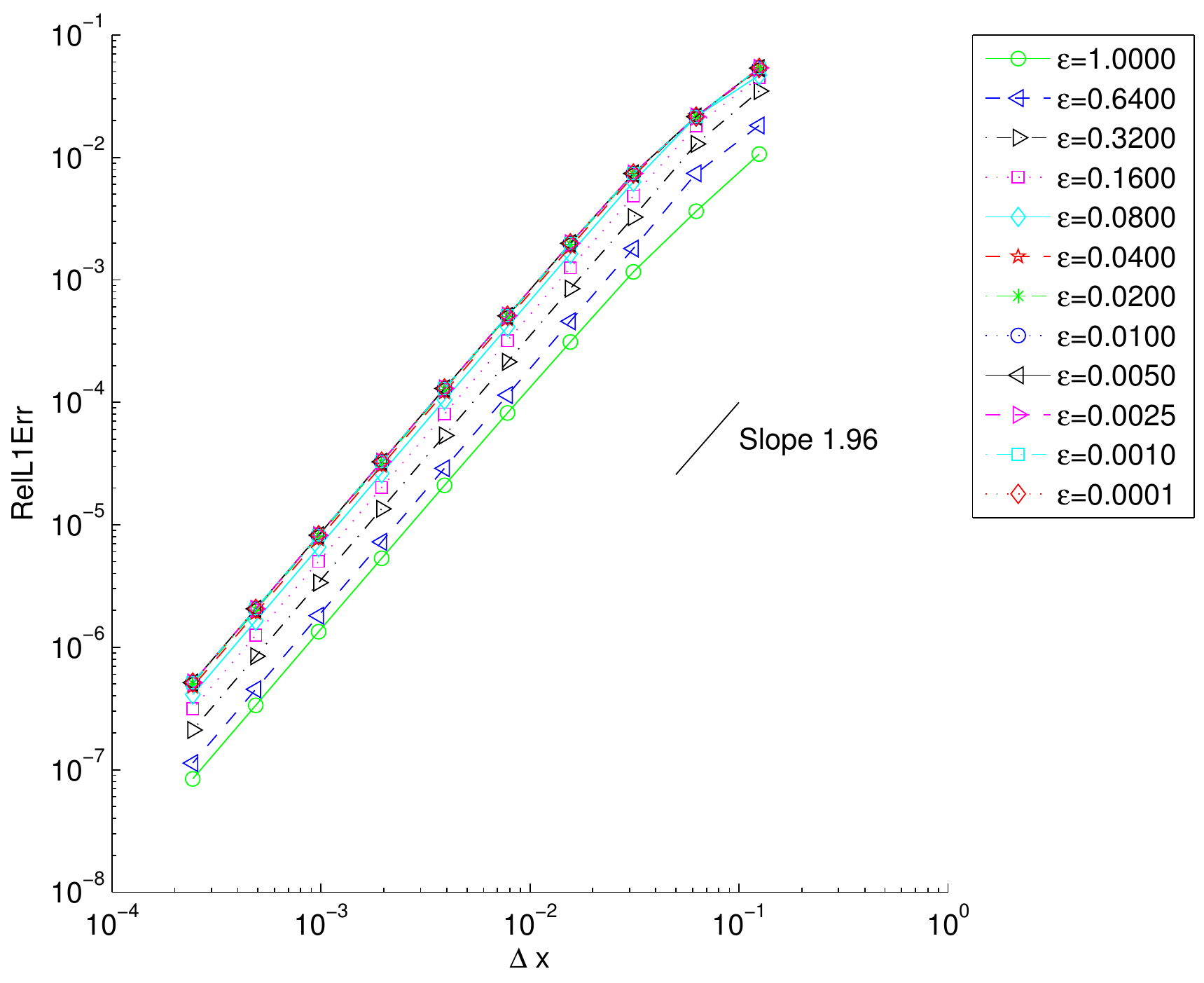}}\qquad
  \subfigure[Error w.r.t $\eps$]{\includegraphics[width=.46\textwidth]{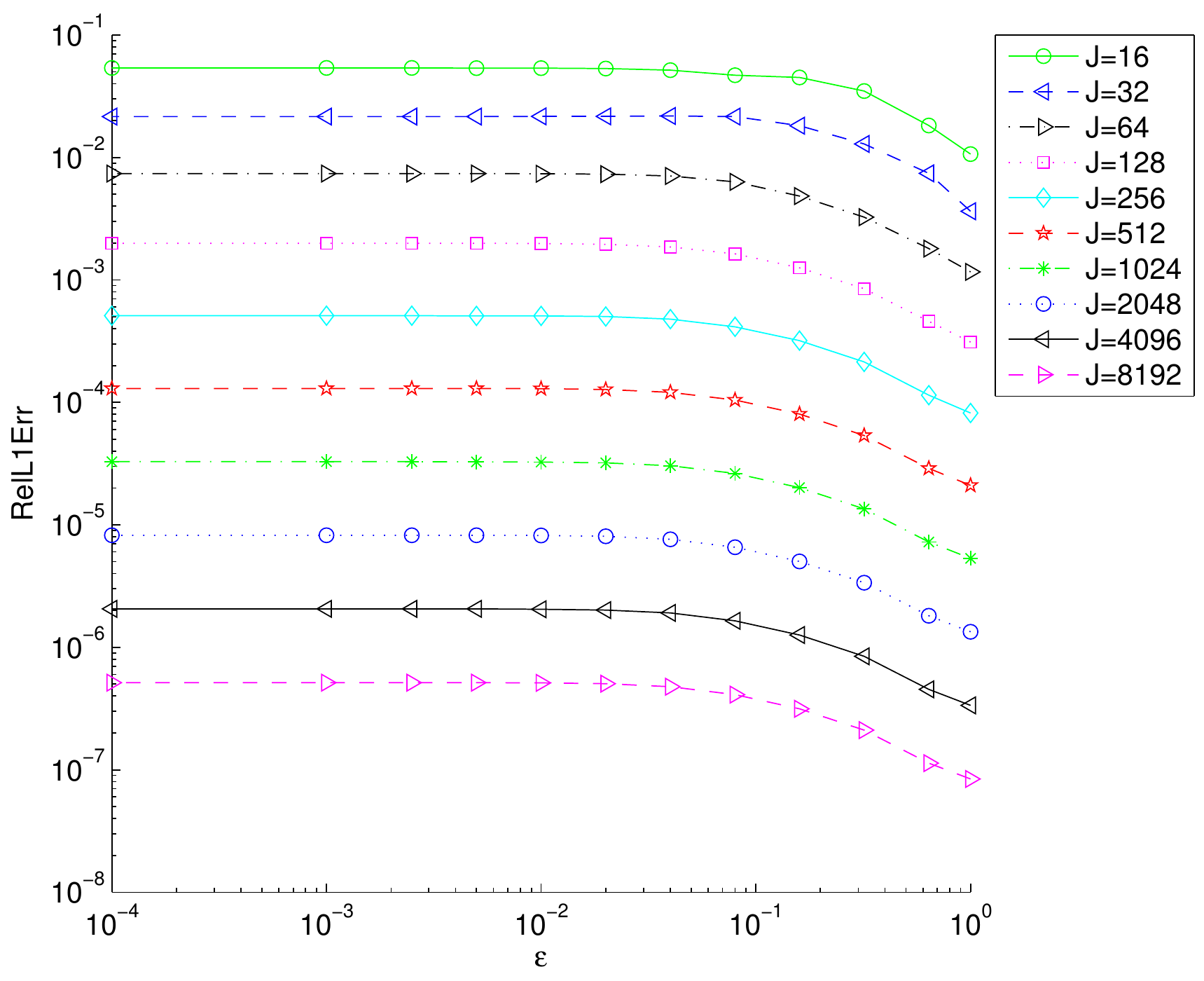}}

  \caption{$\mathrm{err}_{\rho^\eps}(t=0.05)$ for AP scheme}
\label{fig:err_rho_bef_sing_ap}  
\end{figure}

\begin{figure}[!htbp]
  \centering

  \subfigure[Error w.r.t $J=2^M$]{\includegraphics[width=.46\textwidth]{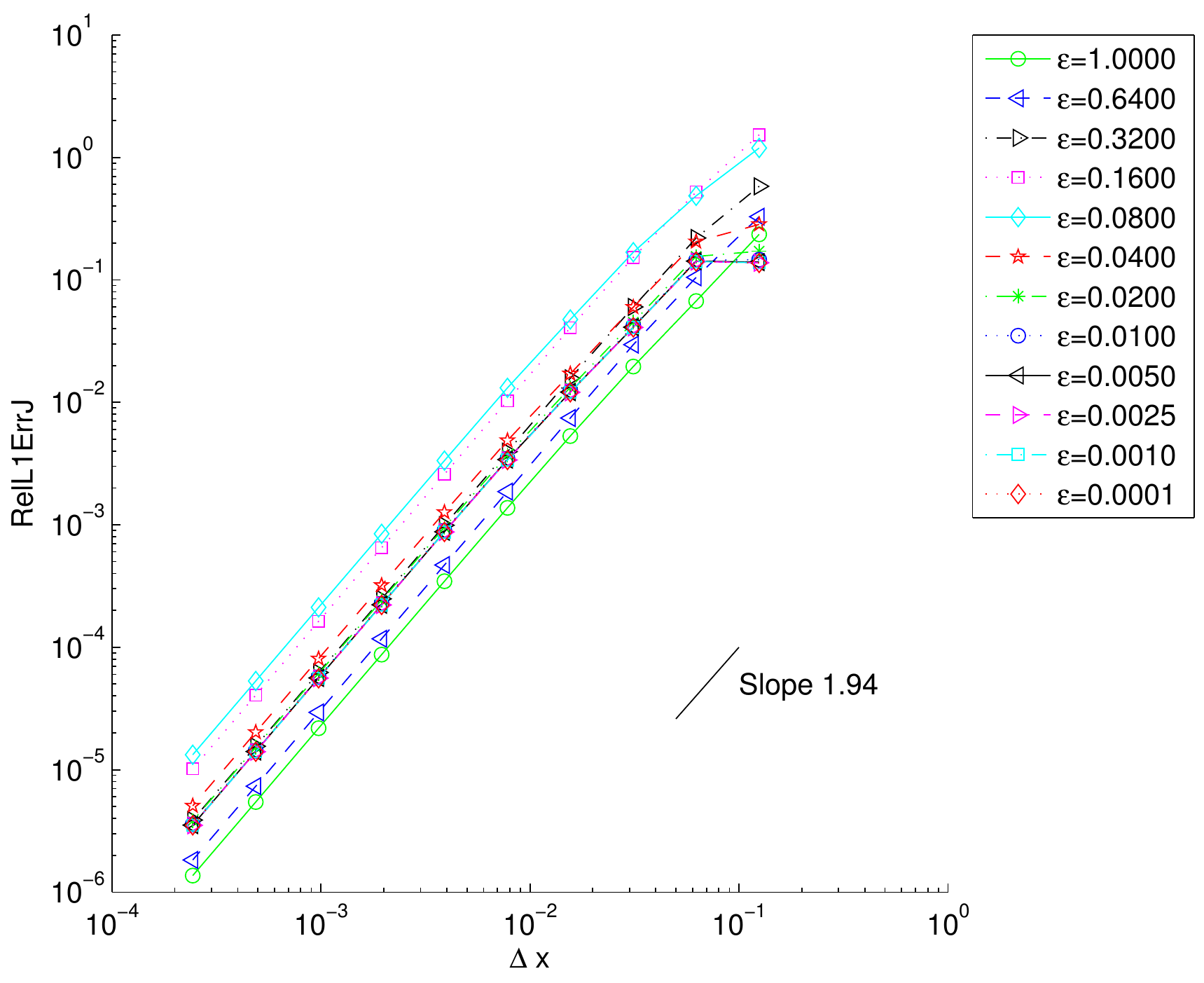}}\qquad
  \subfigure[Error w.r.t $\eps$]{\includegraphics[width=.46\textwidth]{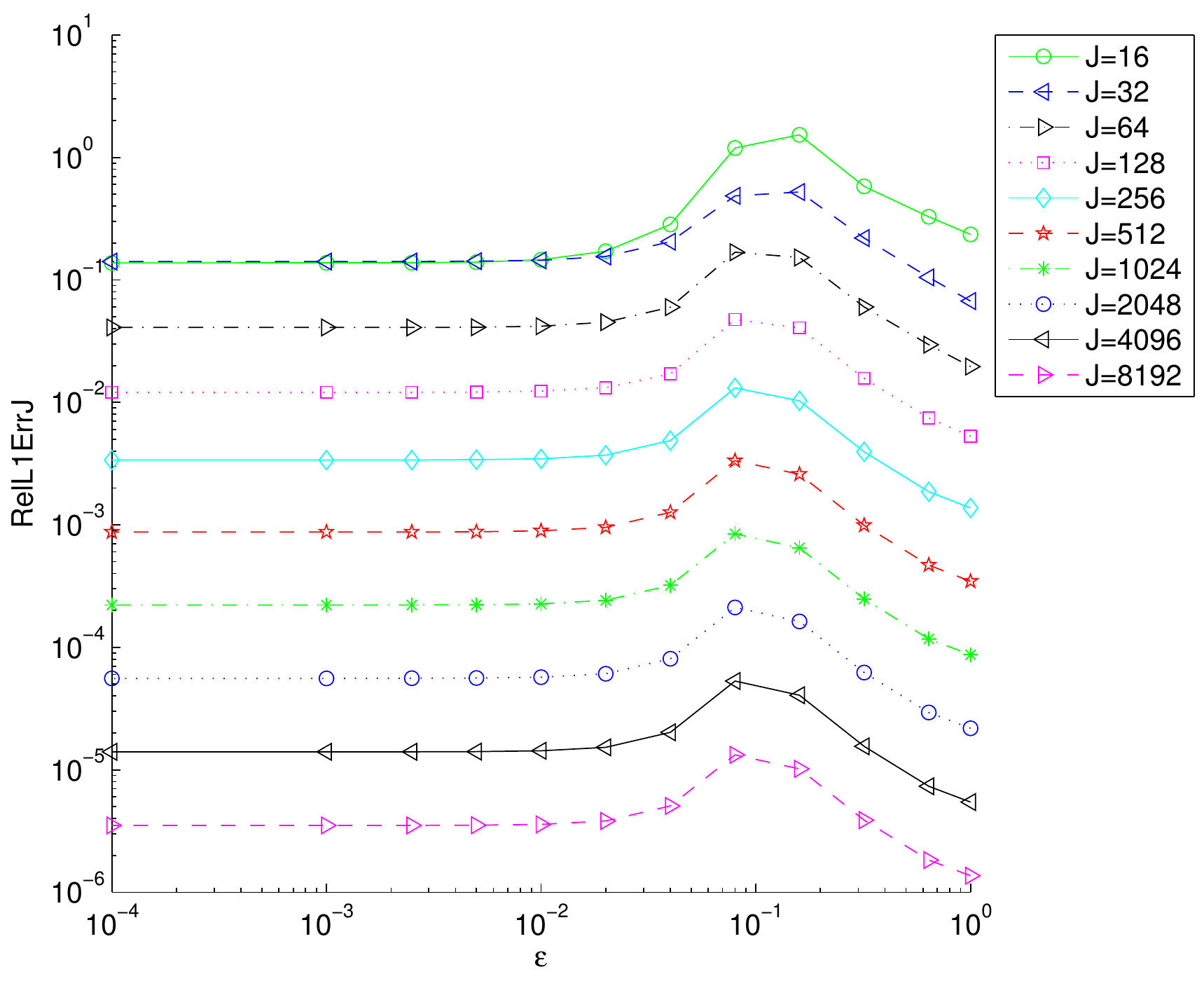}}
  
  \caption{$\mathrm{err}_{\mathbf{j}^\eps}(t=0.05)$ for AP scheme}
\label{fig:err_j_bef_sing_ap}
\end{figure}

\begin{figure}[!htbp]
  \centering
  \subfigure[Error w.r.t $J=2^M$]{\includegraphics[width=.46\textwidth]{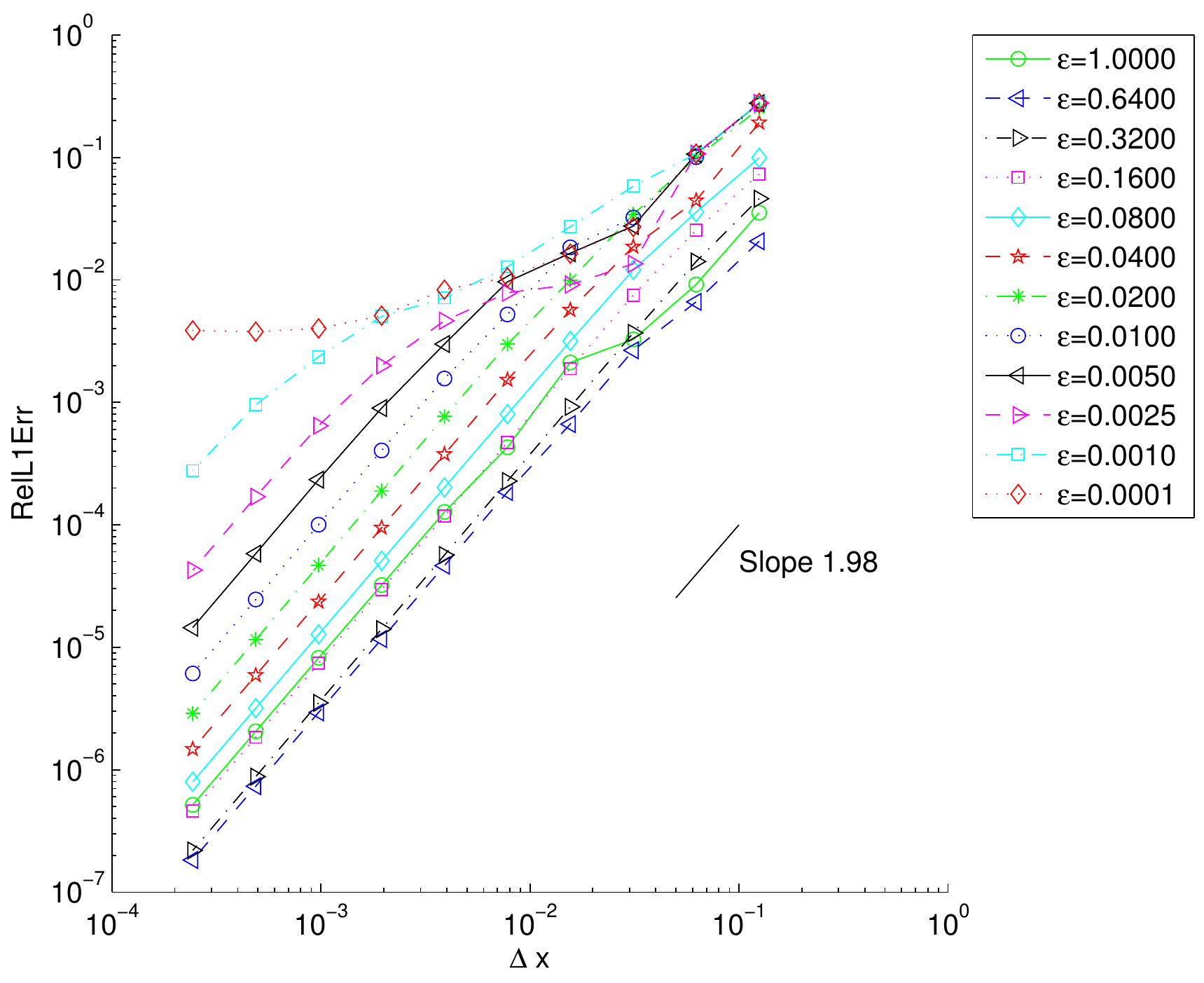}}\qquad
  \subfigure[Error w.r.t $\eps$]{\includegraphics[width=.46\textwidth]{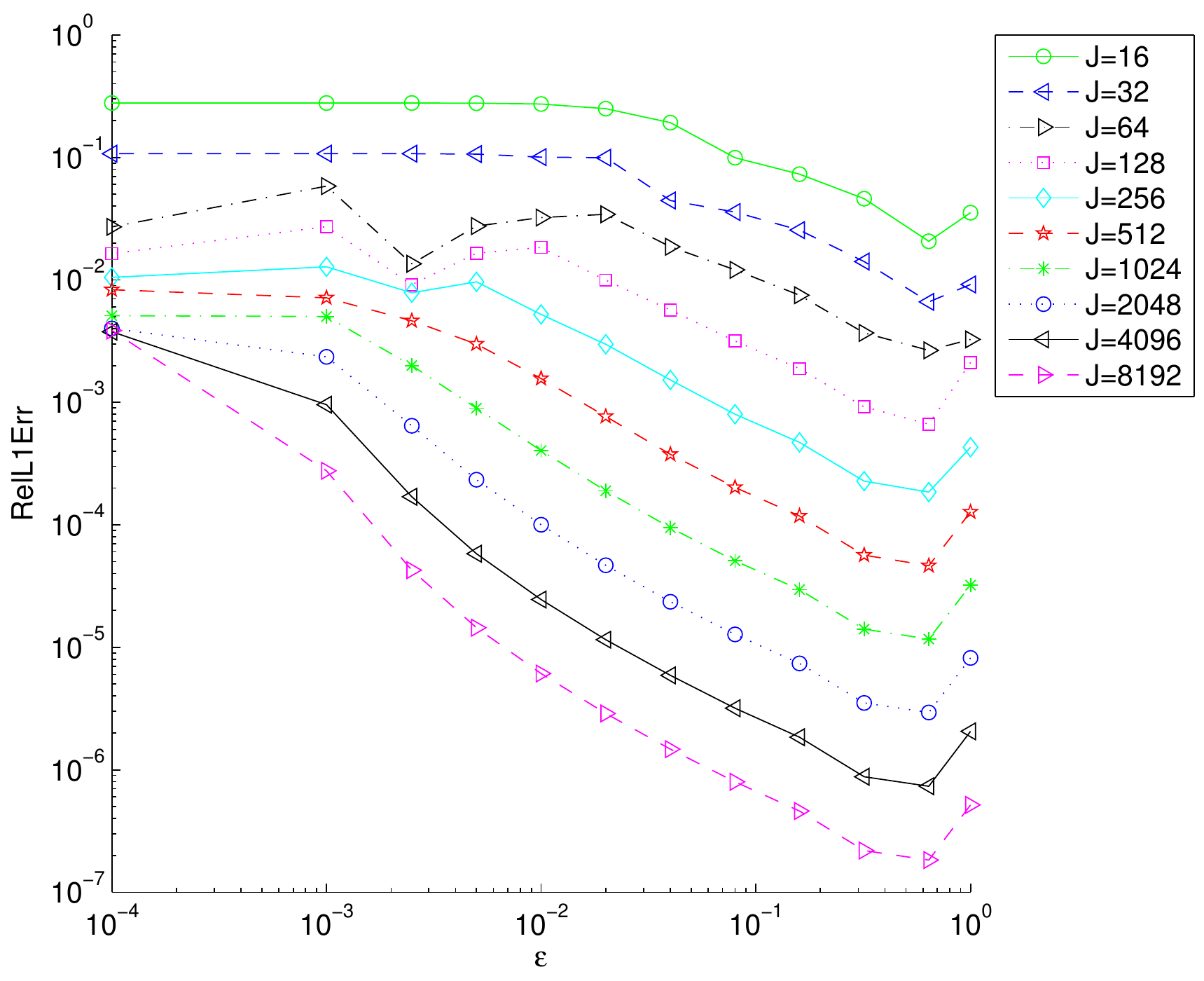}}
  
  \caption{$\mathrm{err}_{\rho^\eps}(t=0.13)$ for AP scheme}
\label{fig:err_rho_aft_sing_ap}  
\end{figure}

\begin{figure}[!htbp]
  \centering
  \subfigure[Error w.r.t $J=2^M$]{\includegraphics[width=.46\textwidth]{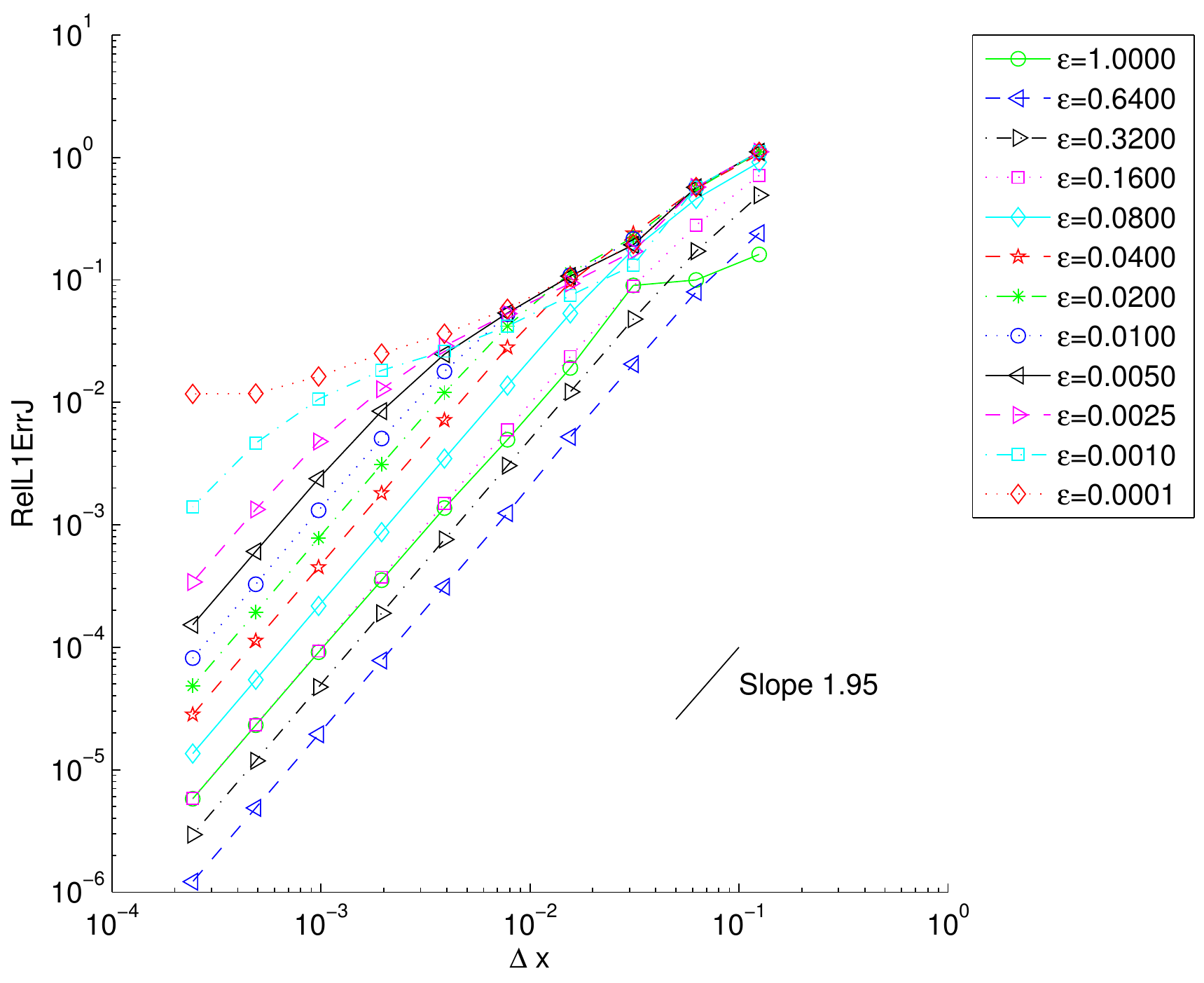}}\qquad
  \subfigure[Error w.r.t $\eps$]{\includegraphics[width=.46\textwidth]{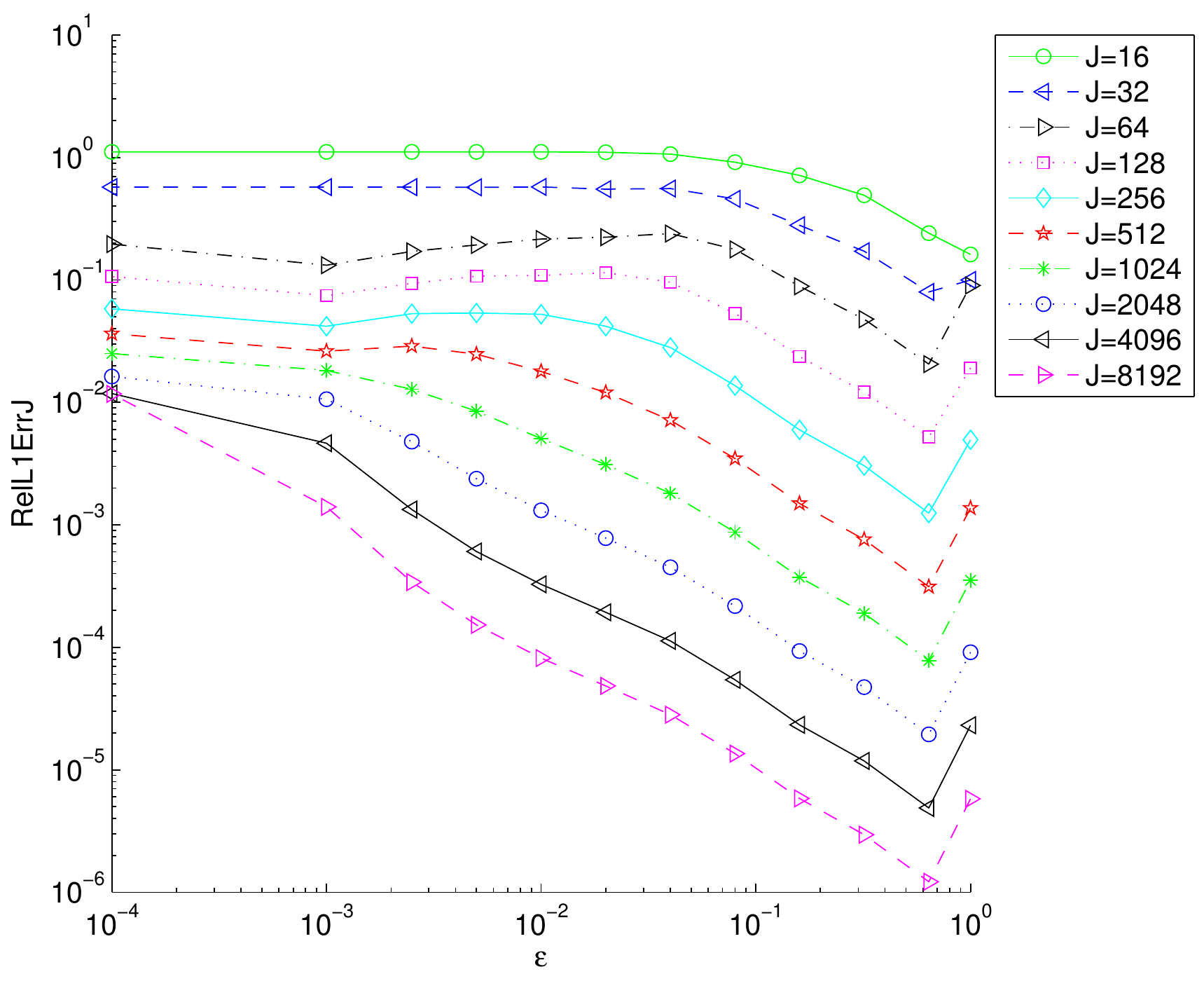}}
  
  \caption{$\mathrm{err}_{\mathbf{j}^\eps}(t=0.13)$ for AP scheme}
\label{fig:err_j_aft_sing_ap}
\end{figure}

In order to show the asymptotic preserving property of our scheme, we
compute relative $\ell^1$ errors at time $t_n$ thanks to the formula 
\begin{align*}
  & \mathrm{err}_{\rho^\eps}(t_n) = \|\rho^\eps_{\rm ref}-\rho^{n,\eps}\|_1/\|\rho^\eps_{\rm ref}\|_1, \quad \rho^{n,\eps}_j=|a_j^{n,\eps}|^2, \\
& \mathrm{err}_{\mathbf{j}^\eps}(t_n)=\|\mathbf{j}^\eps_{\rm ref}-\mathbf{j}^{n,\eps}\|_1/\|\mathbf{j}^\eps_{\rm ref}\|_1, \quad \mathbf{j}^{n,\eps}_j=\eps \IM(\overline{a_j^{n,\eps}} \nabla a^{n,\eps}_j)+\rho^{n,\eps}_j v^{n,\eps}_j,
\end{align*}
where the $\ell^1$ norm is
\begin{equation*}
\|u\|_1=\dx\sum_{k=1}^{J-1}|u_j|.
\end{equation*}
We evaluate the error function for different values of $J=2^M$, where
$M$ is an integer here chosen in $[4,13]$, and various scaled Planck
parameter $\eps$ in $[10^{-4},10^0]$. 
We plot on Figures~\ref{fig:err_rho_bef_sing_ap},
\ref{fig:err_j_bef_sing_ap}, \ref{fig:err_rho_aft_sing_ap} and
\ref{fig:err_j_aft_sing_ap} the relative error on physical observables
$\rho^\eps$ and $\mathbf{j}^\eps$
computed with our AP schemes at time $t=0.05$ and $t=0.13$
respectively before and after the formation of singularities. We
clearly see that the error is proportional to $(\Delta x)^2$,
independently of $\eps$ before formation of singularities. Indeed, the
mean slope of the error curves with respect to $\Delta x$ is close to
the value $2$ (subfigures (a)). The independence with respect to
$\eps$ appears in subfigures (b) where error curves are flat. 
After the formation of shocks, the behavior of our scheme is very different and
the mesh parameters have to be reduced as $\eps \to 0$ to get good accuracy.

\begin{figure}[!htbp]
  \centering
  \subfigure[Error w.r.t $J=2^M$]{\includegraphics[width=.46\textwidth]{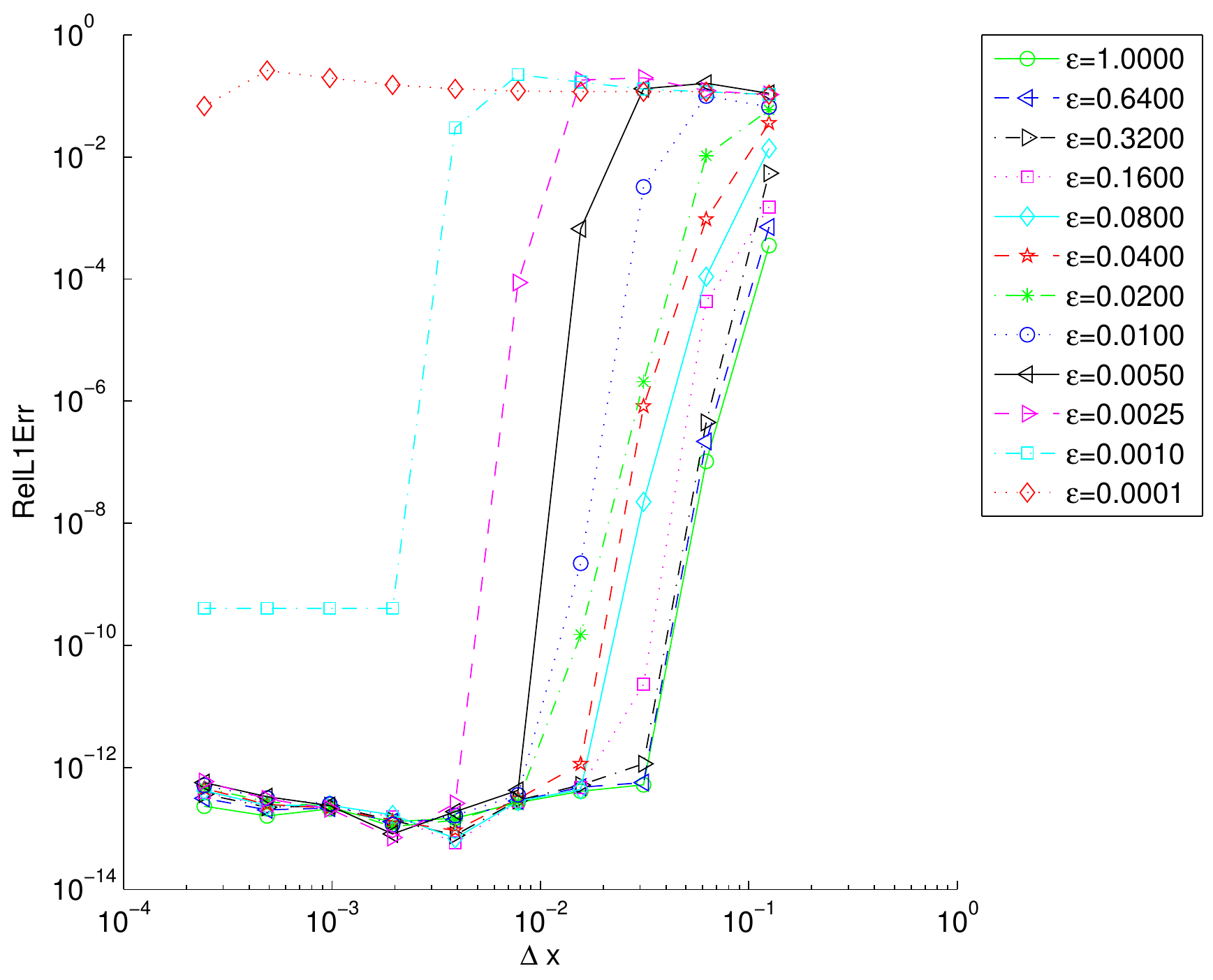}}\qquad
  \subfigure[Error w.r.t $\eps$]{\includegraphics[width=.46\textwidth]{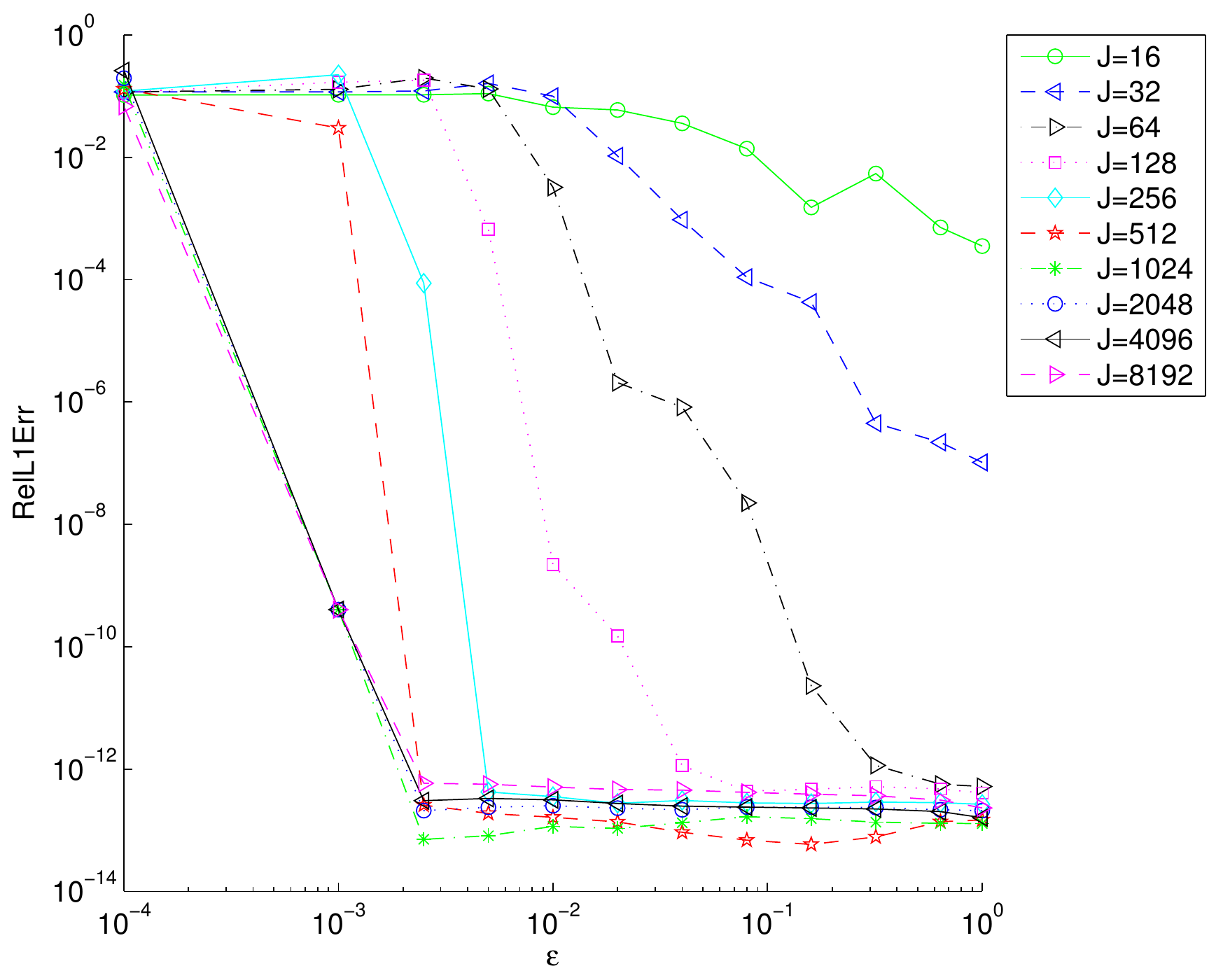}}
  
  \caption{$\mathrm{err}_{\rho^\eps}(t=0.05)$ for time splitting scheme}
\label{fig:err_rho_bef_sing_sp}  
\end{figure}

\begin{figure}[!htbp]
  \centering

  \subfigure[Error w.r.t $J=2^M$]{\includegraphics[width=.46\textwidth]{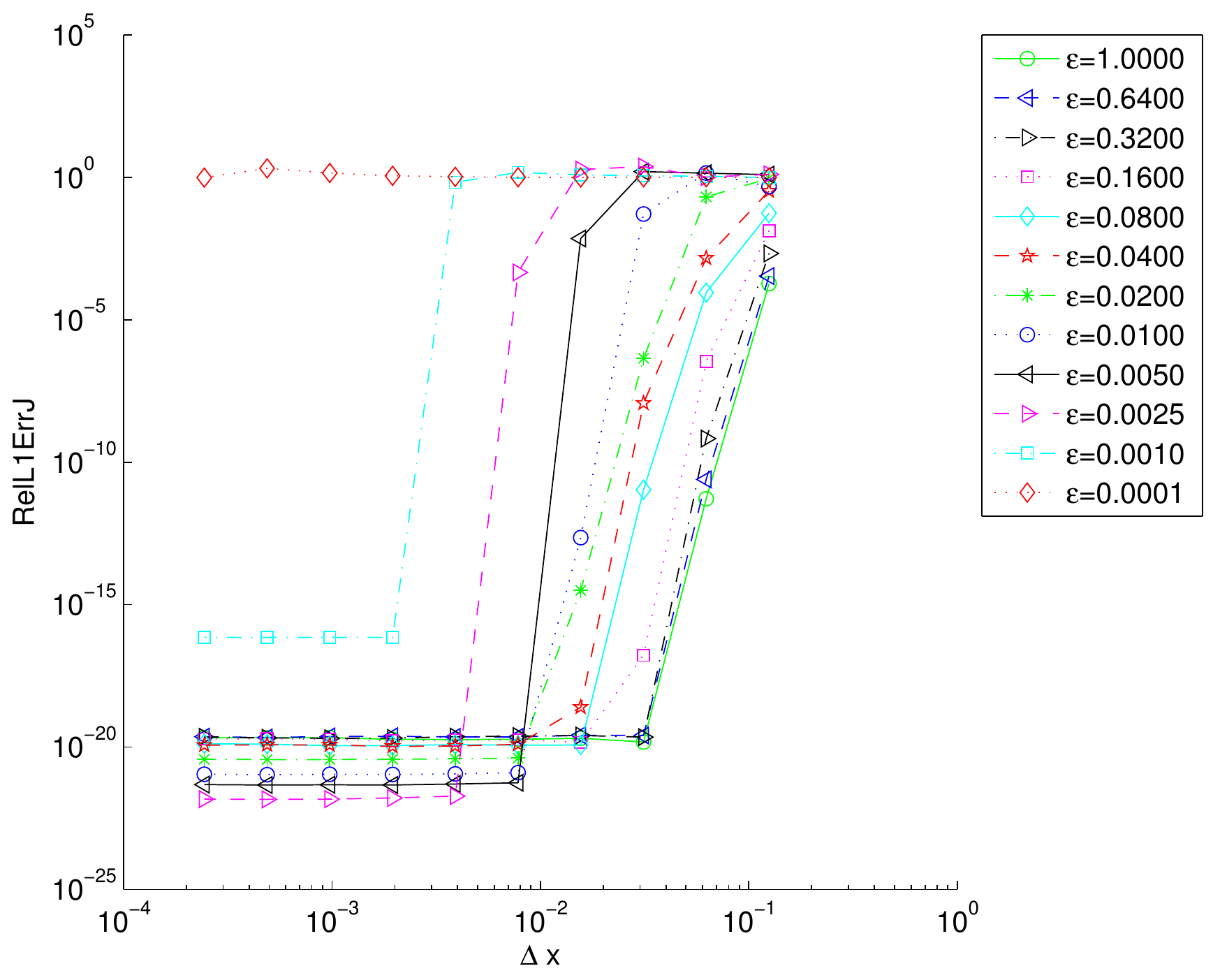}}\qquad
  \subfigure[Error w.r.t $\eps$]{\includegraphics[width=.46\textwidth]{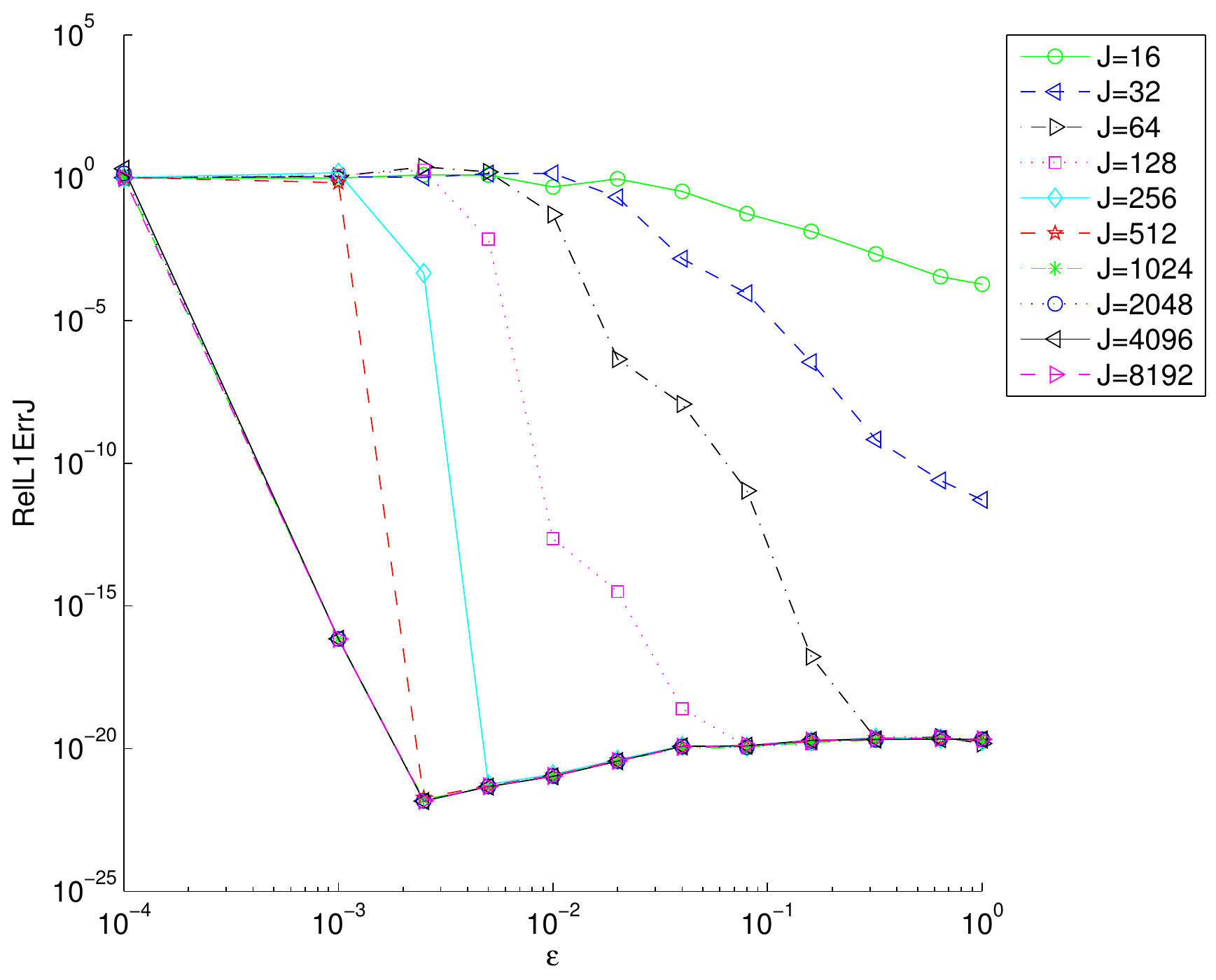}}
  
  \caption{$\mathrm{err}_{\mathbf{j}^\eps}(t=0.05)$ for time splitting scheme}
\label{fig:err_j_bef_sing_sp}  
\end{figure}

\begin{figure}[!htbp]
  \centering
  \subfigure[Error w.r.t $J=2^M$]{\includegraphics[width=.46\textwidth]{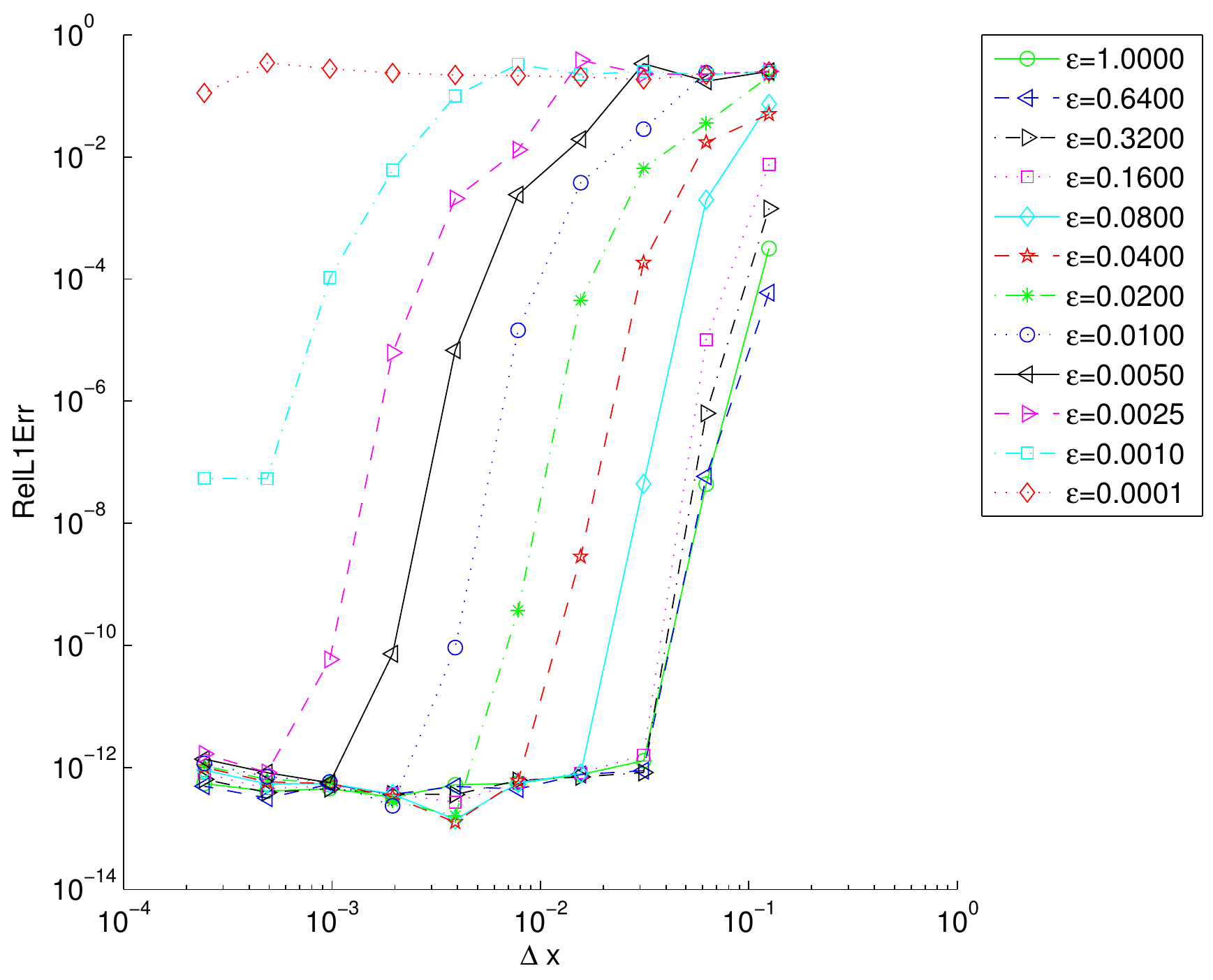}}\qquad
  \subfigure[Error w.r.t $\eps$]{\includegraphics[width=.46\textwidth]{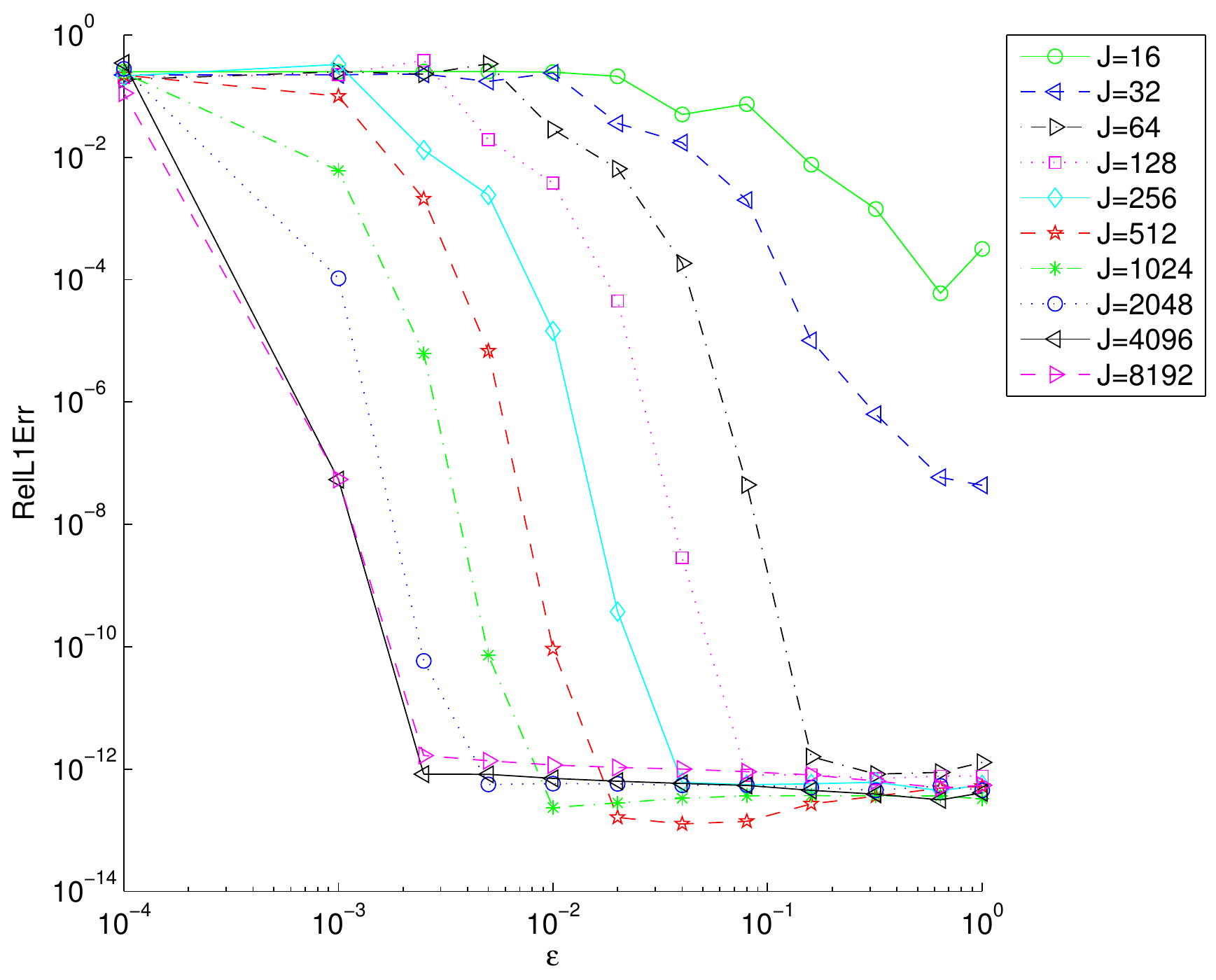}}
  
  \caption{$\mathrm{err}_{\rho^\eps}(t=0.13)$ for time splitting scheme}
\label{fig:err_rho_aft_sing_sp}  
\end{figure}

\begin{figure}[!htbp]
  \centering
  \subfigure[Error w.r.t $J=2^M$]{\includegraphics[width=.46\textwidth]{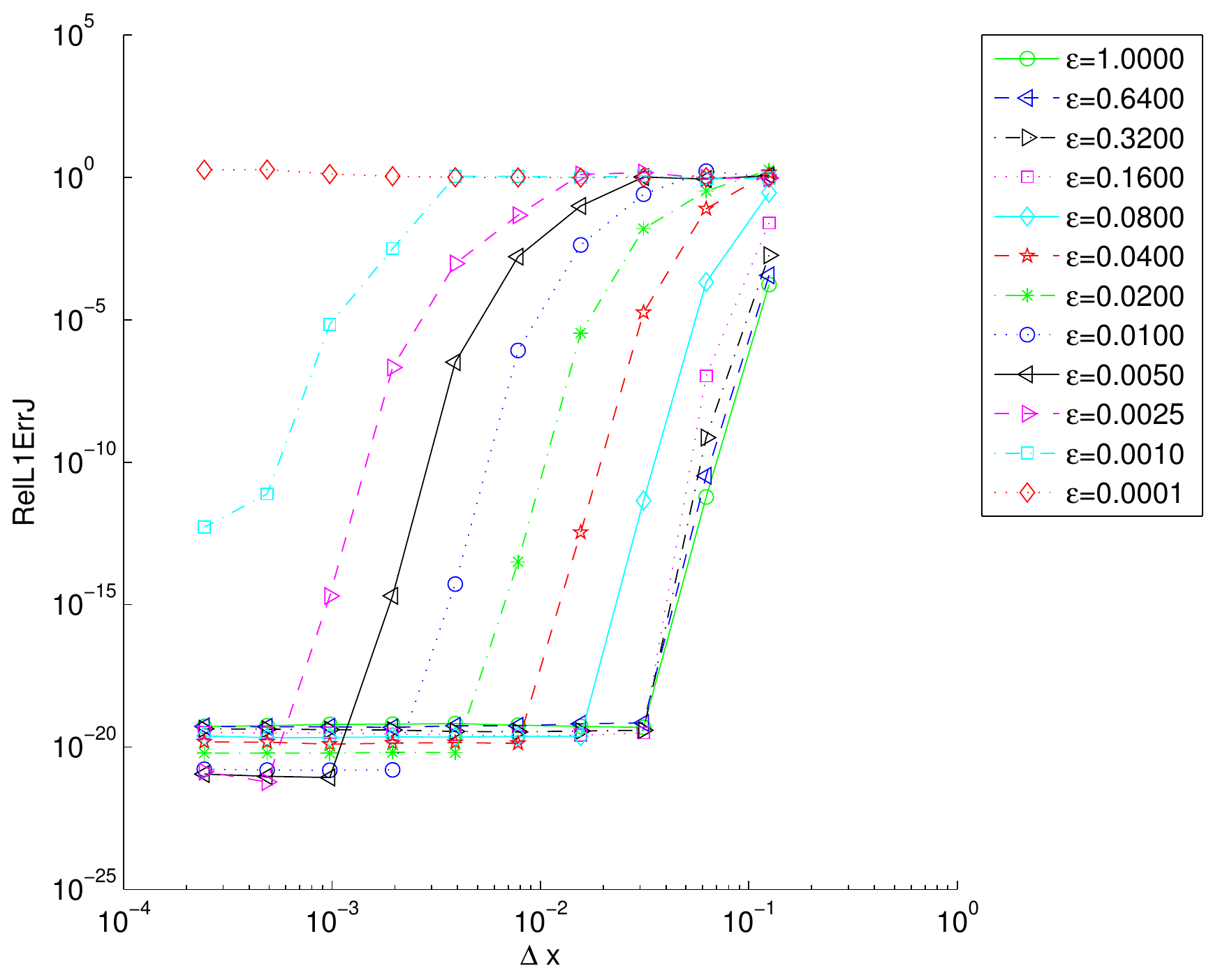}}\qquad
  \subfigure[Error w.r.t $\eps$]{\includegraphics[width=.46\textwidth]{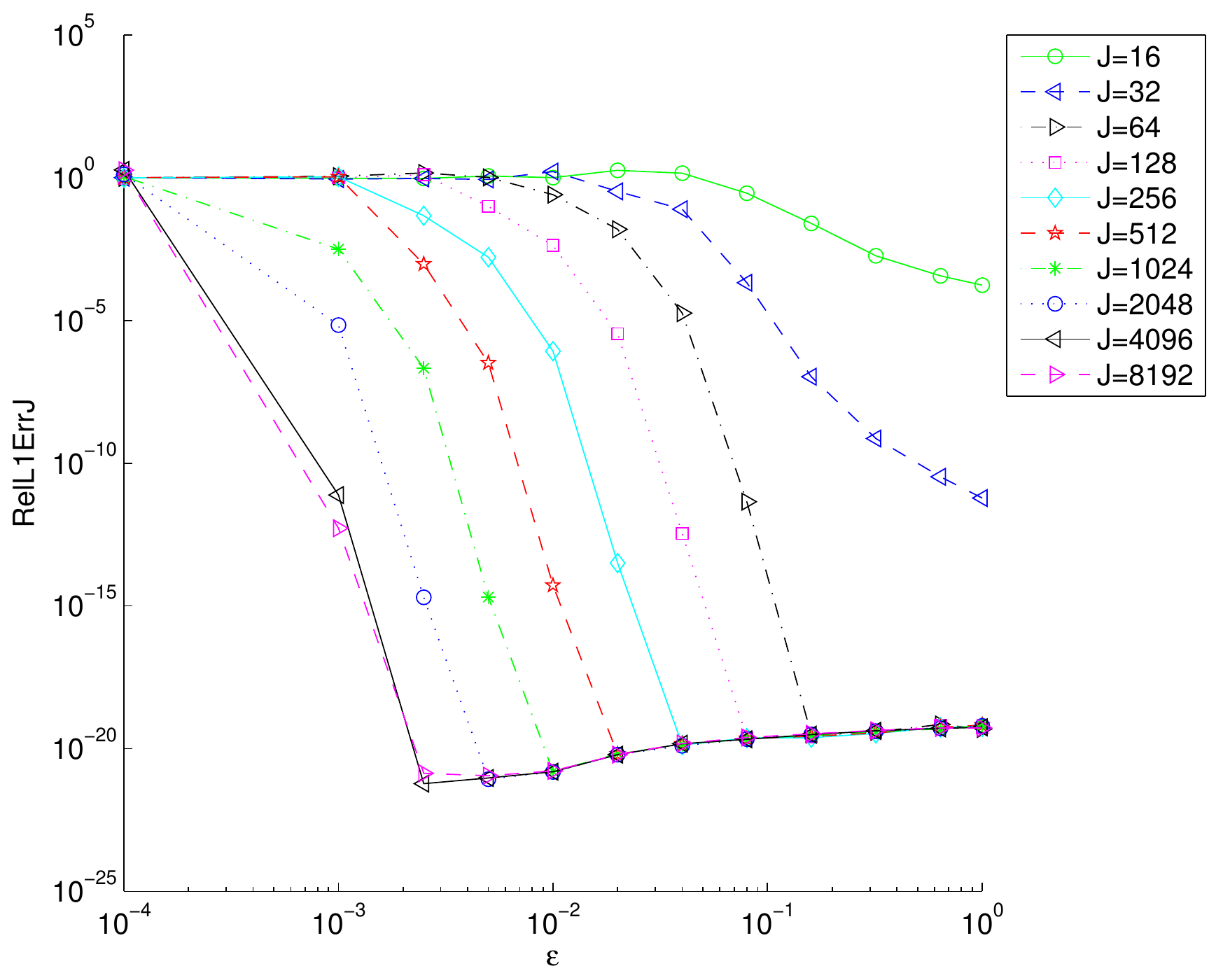}}
  
  \caption{$\mathrm{err}_{\mathbf{j}^\eps}(t=0.13)$ for time splitting scheme}
\label{fig:err_j_aft_sing_sp}  
\end{figure}

For a comparison, we present the same study for the classical time
splitting scheme applied to \eqref{eq:nls}-\eqref{eq:ci} on Figures~\ref{fig:err_rho_bef_sing_sp},
\ref{fig:err_j_bef_sing_sp}, \ref{fig:err_rho_aft_sing_sp},
\ref{fig:err_j_aft_sing_sp}. Thanks to the spectral accuracy of Fast
Fourier Transform, the error levels are smaller than
for our AP scheme which is only second order with respect to time and
space variables. However, contrary to our AP scheme, the error curves
clearly depend on $\eps$, the error being of order $\O(1)$ when $\eps$
is smaller than $2.10^{-3}$. In order to have good accuracy, it is
required to take $\Delta x \lesssim \eps$.

\begin{figure}[!htbp]
  \centering
  \subfigure[$J=128$]{\includegraphics[width=.46\textwidth]{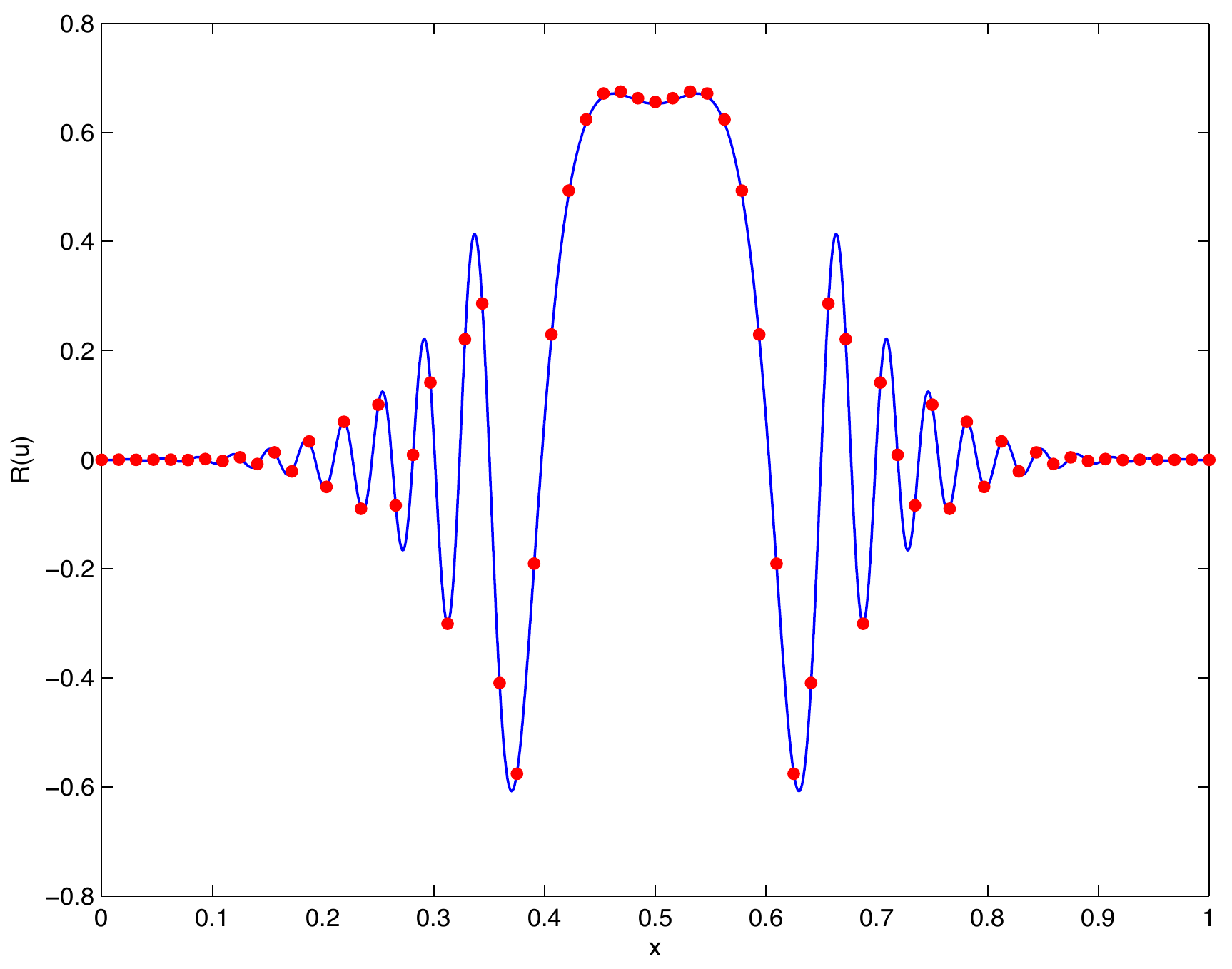}} \qquad
  \subfigure[$J=256$]{\includegraphics[width=.46\textwidth]{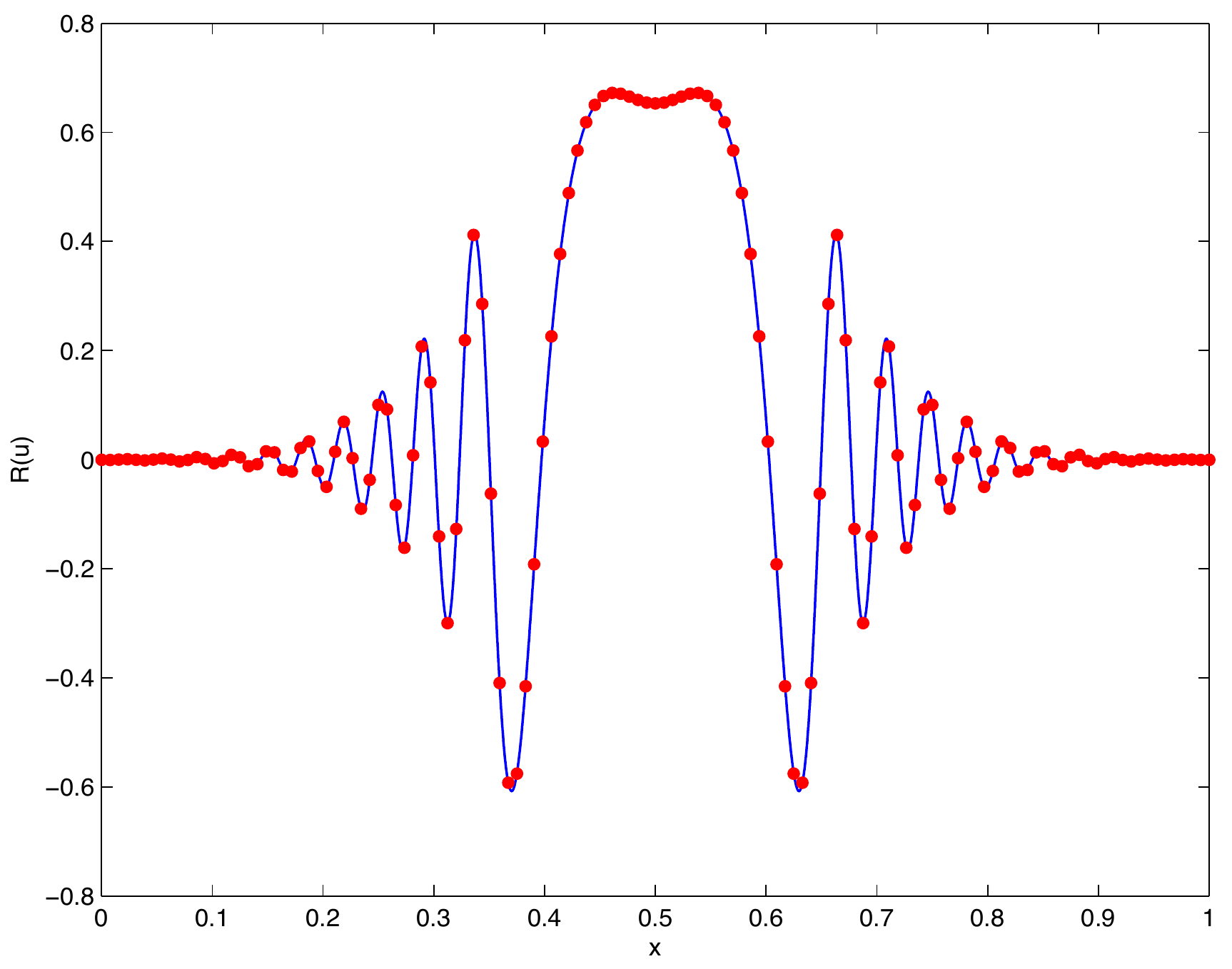}}
  \caption{$\text{Re}(u^\eps)$ at $T=0.05$ for $\eps=0.005$}
  \label{fig:recons_bef_sing}
\end{figure}

\begin{figure}[!htbp]
  \centering
  \subfigure[$J=128$]{\includegraphics[width=.46\textwidth]{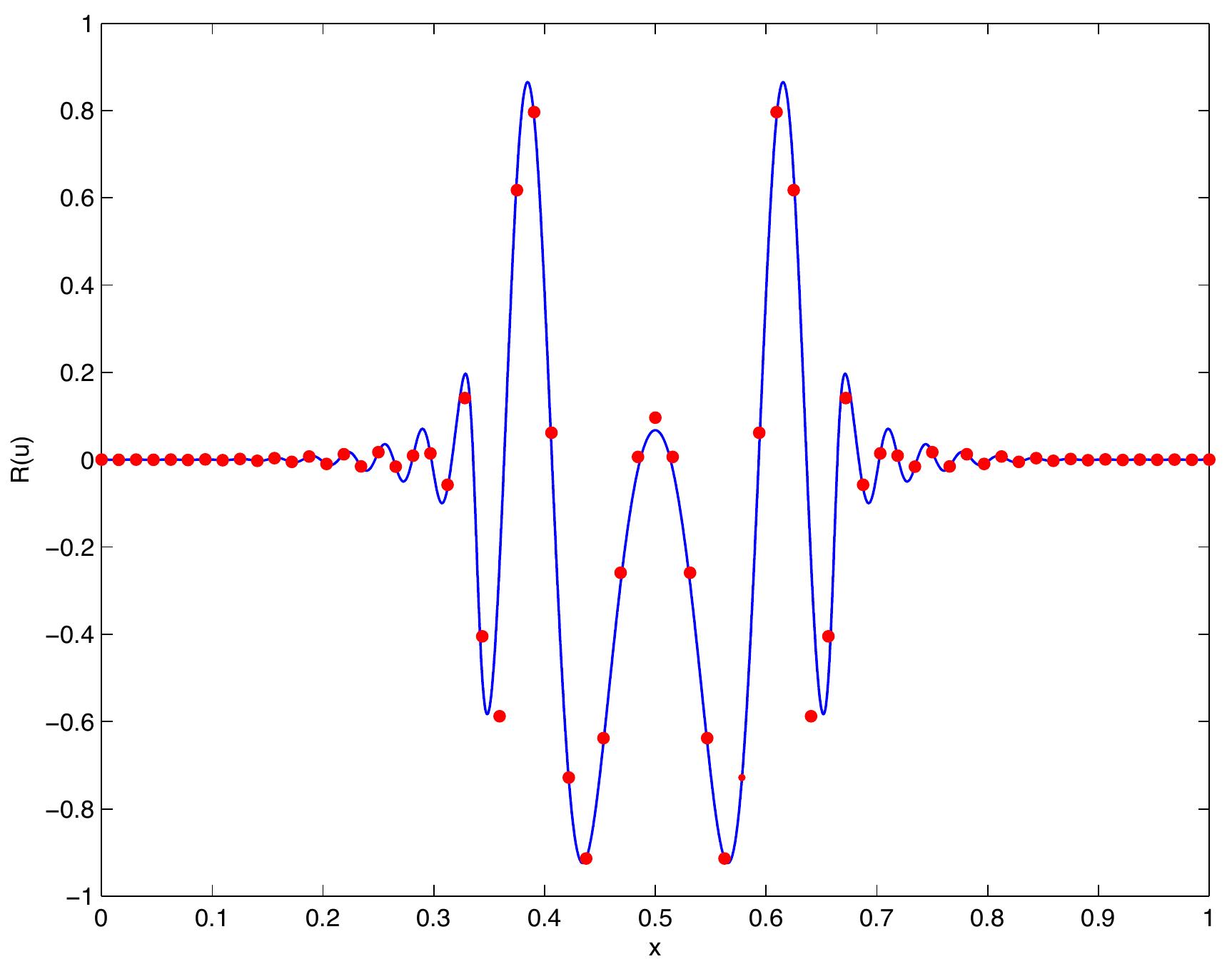}} \qquad
  \subfigure[$J=256$]{\includegraphics[width=.46\textwidth]{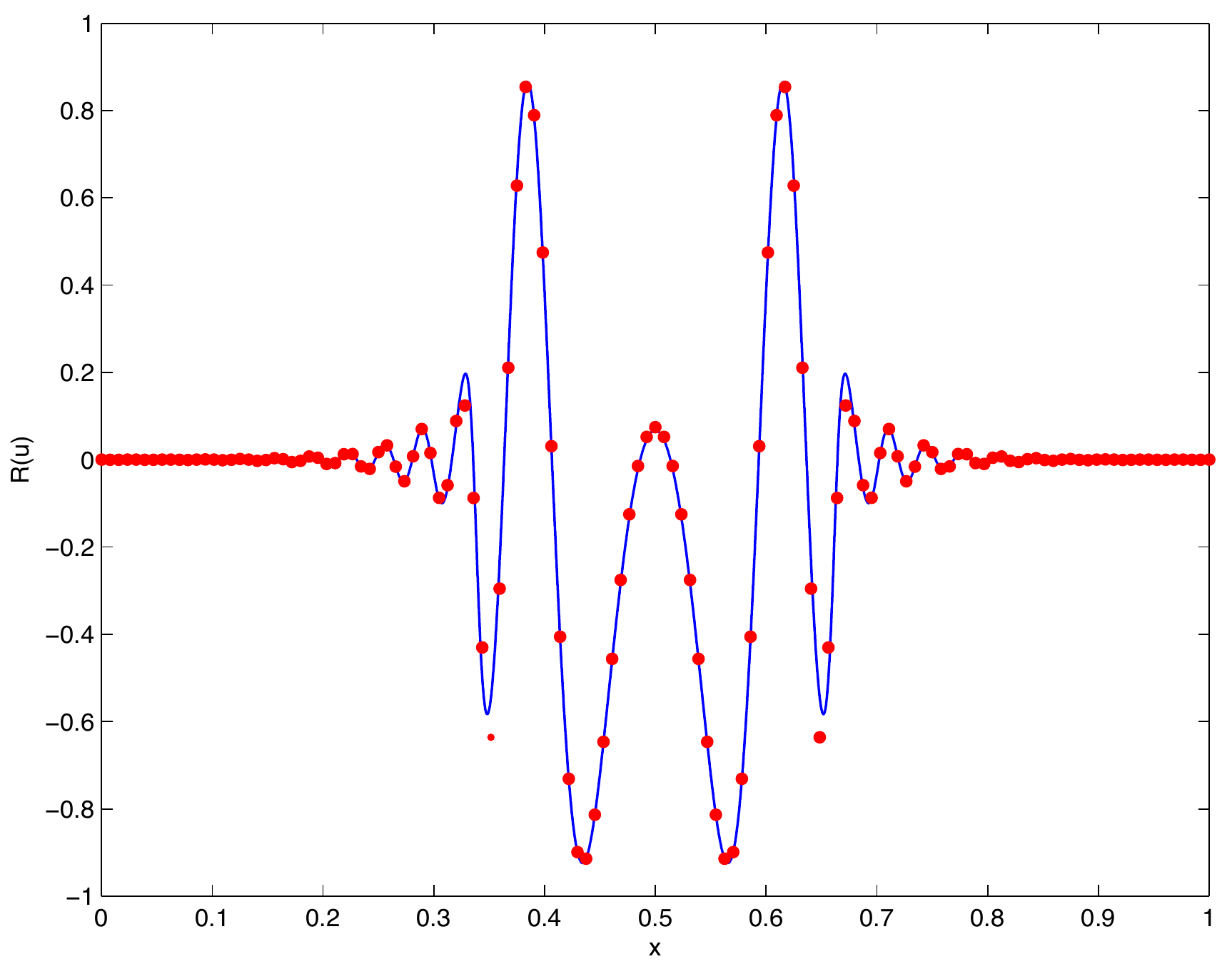}}
  \caption{$\text{Re}(u^\eps)$ at $T=0.13$ for $\eps=0.005$}
  \label{fig:recons_aft_sing}
\end{figure}

\begin{figure}[!htbp]
  \centering
  \subfigure[$J=512$]{\includegraphics[width=.46\textwidth]{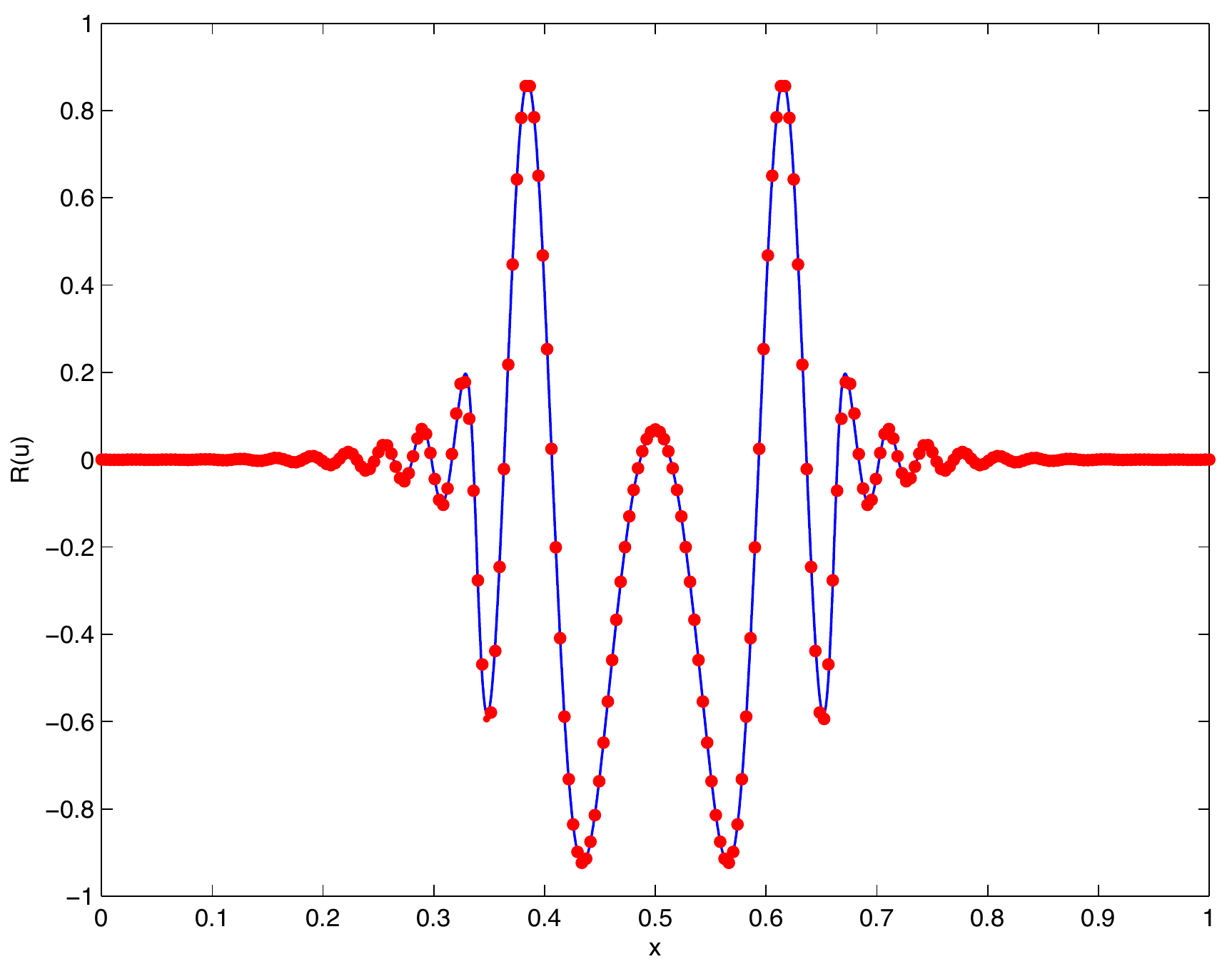}}
  \caption{$\text{Re}(u^\eps)$ at $T=0.13$ for $\eps=0.005$}
  \label{fig:recons_aft_sing_512}
\end{figure}

To make our presentation complete, we present the reconstruction of $u^\eps$ thanks to
equation \eqref{eq:phi} with $\eps=0.005$. We therefore compute the phase $\phi^\eps$
with amplitude $a^\eps$ and velocity $v^\eps$. Since we have access to
values of these quantities at every discrete time $t_n$, we approximate the time
dependent integral thank to a simple rectangular quadrature. We see
that we can recover a very good approximation of the wave function
$u^\eps$ with very few points before the formation of singularities
(see Fig. \ref{fig:recons_bef_sing}). Obviously, one needs more points
to have a good reconstruction of the wave function passed the
singularities (see Figs. \ref{fig:recons_aft_sing} and
\ref{fig:recons_aft_sing_512}).

\clearpage

\subsection{Numerical experiments in dimension 2}
\begin{figure}[!htbp]
    \centering
    
    \subfigure[$\rho_{\mathrm{ref}}$]{\includegraphics[width=.46\textwidth]{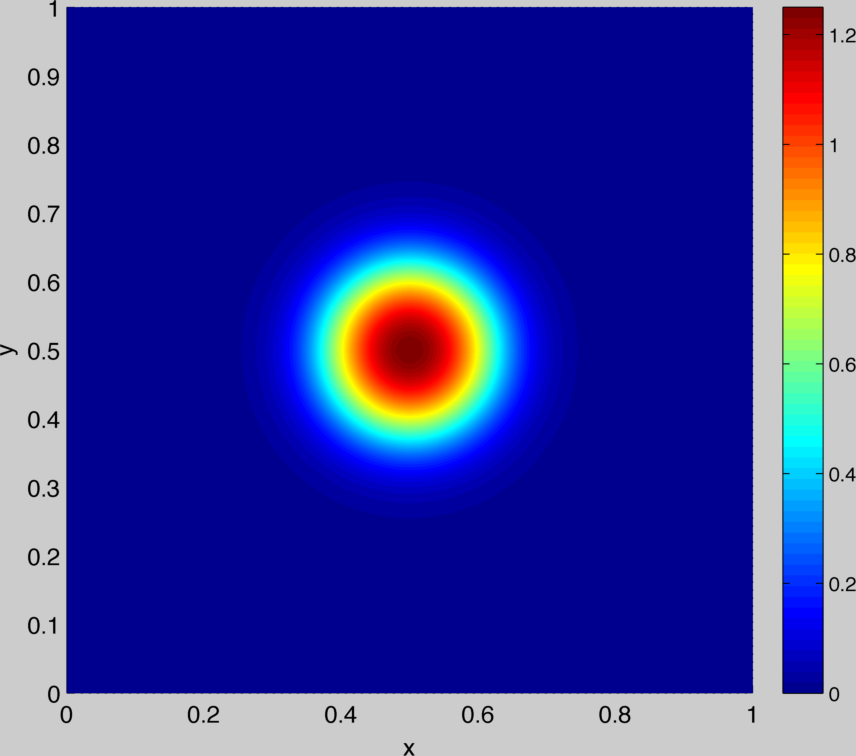}}\qquad
    \subfigure[$\mathbf{j}_{\mathrm{ref}}$]{\includegraphics[width=.46\textwidth]{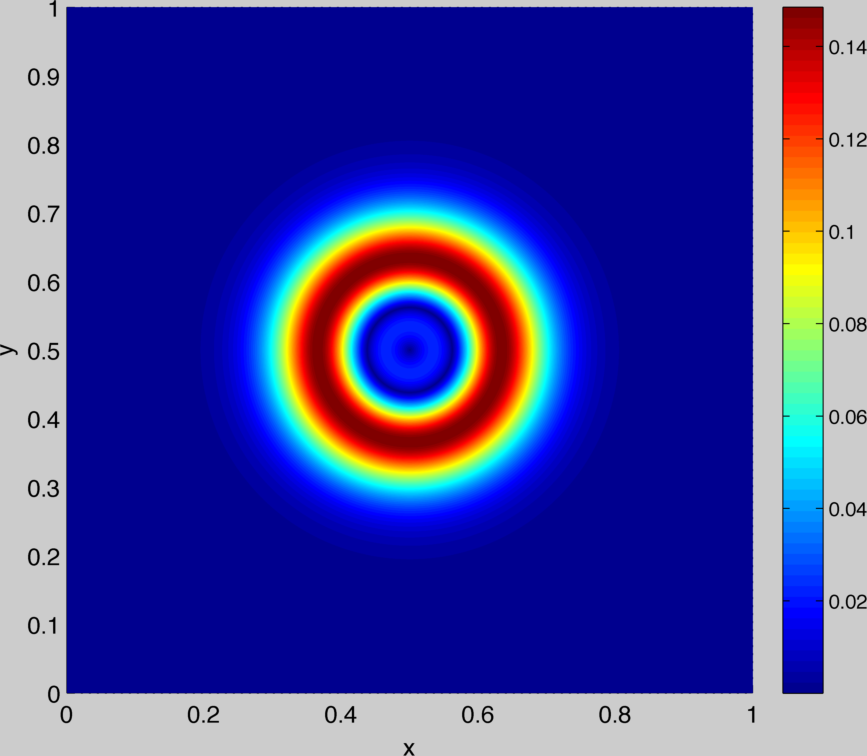}}
    \caption{Contour plot of particle and current densities at time $t=0.05$ and for $\eps=5.10^{-4}$}
    \label{fig:ref2D_bef_sing}  

  \end{figure}

The previous subsection was devoted to one dimensional simulations. We
present here the equivalent error analysis for the two dimensional
case. The computational domain is now 
the square $[-0.5,1.5]^2$ discretized with $J$ points in each
direction. We still need a reference solution generated with Strang
splitting scheme applied to equation \eqref{eq:nls}-\eqref{eq:ci}. We
take $J_{\text{ref}}=2^{13}=8192$, so $\Delta x \sim 2.5\, 10^{-4}$
and $\Delta t =\eps/100$. In a first test, the initial datum is related to the one
chosen for the one dimensional case
$$
\begin{array}{l}
\displaystyle a_0(x,y)=e^{-25\, r^2},\\[2mm]
\displaystyle v_0(x)=-\frac{1}{5}\nabla \ln{(e^{5\, r}+e^{-5\, r})}.
\end{array}
$$
where $r=\sqrt{(x-0.5)^2+(y-0.5)^2}$. We evaluate the error function for different values of $J=2^M$, where
$M$ is chosen in $[4,11]$, and various scaled Planck
parameter $\eps$ in $[5.10^{-4},10^0]$. The formation time of
singularities is as in one dimension situation located around
$0.10$. The reference particle and current densities for $\eps=5.10^{-4}$ are
represented on Figures~\ref{fig:ref2D_bef_sing} and
\ref{fig:ref2D_aft_sing} before and after formation of
singularities. We clearly see that a tiny front appears after shocks
both in particle and current densities. The wave function presents
strong oscillations with new ones created past the shock creation (see
Fig.~\ref{fig:ref2D_aft_sing})

  \begin{figure}[!htbp]
    \centering
    
    \subfigure[$\rho_{\mathrm{ref}}$]{\includegraphics[width=.46\textwidth]{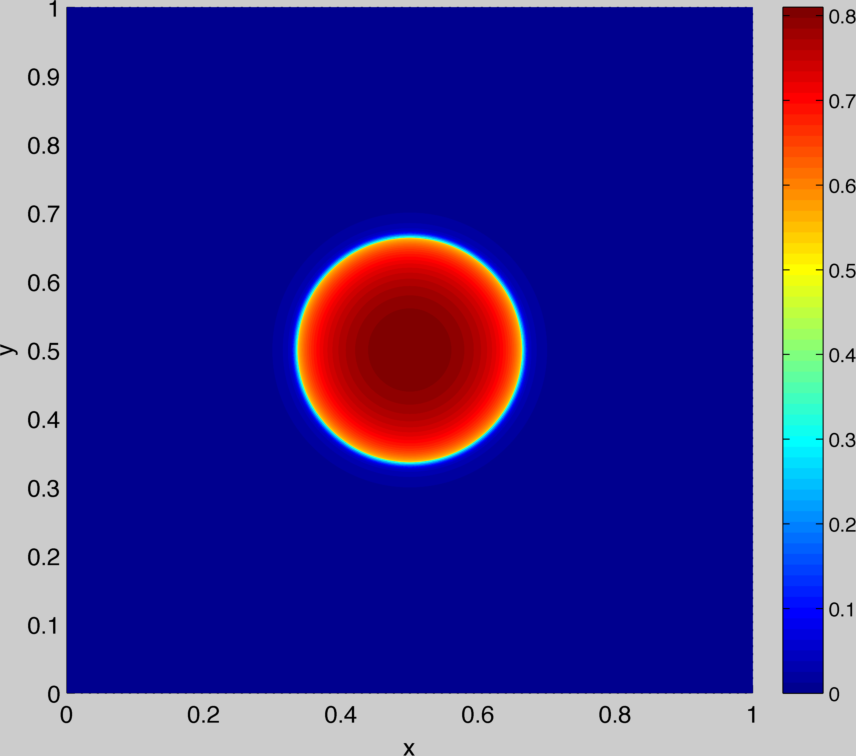}}\qquad
    \subfigure[$\mathbf{j}_{\mathrm{ref}}$]{\includegraphics[width=.46\textwidth]{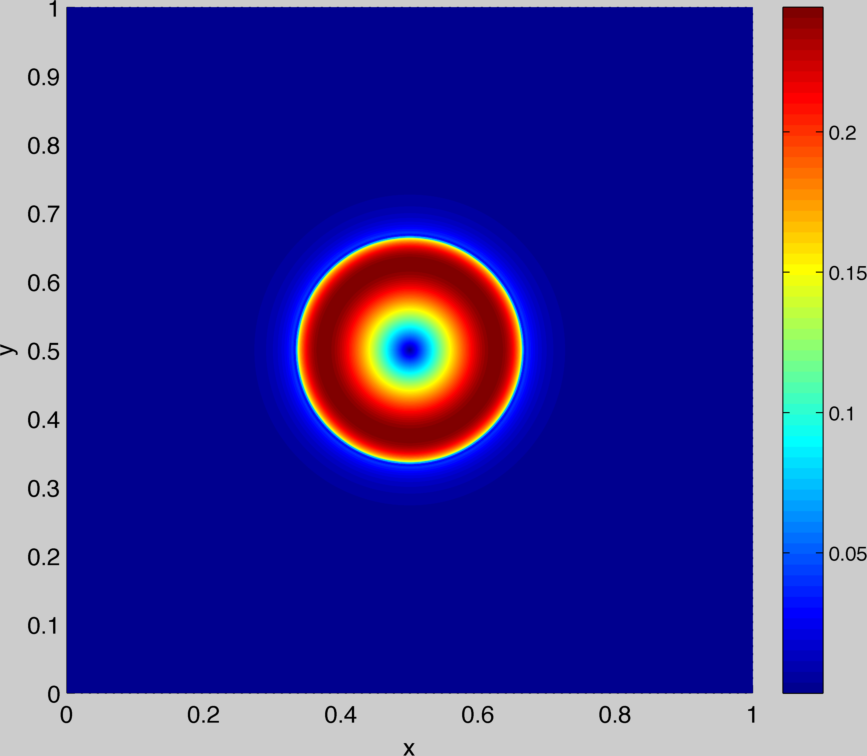}}
    \caption{Contour plot of particle and current densities at time $t=0.13$ and for $\eps=5.10^{-4}$}
    \label{fig:ref2D_aft_sing}  

  \end{figure}

  \begin{figure}[!htbp]
    \centering
    
    \subfigure[$u^\eps$ at time $t=0.05$]{\includegraphics[width=.46\textwidth]{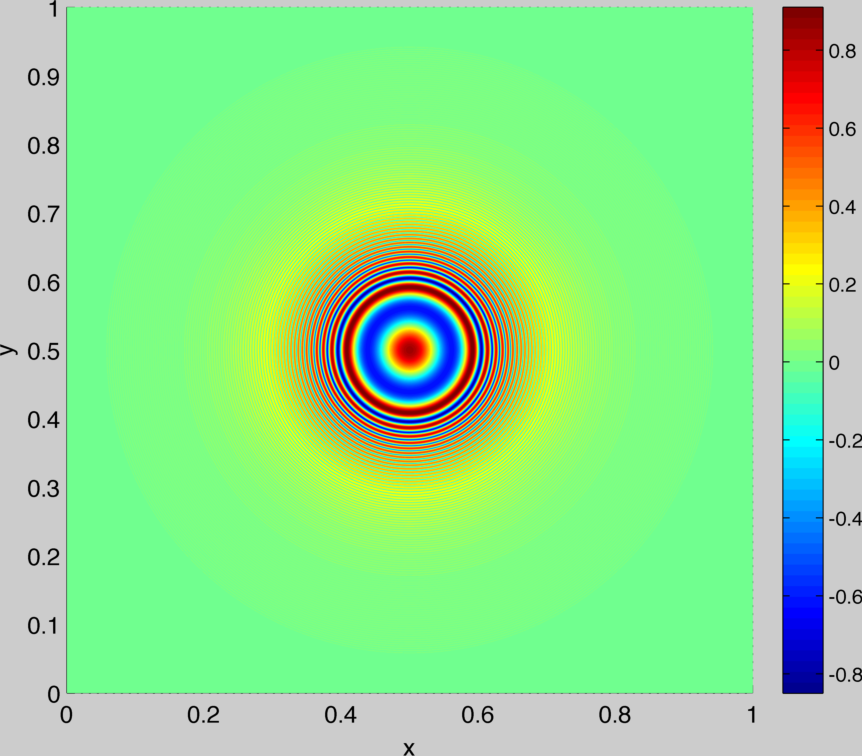}}\qquad
    \subfigure[$u^\eps$ at time $t=0.13$]{\includegraphics[width=.46\textwidth]{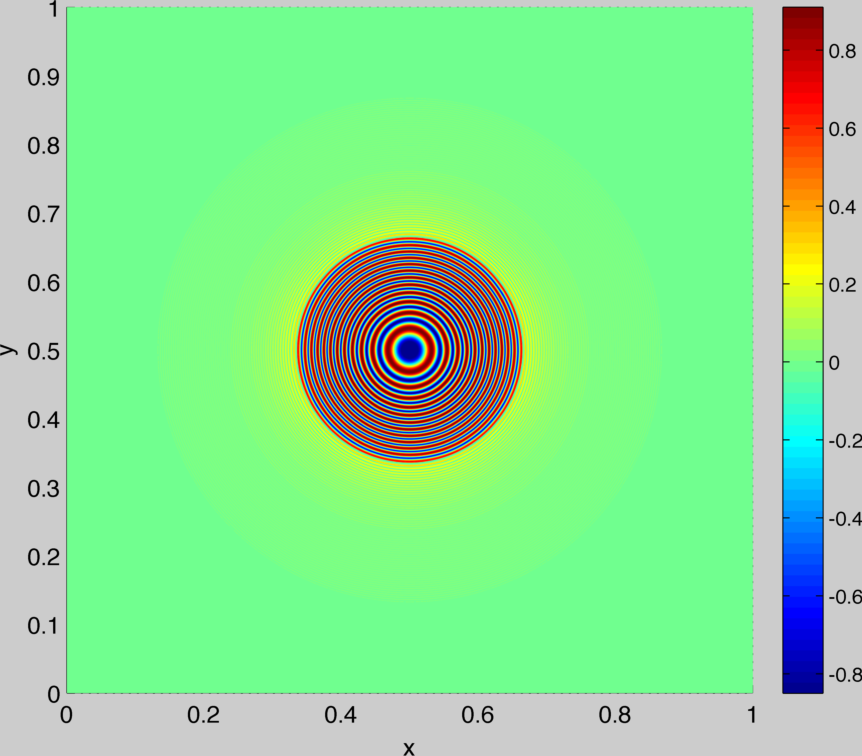}}
    \caption{Contour plot of real part of $u^\eps$ at times $t=0.05$
      and $t=0.13$ and for $\eps=5.10^{-4}$}
    \label{fig:ref2D_aft_sing}  

  \end{figure}
  
\begin{figure}[!htbp]
  \centering

  \subfigure[Error w.r.t $J=2^M$]{\includegraphics[width=.46\textwidth]{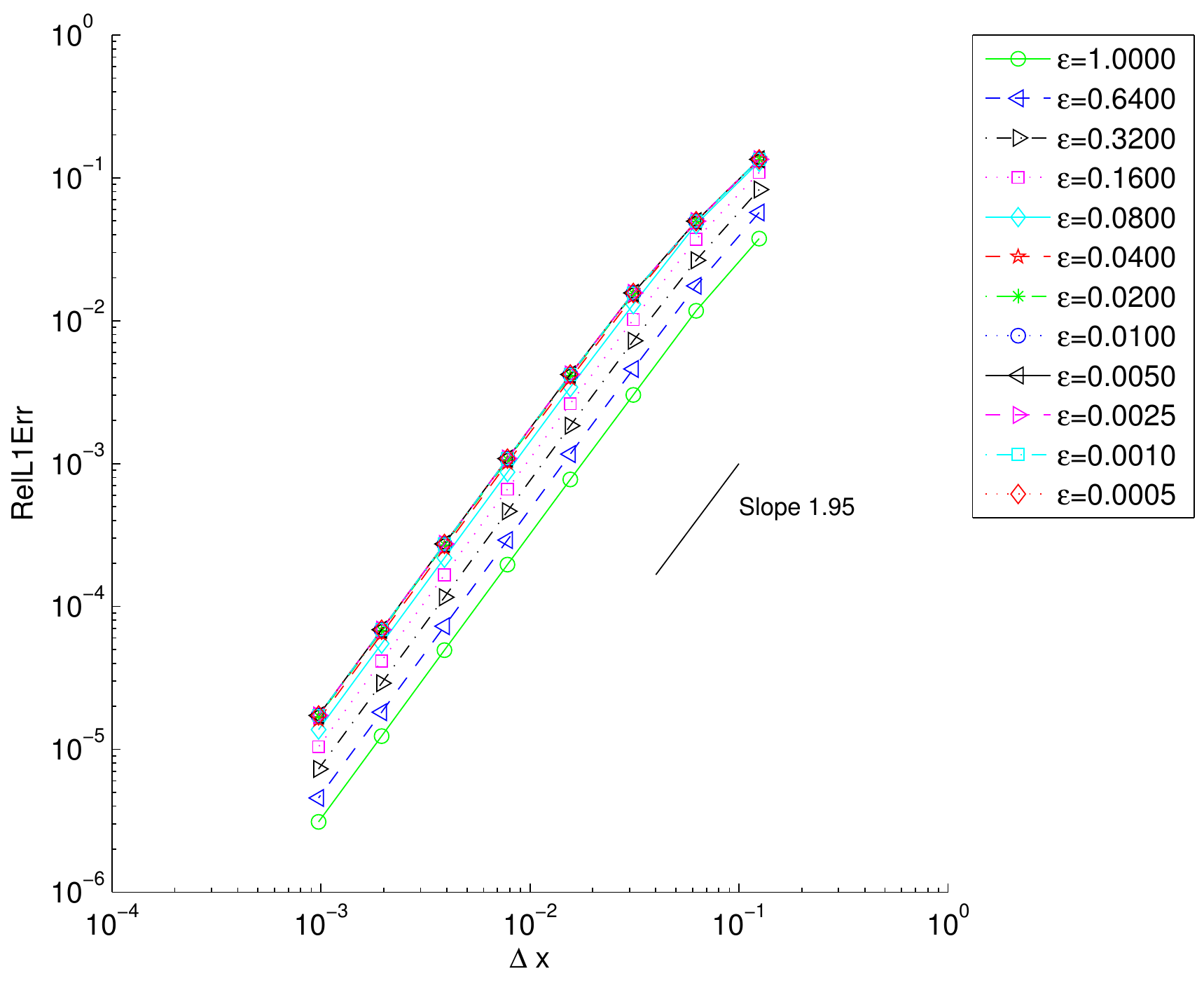}}\qquad
  \subfigure[Error w.r.t $\eps$]{\includegraphics[width=.46\textwidth]{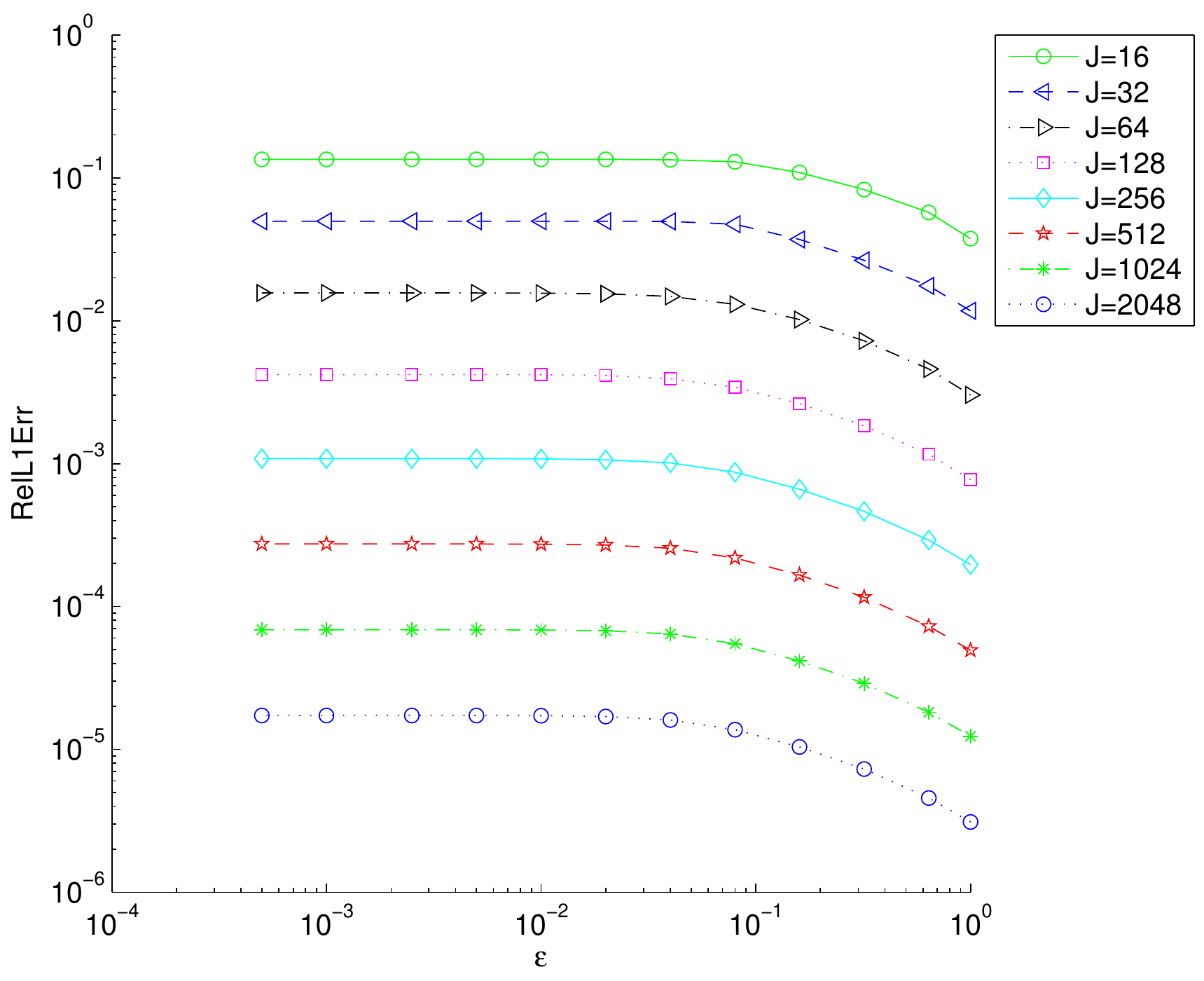}}

  \caption{$\mathrm{err}_{\rho^\eps}(t=0.05)$ for AP scheme}
\label{fig:err2D_rho_bef_sing_ap}  
\end{figure}

\begin{figure}[!htbp]
  \centering

  \subfigure[Error w.r.t $J=2^M$]{\includegraphics[width=.46\textwidth]{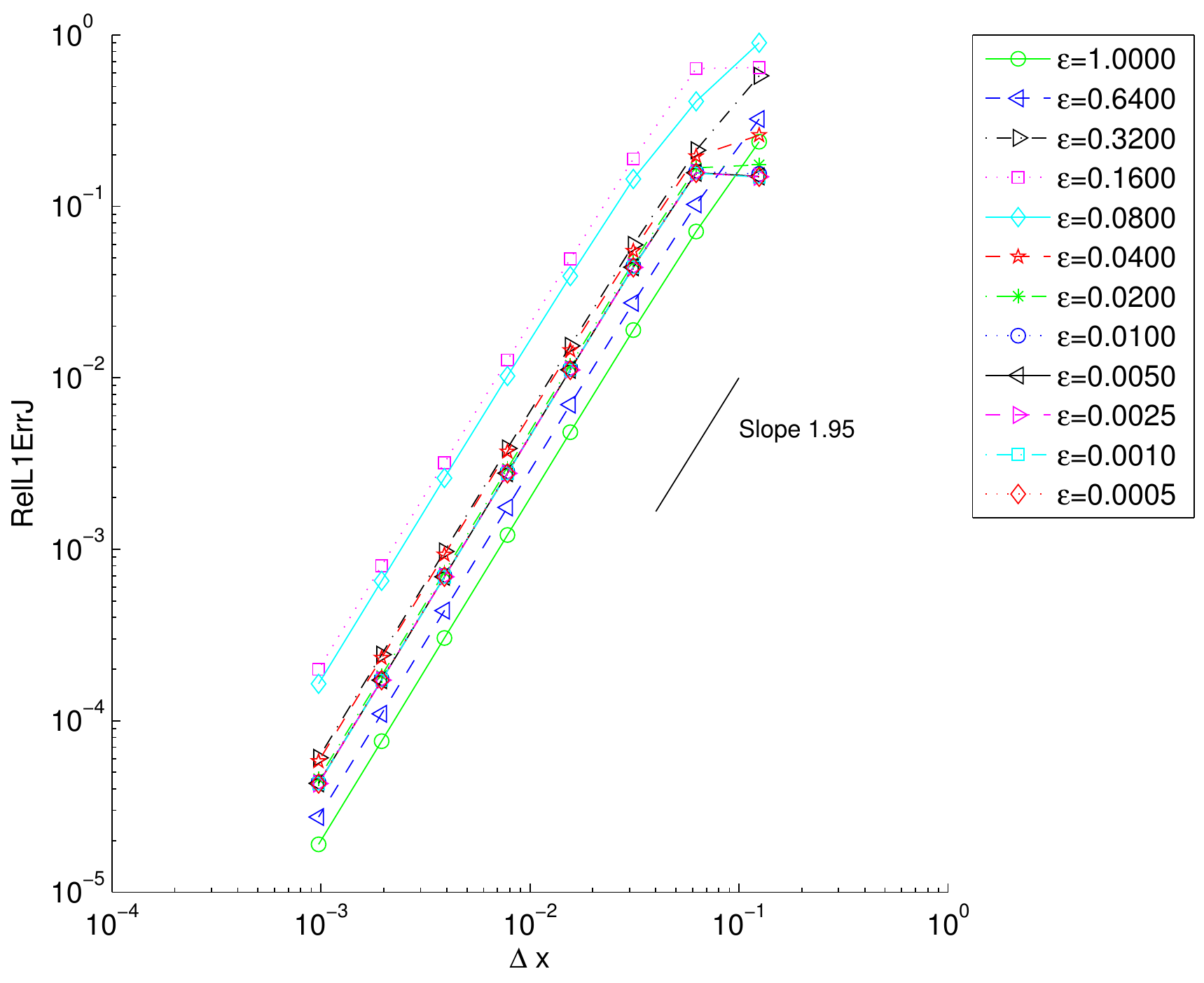}}\qquad
  \subfigure[Error w.r.t $\eps$]{\includegraphics[width=.46\textwidth]{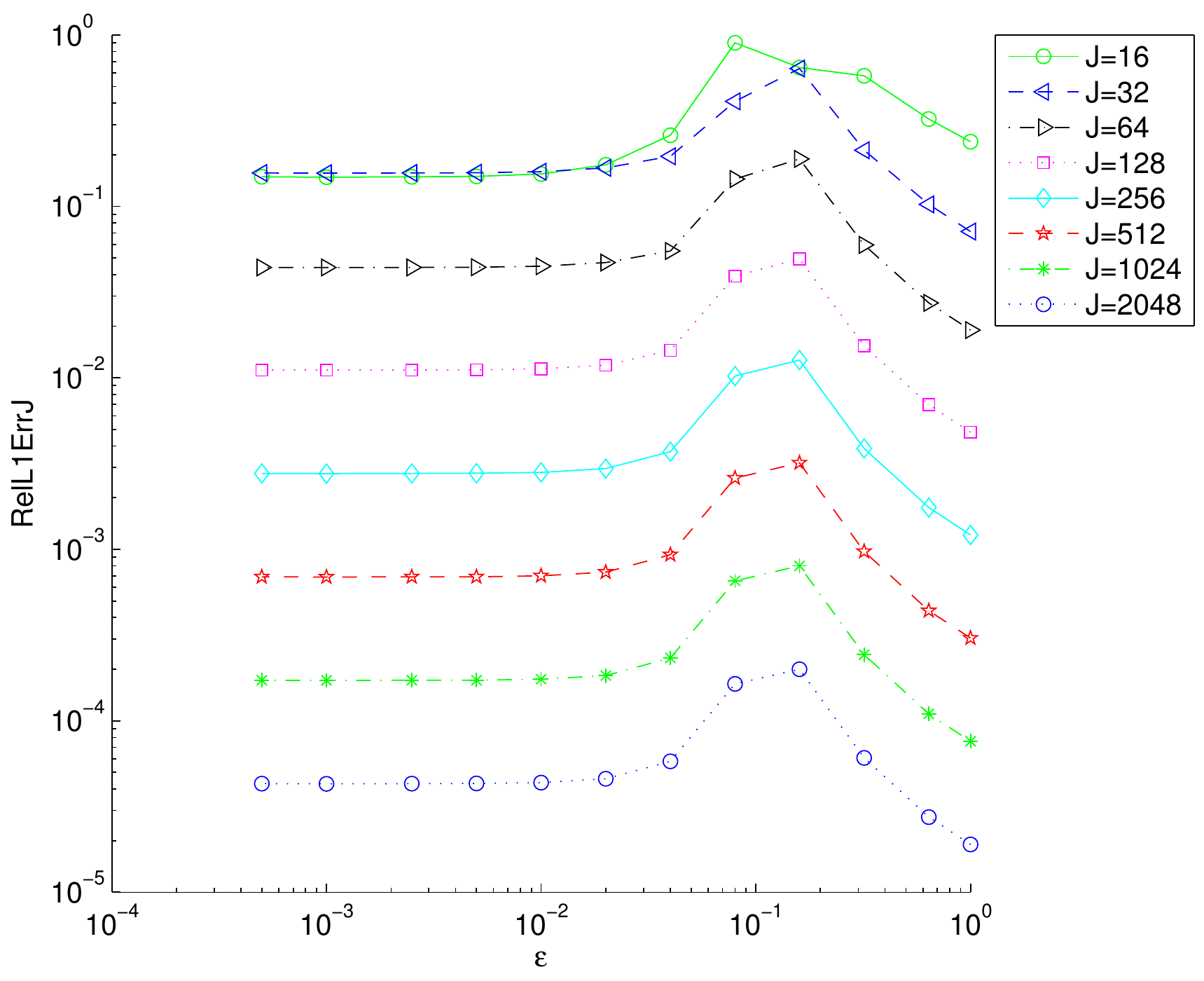}}
  
  \caption{$\mathrm{err}_{\mathbf{j}^\eps}(t=0.05)$ for AP scheme}
\label{fig:err2D_j_bef_sing_ap}
\end{figure}

\begin{figure}[!htbp]
  \centering
  \subfigure[Error w.r.t $J=2^M$]{\includegraphics[width=.46\textwidth]{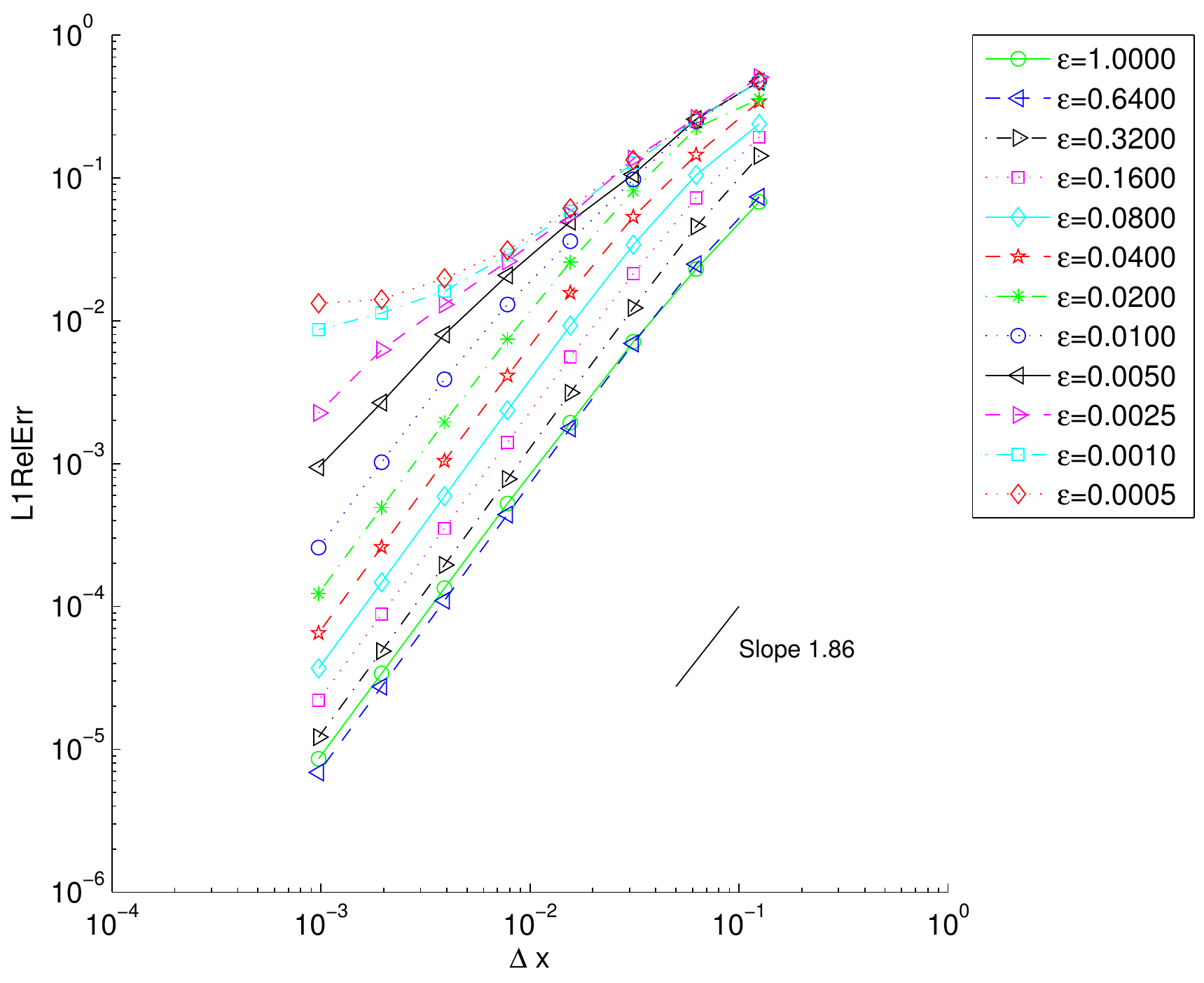}}\qquad
  \subfigure[Error w.r.t $\eps$]{\includegraphics[width=.46\textwidth]{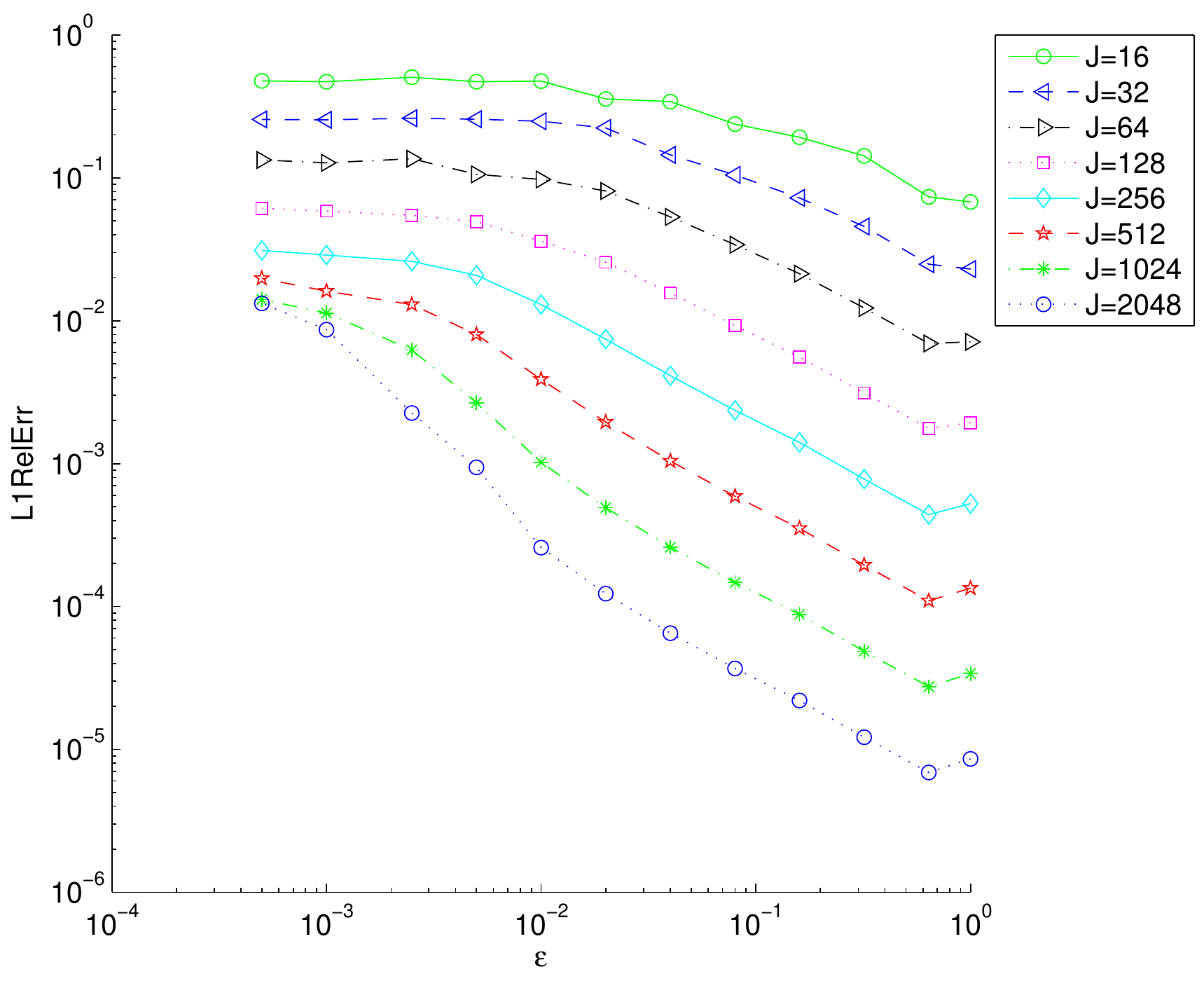}}
  \caption{$\mathrm{err}_{\rho^\eps}(t=0.13)$ for AP scheme}
\label{fig:err2D_rho_aft_sing_ap}  
\end{figure}

\begin{figure}[!htbp]
  \centering
  \subfigure[Error w.r.t $J=2^M$]{\includegraphics[width=.46\textwidth]{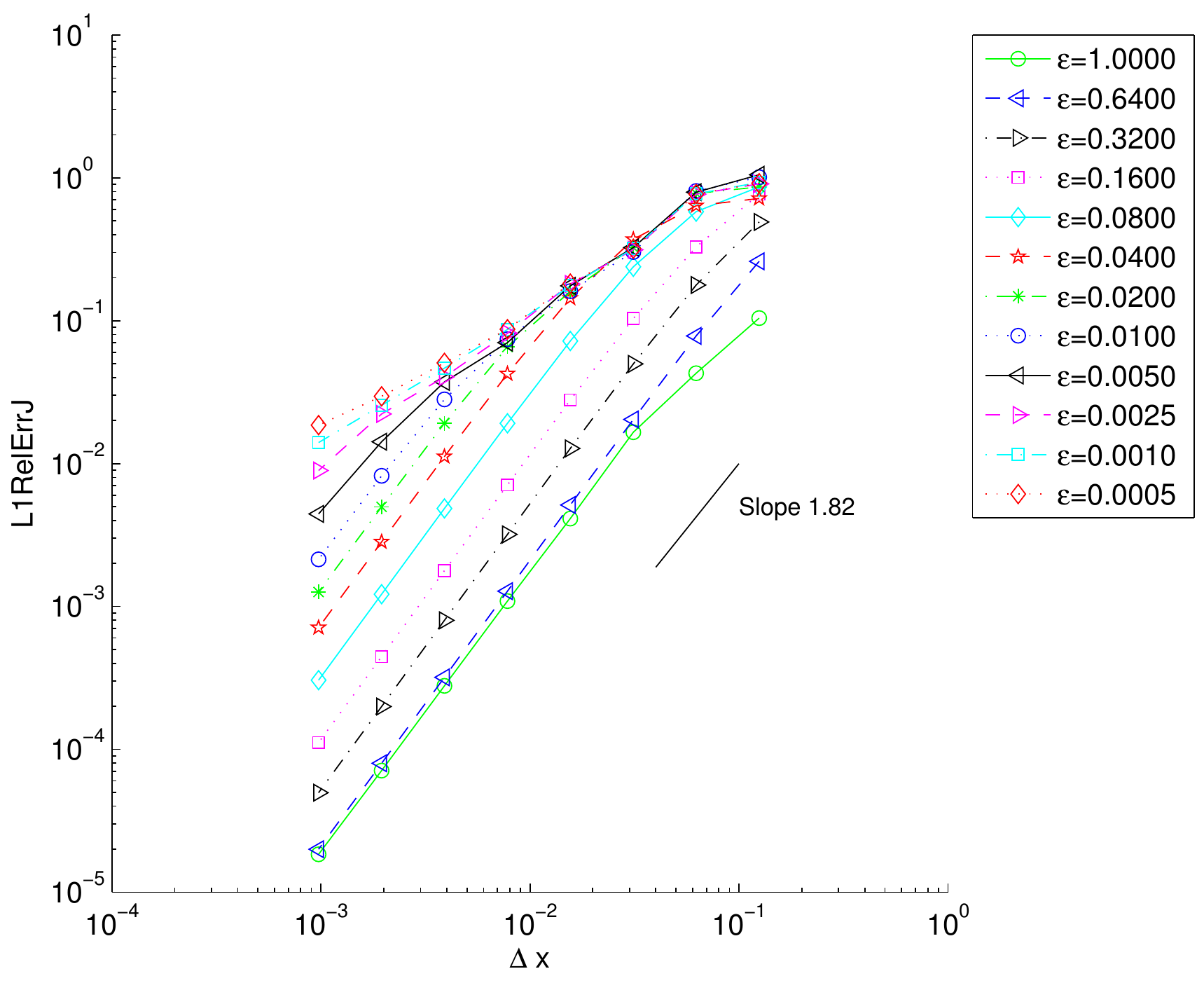}}\qquad
  \subfigure[Error w.r.t $\eps$]{\includegraphics[width=.46\textwidth]{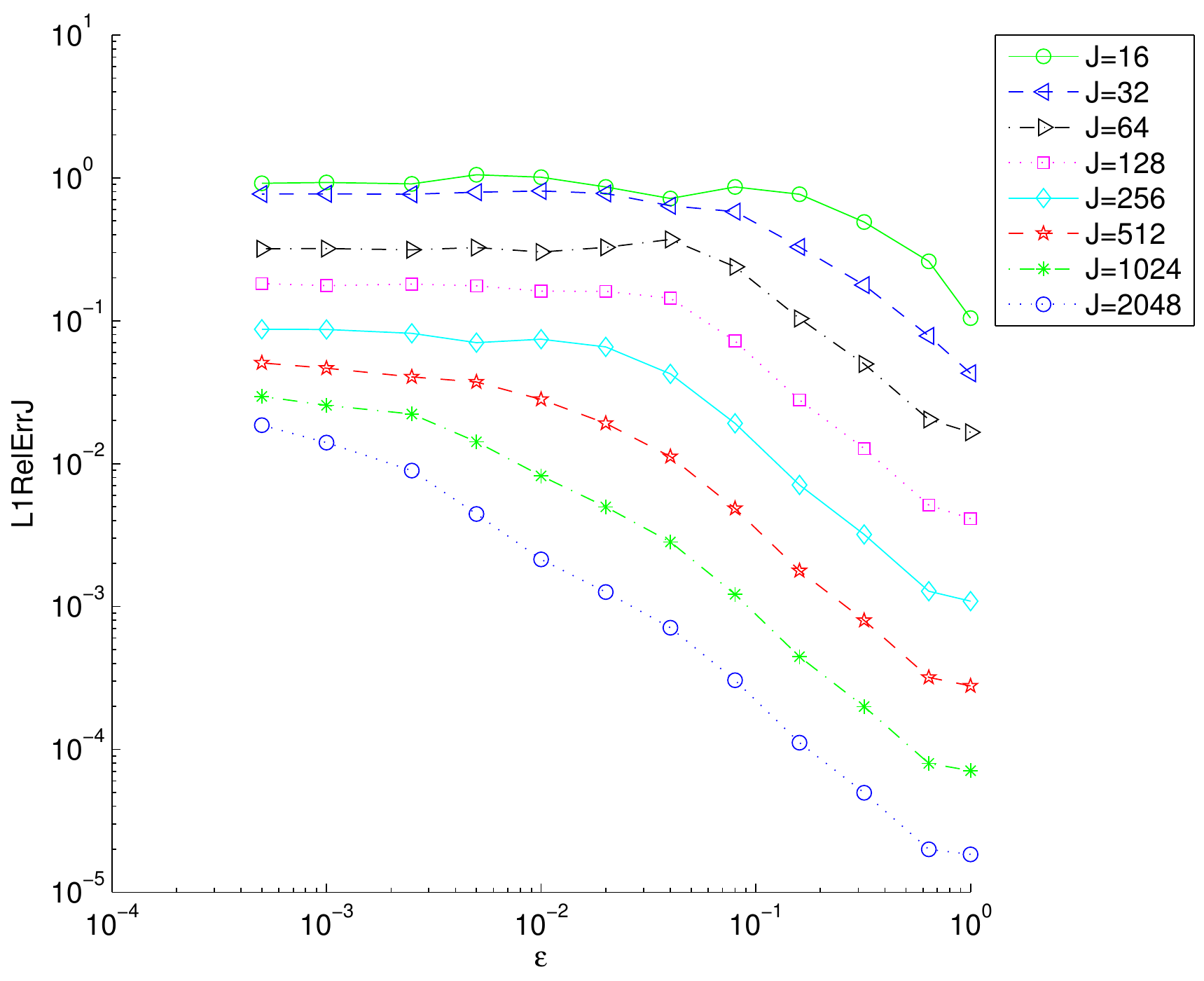}}
  \caption{$\mathrm{err}_{\mathbf{j}^\eps}(t=0.13)$ for AP scheme}
  \label{fig:err2D_j_aft_sing_ap}
\end{figure}

The error analysis results are equivalent to the ones obtained for one
dimensional simulation. We recover an independence of the error with
respect to $\eps$ before the formation of singularities (see Figs
\ref{fig:err2D_rho_bef_sing_ap} and \ref{fig:err2D_j_bef_sing_ap}). We
always need to take finer mesh past the shocks formation (see Figs
\ref{fig:err2D_rho_aft_sing_ap}  and \ref{fig:err2D_j_aft_sing_ap}).

\begin{figure}[!htbp]
  \centering
  \includegraphics[width=.60\textwidth]{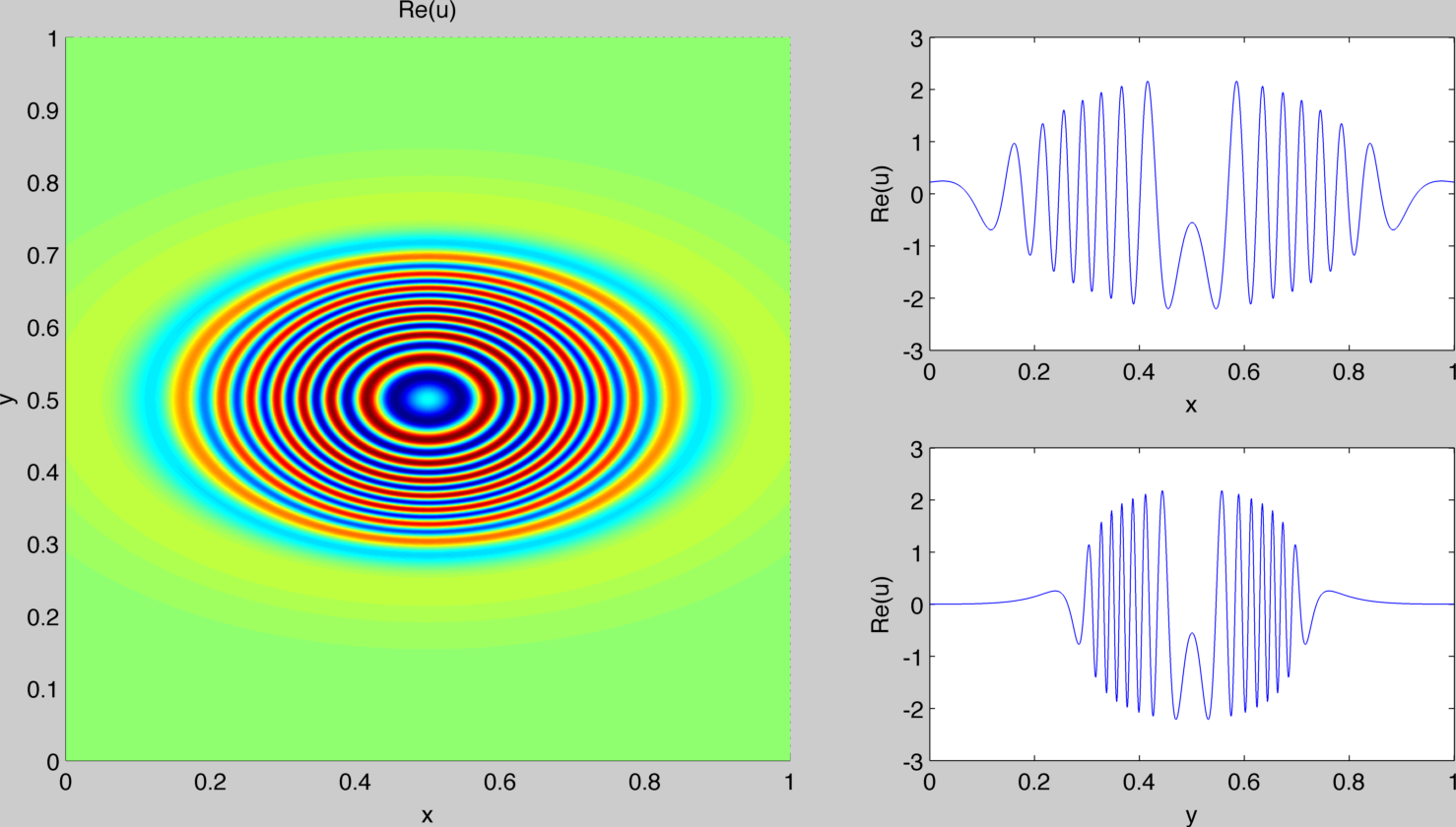}
  \caption{$\text{Re}(u^\eps)(t=0.035)$}
  \label{fig:Re_u_strang_t_0_035}
\end{figure}

\begin{figure}[!htbp]
  \centering
  \includegraphics[width=.60\textwidth]{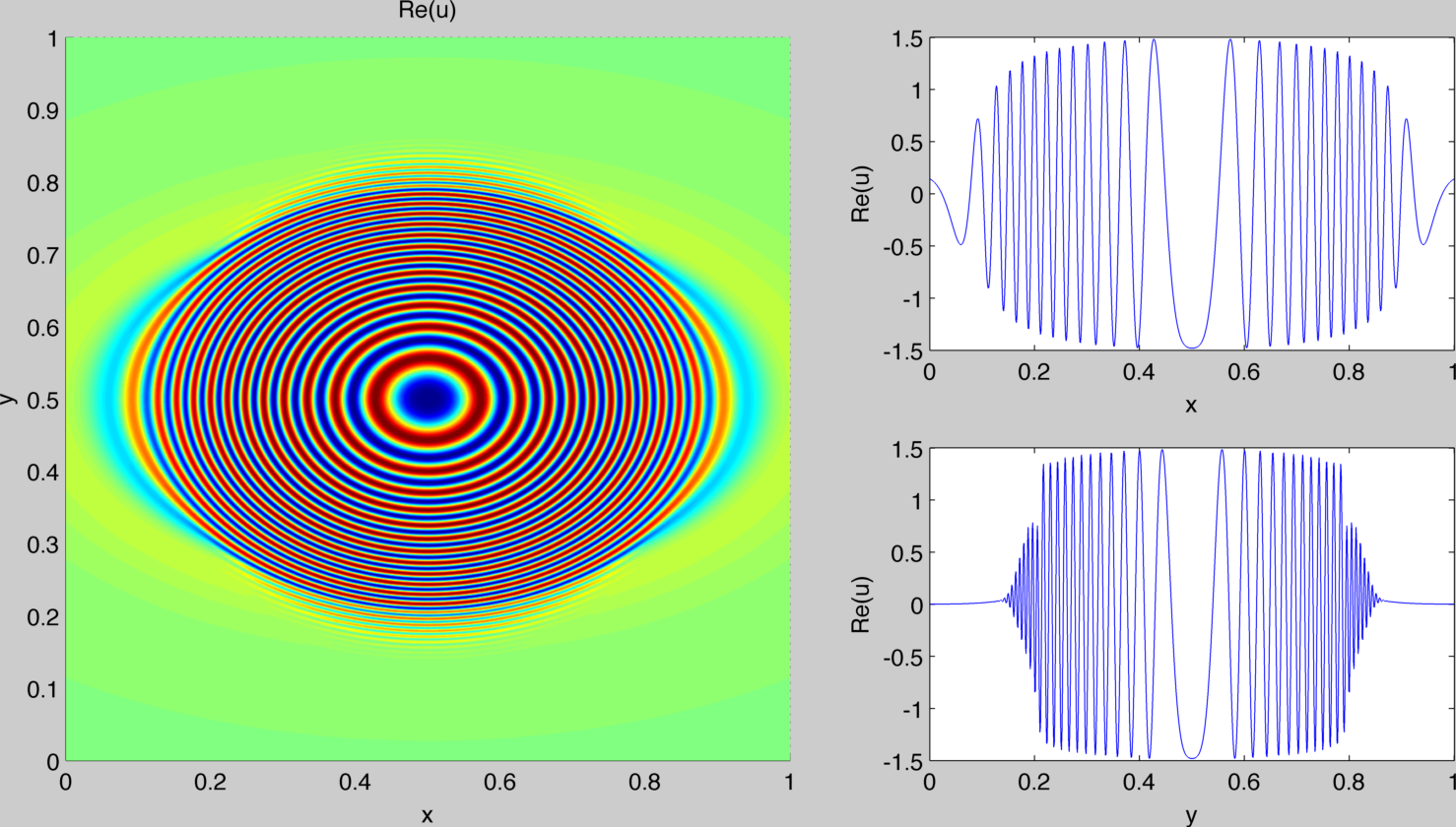}
  \caption{$\text{Re}(u^\eps)(t=0.08)$}
  \label{fig:Re_u_strang_t_0_08}
\end{figure}

The last computation concerns a non-symmetric solution. The aim of
our simulations is to show that we recover good qualitative
properties, even for non radially symmetric data. We consider here as
initial datum a Maxwell distribution with two temperatures
$\theta_1=0.05$ and $\theta_2=0.015$ for the initial amplitude given by
  $$ a_0=\frac{1}{\sqrt{4}} \frac{1}{2 \pi \sqrt{\theta_1}
    \sqrt{\theta_2}} \exp{\left ({-\frac{(x-0.5)^2}{2
          \theta_1}-\frac{(y-0.5)^2}{2 \theta_2}}\right )}. $$
The initial phase is reduced to $\phi_0=0$ (note that from
\eqref{eq:bcm0}, $\d_t \phi_{\mid t=0}=-|a_0|^2\not =0$, so $\phi$
becomes instantaneously non-trivial). The scaled Planck factor
$\eps$ is equal to $\eps=0.005$ and the simulations are performed with
$J=2048$ intervals in each direction. To distinguish the evolution for
each direction $(x,y)$, we make contour plots of wave function and
physical observables but also present slice for plane $x=0.5$ and
$y=0.5$ respectively on bottom right and upper right axes. The
solution has many oscillations in both directions before formation of
singularities (see the real part of the wave function $u^\eps$ on
Fig. \ref{fig:Re_u_strang_t_0_035} at time $t=0.035$). The oscillatory
nature is enhanced at time $t=0.08$ mainly in $y-$direction (see Fig. \ref{fig:Re_u_strang_t_0_08}). 

\begin{figure}[!htbp]
  \centering
  \includegraphics[width=.65\textwidth]{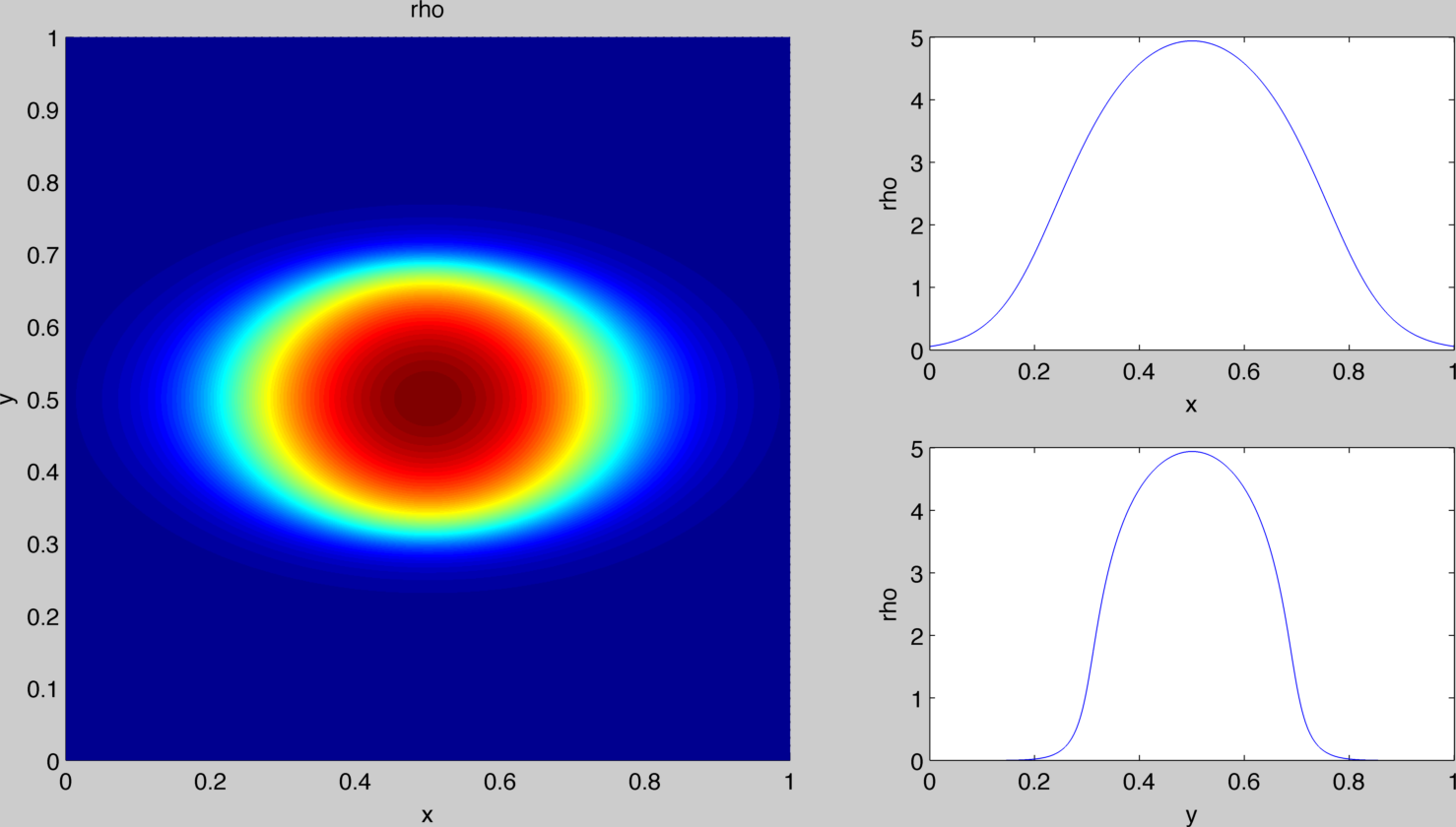}
  \caption{Strang $\rho^\eps(t=0.035)$}
  \label{fig:rho_strang_t_0_035}
\end{figure}

\begin{figure}[!htbp]
  \centering
  \includegraphics[width=.65\textwidth]{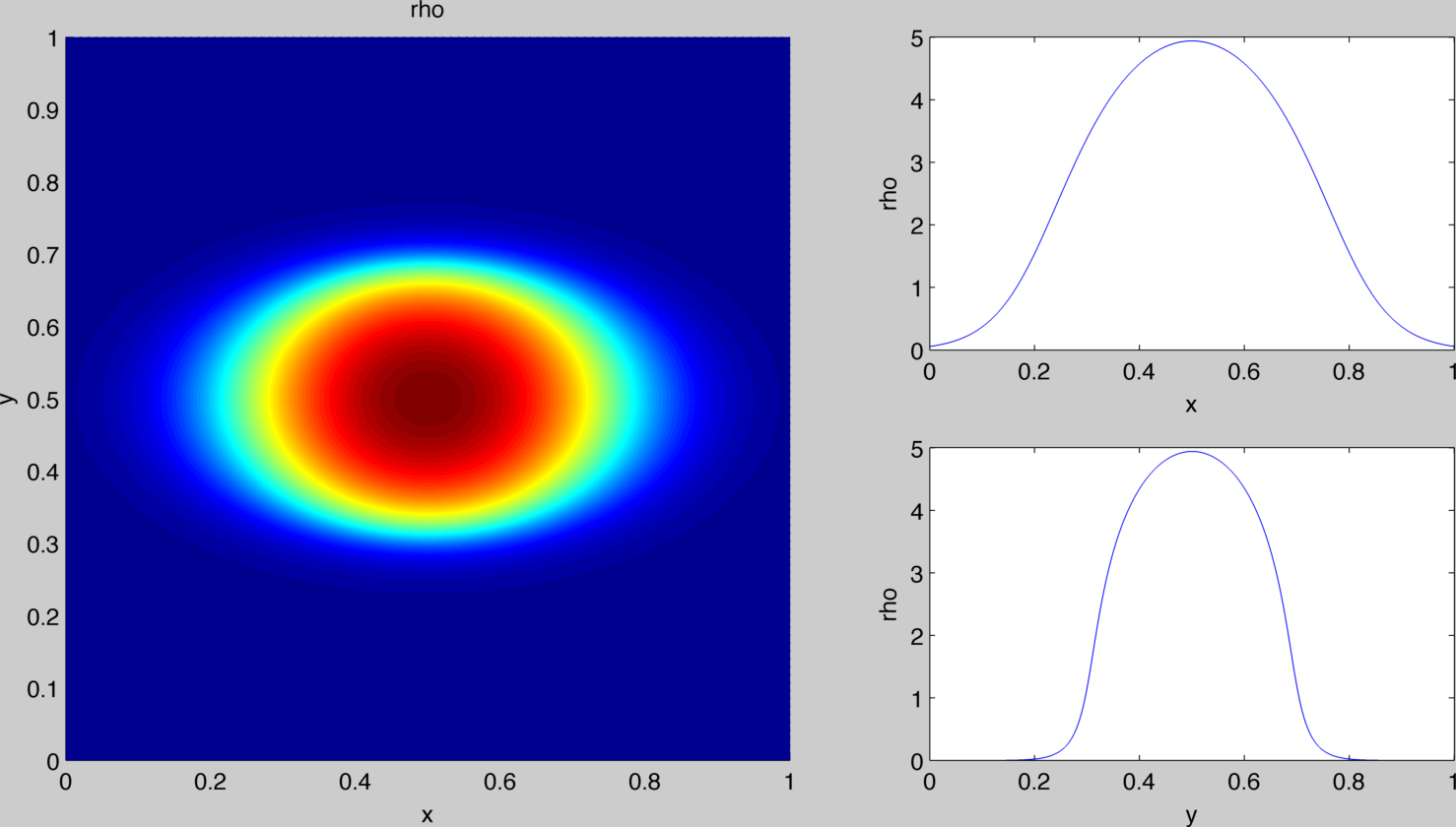}
  \caption{AP $\rho^\eps(t=0.035)$}
  \label{fig:rho_apnls_t_0_035}
\end{figure}

\begin{figure}[!htbp]
  \centering
  \includegraphics[width=.65\textwidth]{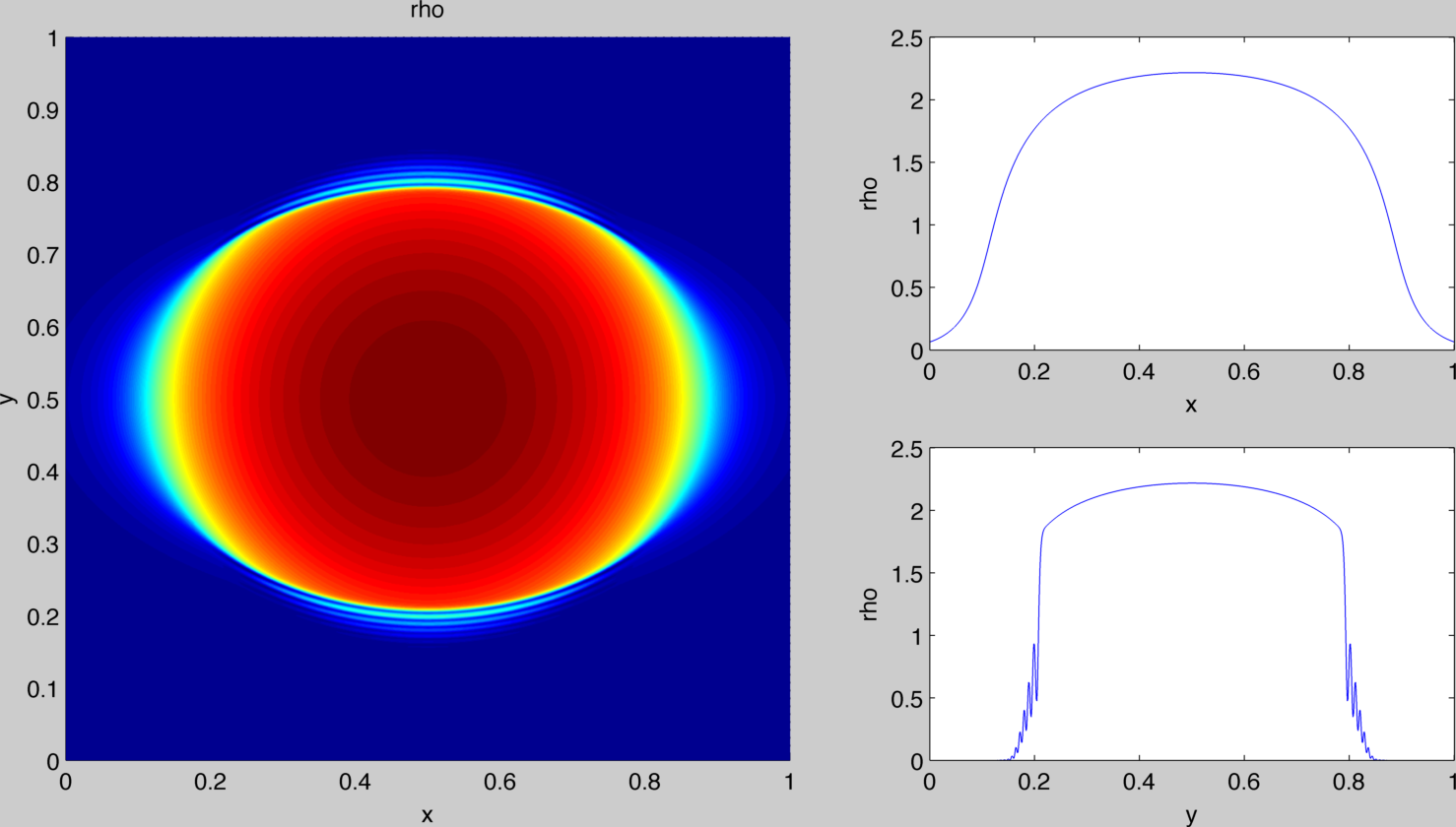}
  \caption{Strang $\rho^\eps(t=0.08)$}
  \label{fig:rho_strang_t_0_08}
\end{figure}

\begin{figure}[!htbp]
  \centering
  \includegraphics[width=.65\textwidth]{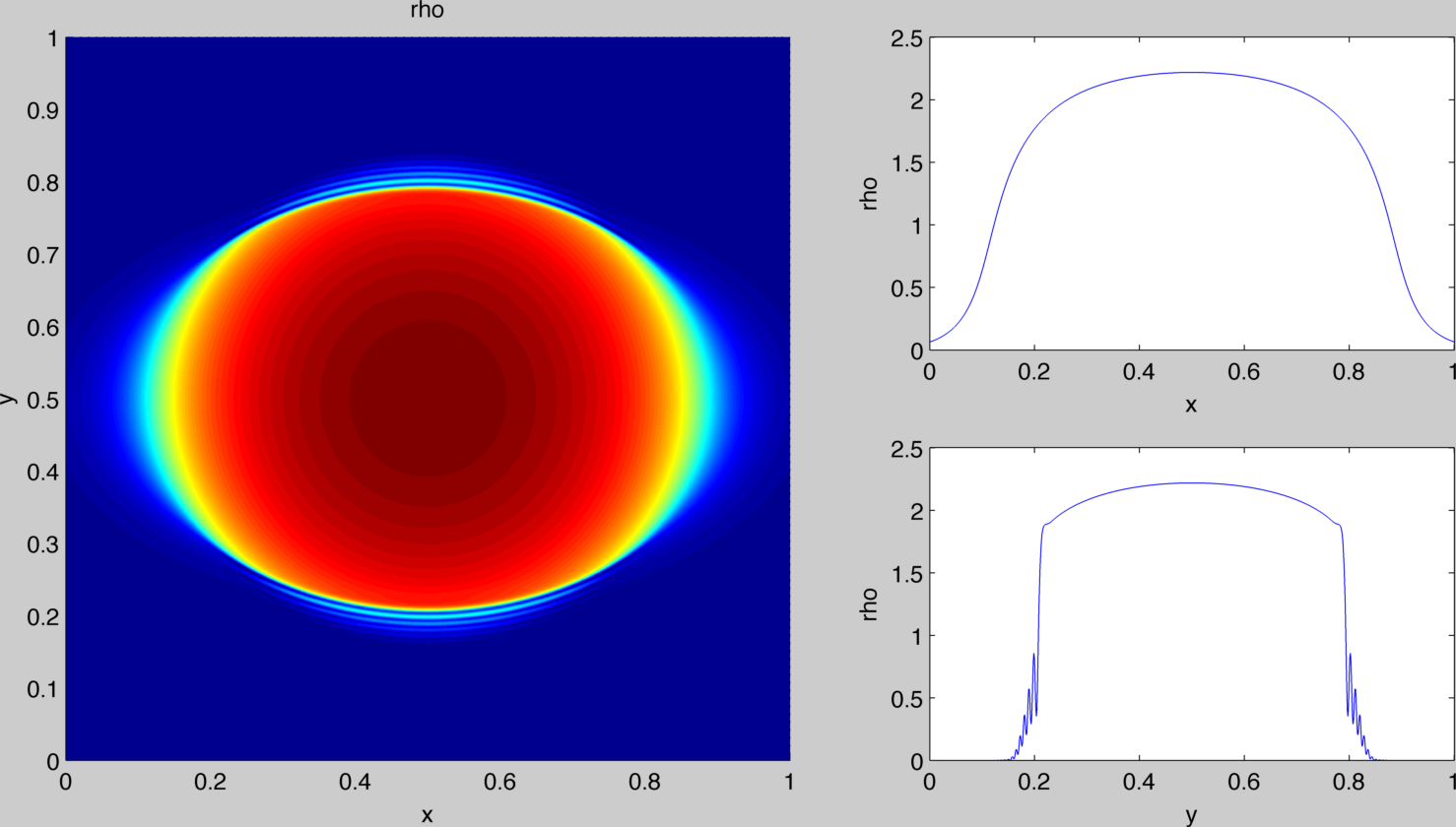}
  \caption{AP $\rho^\eps(t=0.08)$}
  \label{fig:rho_apnls_t_0_08}
\end{figure}

\begin{figure}[!htbp]
  \centering
  \includegraphics[width=.65\textwidth]{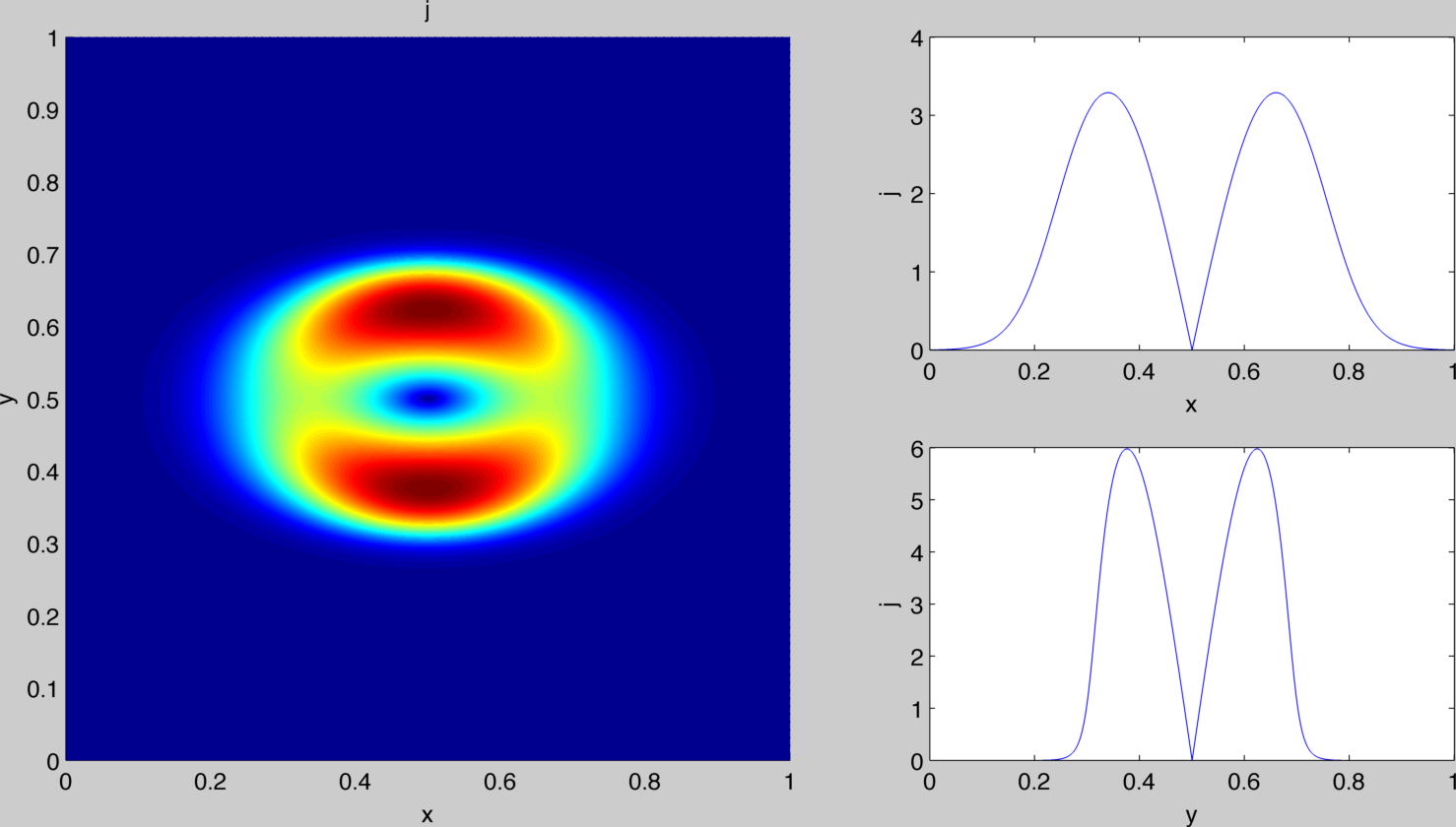}
  \caption{Strang ${\mathbf j}^\eps(t=0.035)$}
  \label{fig:J_strang_t_0_035}
\end{figure}

\begin{figure}[!htbp]
  \centering
  \includegraphics[width=.65\textwidth]{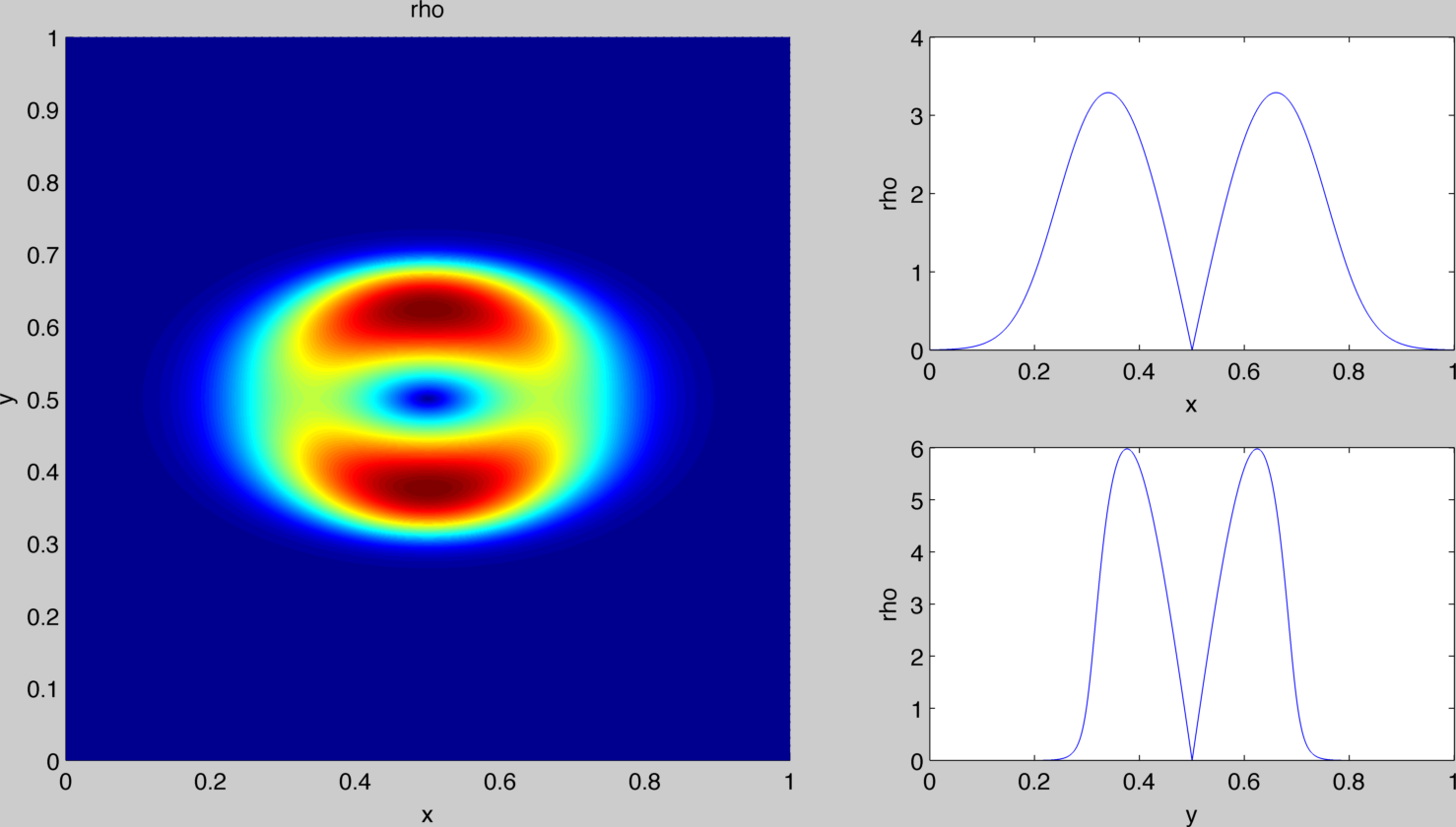}
  \caption{AP ${\mathbf j}^\eps(t=0.035)$}
  \label{fig:J_apnls_t_0_035}
\end{figure}

\begin{figure}[!htbp]
  \centering
  \includegraphics[width=.65\textwidth]{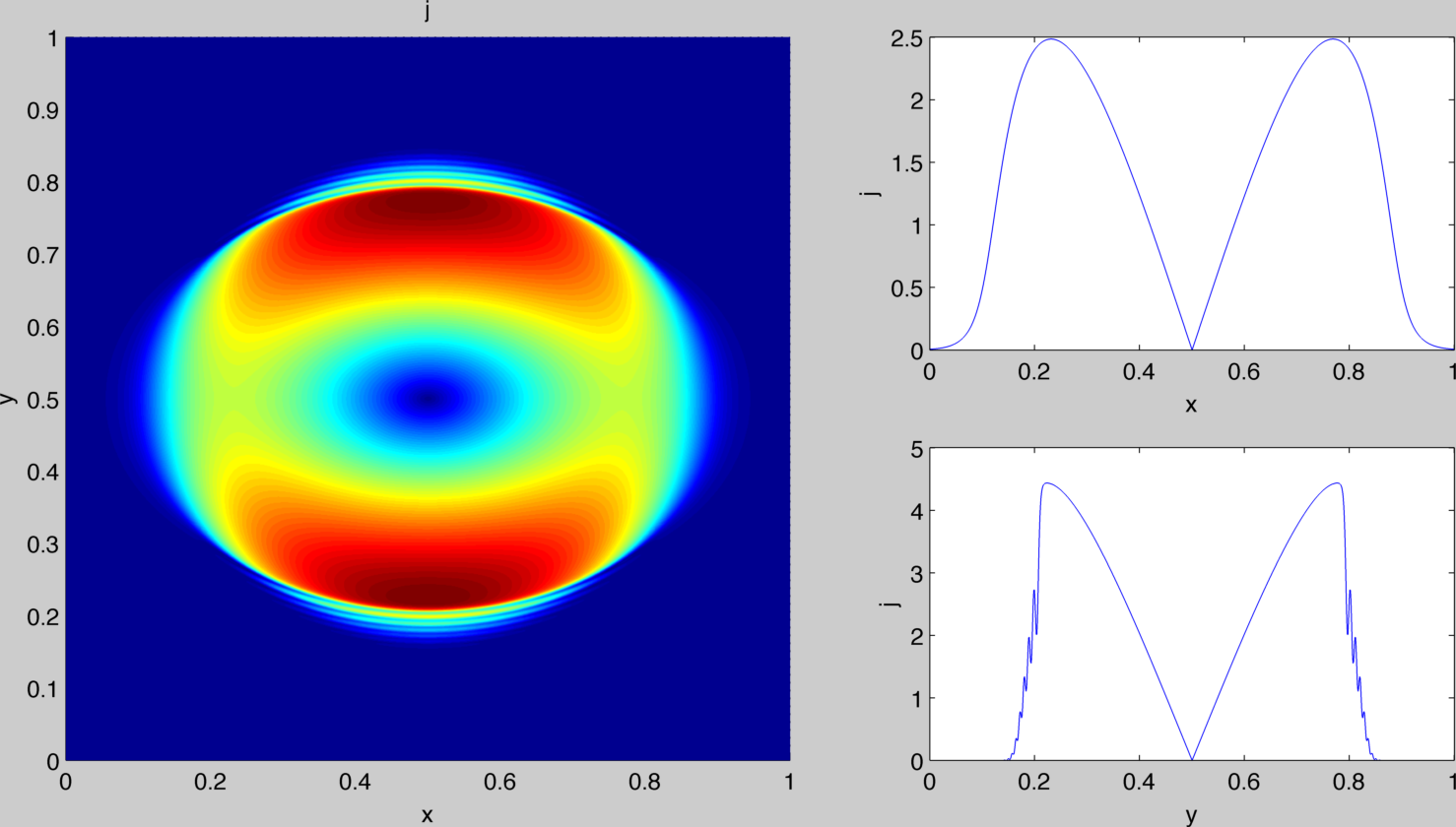}
  \caption{Strang ${\mathbf j}^\eps(t=0.08)$}
  \label{fig:J_strang_t_0_08}
\end{figure}

\begin{figure}[!htbp]
  \centering
  \includegraphics[width=.65\textwidth]{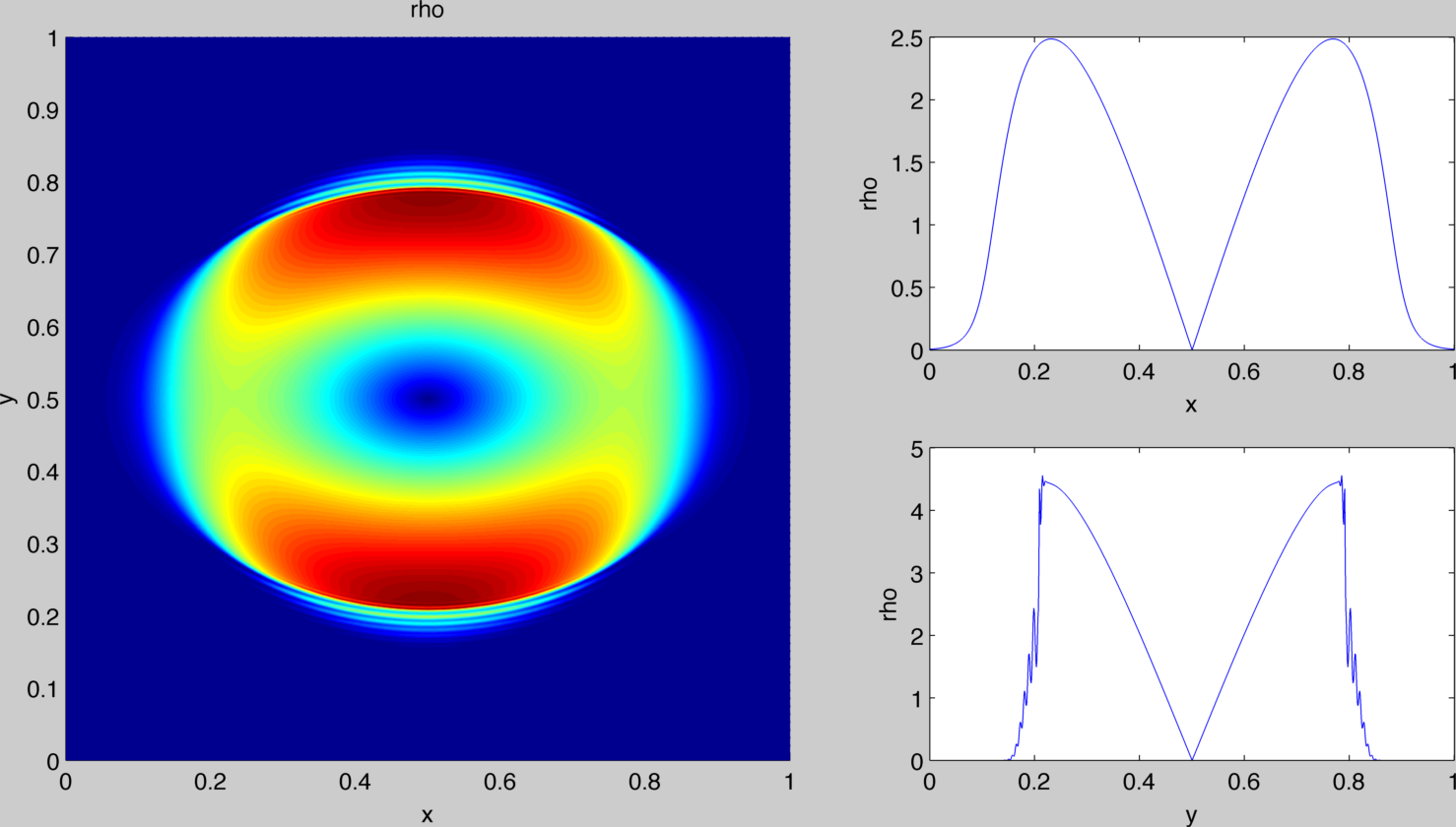}
  \caption{AP ${\mathbf j}^\eps(t=0.08)$}
  \label{fig:J_apnls_t_0_08}
\end{figure}

We recover the good behavior of physical observables before the
formation time of singularities. Actually, we compare the particle
densities computed by Strang splitting scheme and our AP scheme at
time $t=0.035$ on
Figures~\ref{fig:rho_strang_t_0_035} and
\ref{fig:rho_apnls_t_0_035}. The similar representation for
$\mathbf{j}^\eps$ is available on Figures~\ref{fig:J_strang_t_0_035}
and \ref{fig:J_apnls_t_0_035}. Past the formation time of
singularities, the contour plots could be thought as equivalent both
for particle and current densities. The only noticeable differences
can be seen on $x=0.5$ and $y=0.5$ plane slices (see
Fig.~\ref{fig:rho_strang_t_0_08}, \ref{fig:rho_apnls_t_0_08},
\ref{fig:J_strang_t_0_08} and \ref{fig:J_apnls_t_0_08}).

\clearpage

\section{Extension to other frameworks}
\label{sec:extension}

\subsection{The linear case: global solution of the viscous eikonal equation}\label{sec:linear}

%{\tt Commentaire de Florian. \`a faire: d\'ecider ce qu'on fait de
%cette section qui me semble quand m\^eme int\'eressante ici. Est-ce
%qu'on montre un th\'eor\`eme complet ? Si non, on peut motiver cette
%section encore par le num\'erique, car il est beaucoup plus
%int\'eressant de simuler HJ visqueuse. Je n'ai pas touch\'e \`a cette
%preuve fournie par R\'emi.} 

The advantage of introducing an artificial diffusion is more striking
in the linear case. Consider 
\begin{equation}
  \label{eq:linear}
  i\eps\d_t u^\eps +\frac{\eps^2}{2}\Delta u^\eps = V_{\rm    ext}
  u^\eps;\quad u^\eps_{\mid t=0} =a_0 e^{i\phi_0/\eps},
\end{equation}
with $V_{\rm    ext}=V_{\rm    ext}(t,x)$ real-valued. As noticed in
\cite{DGM}, using the Madelung transform in this context is
interesting only if vacuum can be avoided. Having the idea of Grenier
\cite{Grenier98} in mind, the counterpart of \eqref{eq:systemmanuel}
reads
\begin{equation*}
   \left\{
    \begin{aligned}
     &\d_t a^\eps +\nabla \phi^\eps \cdot \nabla a^\eps
      +\frac{1}{2}a^\eps \Delta \phi^\eps =i\frac{\eps}{2}\Delta
      a^\eps,\qquad  && a^\eps_{\mid
      t=0}=a_0,\\
      & \d_t \phi^\eps +\frac{1}{2}|\nabla \phi^\eps|^2+
       V_{\rm    ext} = 0,\qquad &&
      \phi^\eps_{\mid t=0}=\phi_0.  \end{aligned}
\right.
\end{equation*}
We note that the two equations are decoupled: the second equation is
of Hamilton-Jacobi type (known as the \emph{eikonal equation} in this
context), and the solution must be expected to develop 
singularities in finite time (see e.g. \cite{bookRemi}), so this
approach is doomed to fail for large time. On the other hand,
introducing the artificial diffusion as in \eqref{eq:bcm0}, we obtain
\begin{equation}\label{eq:bcmlin}
   \left\{
    \begin{aligned}
     &\d_t a^\eps +\nabla \phi^\eps \cdot \nabla a^\eps
      +\frac{1}{2}a^\eps \Delta \phi^\eps =i\frac{\eps}{2}\Delta
      a^\eps-i\eps a^\eps \Delta \phi^\eps,\qquad  && a^\eps_{\mid
      t=0}=a_0,\\
      & \d_t \phi^\eps +\frac{1}{2}|\nabla \phi^\eps|^2+
       V_{\rm    ext} = \eps^2\Delta \phi^\eps,\qquad &&
      \phi^\eps_{\mid t=0}=\phi_0.  \end{aligned}
\right.
\end{equation}
The system is still decoupled, but the good news is that the diffusion
introduced into the Hamilton-Jacobi equation smoothes the solution
out. Therefore, this system can be considered as a good candidate to
obtain an AP scheme for the linear Schr\"odinger equation, regardless
of the presence of vacuum. 
\smallbreak

Instead of giving full details of the analogue of Theorems~\ref{theo:local},
\ref{theo:cveuler} and \ref{theo:global}, as well as of their proofs,
we shall simply give a functional framework where the viscous eikonal
equation can by solved globally in time. Since we consider the case
$\eps>0$ only, we drop out the $\eps^2$ term in the viscous eikonal
equation, and consider
\begin{equation}\label{eq:eikonalv}
  \d_t \eik +\frac{1}{2}|\nabla \eik|^2 +V_{\rm    ext}=\Delta \eik;\quad
  \eik(0,x)=\phi_0(x). 
\end{equation}
For $k\ge 1$, let
\begin{equation*}
  {\tt SL}_k = \{f\in C^k(\R^d;\R),\quad \d^\alpha f\in L^\infty(\R^d),
  \ \forall 1\le |\alpha|\le k\}.
\end{equation*}
For $k\ge 2$, let
\begin{equation*}
  {\tt SQ}_k = \{f\in C^k(\R^d;\R),\quad \d^\alpha f\in L^\infty(\R^d),
  \ \forall 2\le |\alpha|\le k\}.
\end{equation*}
The key result is the following:
\begin{lemma}\label{lem:eikonalv}
  Let $k\ge 2$ and  $\phi_0\in {\tt SL_k}$, $V_{\rm    ext}\in L^\infty_{\rm
    loc}(\R_+;{\tt SL}_k)$. Then
  \eqref{eq:eikonalv} has a unique solution $\eik\in C(\R_+;{\tt
    SL}_{k-1})$. 
\end{lemma}
\begin{remark}
No such global result must be expected in the larger class ${\tt
    SQ}_k$. Indeed, suppose $d=1$, $V_{\rm ext}=0$ and $\phi_0(x) = -x^2/2$. If
  \eqref{eq:eikonalv} had a solution $\eik \in C(\R_+;{\tt
    QL}_{2})$, then $w=\d_x^2 \eik\in C(\R_+;L^\infty)$ would solve
  \begin{equation*}
    \d_t w + v\d_x w + w^2 = \d_x^2 w;\quad w(0,x)=-1.
  \end{equation*}
The solution of this equation does not depend on $x$. It is given by
$ w(t,x) = 1/(t-1)$,  hence a contradiction. 
  However, 
  in the periodic case $x\in \T^d$, this technical issue disappears
  (see Remark~\ref{rem:wiener} below). 
\end{remark}
\begin{proof}
  The first step consists in constructing the gradient of
  $\eik$. Differentiating \eqref{eq:eikonalv} in space, we have to
  solve
  \begin{equation}
    \label{eq:veik}
    \d_t v + v\cdot\nabla v + \nabla V_{\rm    ext} = \Delta v;\quad
    v(0,x)=\nabla \phi_0(x). 
  \end{equation}
This equation is solved locally in time by a fixed point argument by
using Duhamel's formula
\begin{equation*}
  v(t) = e^{t\Delta} \nabla \phi_0 -\int_0^t e^{(t-s)\Delta}(v\cdot
  \nabla v)(s)ds - \int_0^t e^{(t-s)\Delta}(\nabla V_{\rm
    ext})(s)ds. 
\end{equation*}
The right hand side is a contraction on 
\begin{equation*}
\left\{ w\in  C\([0,T];W^{k-1,\infty}\) ,\quad
  \|w\|_{L^\infty([0,T];W^{k-1,\infty})}\le 2 \|\nabla
  \phi_0\|_{W^{k-1,\infty}}\right\}, 
\end{equation*}
provided that $T>0$ is sufficiently small. This follows from the
properties of the heat kernel: for $t>0$, 
\begin{equation*}
  \|e^{t\Delta}f\|_{L^\infty}\le \|f\|_{L^\infty};\quad
  \|e^{t\Delta}\nabla f\|_{L^\infty}\le \frac{C}{\sqrt t}\|f\|_{L^\infty} .
\end{equation*}
The solution to \eqref{eq:veik} is then global, thanks to the maximum
principle (see e.g. \cite[Proposition~52.8]{QuSo07}, which implies
that $v_-\le v\le v_+$, where 
\begin{equation*}
  \d_t v_\pm \mp \|\nabla V_{\rm ext}\|_{L^\infty_x} =0 ;\quad
  v_\pm(0,x)=\pm \|\nabla\phi_0\|_{L^\infty}. 
\end{equation*}
Once such a solution $v\in C(\R_+;W^{k-1,\infty})$ is constructed, set
\begin{equation*}
  \eik (t)=\phi_0 -\int_0^t \(\frac{1}{2}|v(s)|^2 + V_{\rm
    ext}(s)-\DIV v(s)\)ds. 
\end{equation*}
We check that $\d_t (\nabla \eik-v)=\nabla \d_t\eik -\d_t
  v=0$. Therefore, $v=\nabla \eik$, and $\eik$ solves
  \eqref{eq:eikonalv}. 
Rewriting \eqref{eq:eikonalv} as 
\begin{equation*}
  \d_t \eik +\frac{1}{2}|v|^2 +V_{\rm    ext}=\Delta \eik;\quad
  \eik(0,x)=\phi_0(x),  
\end{equation*}
Duhamel's formula reads
\begin{equation*}
  \eik(t)= e^{t\Delta }\phi_0 -\int_0^t e^{(t-s)\Delta}(|v(s)|^2)ds -
  \int_0^t e^{(t-s)\Delta} V_{\rm 
    ext}(s)ds, 
\end{equation*}
from which we conclude that $\eik \in C(\R_+;{\tt
    SL}_{k-1})$.
\end{proof}
\begin{remark}\label{rem:wiener}
  In the periodic setting $x\in \T^d= (\R/(2\pi \Z))^d$, it is natural to work in the
  Wiener algebra
  \begin{equation*}
    W = \left\{f :\T^d\to \C,\quad 
  f(x) = \sum_{j\in \Z^d}b_j e^{ij\cdot x}
\text{ with }(b_j)_{j\in \Z^d}\in \ell^1(\Z^d)\right\},
  \end{equation*}
and for $k\ge 0$,
\begin{equation*}
  W^k = \{f :\T^d\to \C,\quad \d^\alpha f\in W,\ \forall |\alpha|\le
  k\}. 
\end{equation*}
Then Lemma~\ref{lem:eikonalv} is easily adapted by replacing ${\tt
  SL}_{k}$ with $W^k$. 
\end{remark}

\subsection{Other nonlinearities}
\label{sec:otherNL}
One might believe that Theorem~\ref{theo:global} is bound to the cubic
one-dimensional Schr\"odinger equation, which is completely
integrable, and that the two aspects are related. We show that this is
not the case, since Theorem~\ref{theo:global} remains valid for other
nonlinearities (with $d=1$). Consider now
\begin{equation}
  \label{eq:nlsgen}
  i\eps\d_t u^\eps +\frac{\eps^2}{2}\Delta u^\eps = f\( |u^\eps|^2\)
  u^\eps,\quad (t,x)\in \R_+\times \R^d.
\end{equation}
We suppose that there exists $\delta>0$ such that $f'(y)\ge \delta$,
for all $y\ge 0$. 

Following \cite{Grenier98}, the only modification we have to make to
recover the results in Section~\ref{sec:local} 
consists in replacing the symmetrizer \eqref{eq:symetriseur} with
\begin{equation*}
  S=\begin{pmatrix}
     I_2 & 0\\
     0& \frac{1}{4f'(a_1^2+a_2^2)}I_d
    \end{pmatrix}.
\end{equation*}
\mbox{}From our assumption on $f$, $S$ and its inverse $S^{-1}$ are uniformly
bounded, provided that $a$ is bounded. As a matter of fact, the exact
assumption in \cite{Grenier98} is $f'>0$, and all the results in
Section~\ref{sec:local} remain valid under this assumption. 

The more precise assumption $f'(y)\ge \delta$ becomes useful to prove
that in dimension $d=1$, the solution is global. The main point to
notice is that the conservation of mass and momentum are the same as
for \eqref{eq:nls}, like in
Proposition~\ref{prop:nls}. On the other hand, the conservation of
energy becomes
\begin{equation*}
  \frac{d}{dt}\(\|\eps\nabla
u^\eps(t)\|_{L^2}^2+\int_{\R^d}F\(|u^\eps(t,x)|^2\)dx\)=0,
\end{equation*}
where $F$ is an anti-derivative of $f$,
\begin{equation*}
  F(y)= \int_0^y f(r)dr. 
\end{equation*}
By assumption,
\begin{equation*}
  F(y) \ge \frac{\delta}{2}y^2 +f(0)y, 
\end{equation*}
so the potential energy controls the $L^4$-norm:
\begin{align*}
  \|u^\eps(t)\|_{L^4}^4& \le
  \frac{2}{\delta}\(\int_{\R^d}F\(|u^\eps(t,x)|^2\)dx -
  f(0)\|u^\eps(t)\|_{L^2}^2\)\\
& \le
  \frac{2}{\delta}\(\int_{\R^d}F\(|u^\eps(t,x)|^2\)dx - 
  f(0)\|u^\eps(0)\|_{L^2}^2\) ,
\end{align*}
where we have used the conservation of mass in the last
inequality. This implies an estimate of the form
\begin{equation*}
  \|\eps\nabla
u^\eps(t)\|_{L^2}^2+\|u^\eps(t)\|_{L^4}^4\le C,
\end{equation*}
with $C$ independent of $t\ge 0$ and $\eps\in (0,1]$. Therefore, the
analysis presented in Section~\ref{sec:global} can be repeated line by
line. 
\begin{example}
The assumption $f'\ge \delta>0$ is satisfied in the following cases:
  \begin{itemize}
  \item Cubic-quintic nonlinearity: $f(y)=y+\lambda y^2$, $\lambda \ge
    0$.
\item Cubic plus saturated nonlinearity: $f(y)=\delta y
  +\eta\frac{y}{1+\lambda y}$, 
  $\delta,\eta,\lambda>0$.  
  \end{itemize}
\end{example}

\subsection*{Acknowledgements}
\ni
C. Besse was partially supported by the French ANR fundings under the
project MicroWave NT09\_460489. R.~Carles was supported by the French
ANR project
  R.A.S. (ANR-08-JCJC-0124-01). C. Besse and F. M\'ehats were partially supported
by the ANR-FWF Project Lodiquas  (ANR-11-IS01-0003).

 \end{document}